\documentclass[11pt,oneside]{amsart}
\usepackage[top=25mm, bottom=25mm, left=20mm, right=20mm]{geometry}
\usepackage{amsfonts,amssymb,amsmath,amsthm,amsxtra}
\usepackage[T1]{fontenc}
\usepackage[utf8]{inputenc}
\usepackage{bm}
\usepackage{mathtools}
\usepackage{dsfont}
\usepackage{mathrsfs}
\usepackage{enumitem}
\allowdisplaybreaks

\usepackage{graphics}
\usepackage{hyperref}

\numberwithin{equation}{section}
\newtheorem{theorem}[equation]{Theorem}

\newtheorem{lemma}[equation]{Lemma}
\newtheorem{proposition}[equation]{Proposition}
\theoremstyle{definition}
\newtheorem{remark}[equation]{Remark}
\newtheorem{definition}[equation]{Definition}

\newtheorem*{claim}{Claim}

\DeclareMathOperator\Imm{Im}
\DeclareMathOperator\Ree{Re}
\author{Erik Bahnson}
\address[Erik Bahnson]{
Department of Mathematics,
Rutgers University,
Piscataway, NJ 08854-8019, USA}
\email{eeb118@scarletmail.rutgers.edu}

\author{Leonidas Daskalakis} 
\address[Leonidas Daskalakis]{Institute of Mathematics,
Polish Academy of Sciences,
\'Sniadeckich 8,
00-656 Warszawa, Poland}
\email{ldaskalakis@impan.pl}

\author{Abbas Dohadwala}
\address[Abbas Dohadwala]{Department of Mathematics, Purdue University, 150 N. University Street, West
Lafayette, IN 47907-2067, USA}
\email{adohadwa@purdue.edu}

\author{Ish Shah} 
\address[Ish Shah]{Department of Mathematics,
Rutgers University,
Piscataway, NJ 08854-8019, USA}
\email{irs51@scarletmail.rutgers.edu}

\begin{document}

\title{Pointwise ergodic theorems along fractional powers of primes}
\maketitle
\begin{abstract}We establish pointwise convergence for nonconventional ergodic averages taken along $\lfloor p^c\rfloor$, where $p$ is a prime number and $c\in(1,4/3)$ on $L^r$, $r\in(1,\infty)$. In fact, we consider averages along more general sequences $\lfloor h(p)\rfloor$, where $h$ belongs in a wide class of functions, the so-called $c$-regularly varying functions. We also establish uniform multiparameter oscillation estimates for our ergodic averages and the corresponding multiparameter pointwise ergodic theorem in the spirit of Dunford and Zygmund. A key ingredient of our approach are certain exponential sum estimates, which we also use for establishing a Waring-type result. Assuming that the Riemann zeta function has any zero-free strip upgrades our exponential sum estimates to polynomially saving ones and this makes a conditional result regarding the behavior of our ergodic averages on $L^1$ to not seem entirely out of reach. 
\end{abstract}
\section{Introduction}
\subsection{One-parameter ergodic theorem along $c$-regular sequences evaluated at primes}One of the main results of the present work is the following pointwise ergodic theorem along fractional powers of primes.\begin{theorem}[Pointwise ergodic theorem along $(\lfloor p^c\rfloor)_{p\in\mathbb{P}}$]\label{basic}
Let $c\in(1,4/3)$, and assume that $(X,\mathcal{B},\mu)$ is a $\sigma$-finite measure space and $T\colon X\to X$ an invertible $\mu$-invariant transformation. Then for every $r\in(1,\infty)$ and every $f\in L_{\mu}^r(X)$ we have that
\[
\lim_{N\to\infty}\frac{1}{|\mathbb{P}\cap[1,N]|}\sum_{p\in\mathbb{P}\cap[1,N]}f(T^{\lfloor p^c\rfloor}x)\quad\text{exists for $\mu$-a.e. }x\in X\text{.}
\]
\end{theorem} 
This result will be obtained by a more general theorem allowing us to establish pointwise convergence along any of the following orbits 
\begin{equation}\label{examples1}
\lfloor p^c\log^{a}p\rfloor,\quad \lfloor p^ce^{a\log^b p}\rfloor,\quad \lfloor p^c\underbrace{\log\circ \cdots \circ \log p}_{k \text{ times}}\rfloor\text{,}
\end{equation}
where $a\in\mathbb{R}$, $b\in(0,1)$, $k\in\mathbb{N}$ are fixed and $p\in\mathbb{P}$. To formulate this theorem, we begin by introducing the so-called $c$-regularly varying functions.
\begin{definition}[$c$-regularly varying functions] Let $c\in(1,2)$ and $x_0\ge 1$. We define the class of $c$-regularly varying functions $\mathcal{R}_c$ as the set of all functions $h\in\mathcal{C}^3\big([x_0,\infty)\to[1,\infty)\big)$ such that the following conditions hold:
\begin{itemize}
\item[i)]$h'>0$, $h''> 0$\text{ and }
\item[ii)] the function $h$ is of the form 
\[
h(x)=Cx^c\exp\bigg(\int_{x_0}^x\frac{\vartheta(t)}{t}dt\bigg)\text{,}
\]
where $C$ is a positive constant and $\vartheta\in\mathcal{C}^2\big([x_0,\infty)\to\mathbb{R}\big)$ satisfies
\[
\vartheta(x)\to 0\,\,\text{, }x\vartheta'(x)\to 0\,\,\text{, }x^2\vartheta''(x)\to 0\,\,\text{ as }x\to \infty\text{.}
\]
\end{itemize}
\end{definition}
This family has been introduced in \cite{MMR} and, loosely speaking, one may think of the family $\mathcal{R}_c$ as pertubations of the function $x^c$. We now state one of the main theorems of the present work.
\begin{theorem}[Pointwise ergodic theorem along $(\lfloor h(p)\rfloor)_{p\in\mathbb{P}}$]\label{MT}
Let $c\in(1,4/3)$ and $h\in\mathcal{R}_c$, and assume $(X,\mathcal{B},\mu)$ is a $\sigma$-finite measure space and $T\colon X\to X$ an invertible $\mu$-invariant transformation. Then for every $r\in(1,\infty)$ and every $f\in L_{\mu}^r(X)$ we have that
\begin{equation}\label{ErgAv}
\lim_{N\to\infty}\frac{1}{|\mathbb{P}\cap[1,N]|}\sum_{p\in\mathbb{P}\cap[1,N]}f(T^{\lfloor h(p)\rfloor}x)\quad\text{exists for $\mu$-a.e. }x\in X\text{.}
\end{equation}
\end{theorem} 
It is not difficult to see that Theorem~$\ref{MT}$ implies Theorem~$\ref{basic}$ and it also establishes pointwise convergence along any of the orbits mentioned in $\ref{examples1}$ since $\mathcal{R}_c$ does include all the corresponding functions.

We derive Theorem~$\ref{MT}$ by establishing uniform 2-oscillation $L^r_{\mu}$-estimates, see $\ref{U2OE}$ in Theorem~$\ref{OM}$, which the reader may find of independent interest and are also exploited in the sequel. Before making precise formulations, let us introduce some notation. For any $Y\subseteq X\subseteq \mathbb{R}$, with $|X|>2$, let $\mathfrak{S}_J(X)=\{\{I_0<\dots<I_J\}\subseteq X\}$,
and for any family of complex-valued functions $(a_t(x):\,t\in X)$ we have
\[
O^2_{I,J}(a_t(x):\,t\in Y)=\Big(\sum_{j=0}^{J-1}\sup_{t\in [I_j,I_{j+1})\cap Y}|a_t(x)-a_{I_j}(x)|^2 \Big)^{1/2}\text{,}
\]
see subsections~2.6 and 2.7 from \cite{MOE} for the basic properties of oscillations. 
\begin{theorem}[Uniform 2-oscillation and vector-valued maximal estimates]\label{OM} Let $c\in(1,4/3)$ and $h\in\mathcal{R}_c$, and assume $(X,\mathcal{B},\mu)$ is a $\sigma$-finite measure space and $T\colon X\to X$ an invertible $\mu$-invariant transformation. Let
\[
A_Nf(x)=\frac{1}{|\mathbb{P}\cap[1,N]|}\sum_{p\in\mathbb{P}\cap[1,N]}f(T^{\lfloor h(p) \rfloor}x)\text{.}
\]
Then for all $r\in(1,\infty)$, there exists a positive constant $C_{r,h}$ such that
\begin{equation}\label{U2OE}
\sup_{J\in\mathbb{N}}\sup_{I\in\mathfrak{S}_J(\mathbb{N})}\|O^2_{I,J}(A_n f:\,n\in\mathbb{N})\|_{L_{\mu}^r(X)}\le C_{r,h}\|f\|_{L_{\mu}^r(X)}\text{  for any $f\in L_{\mu}^r(X)$}
\end{equation}
and such that  
\begin{equation}\label{VVMFE}
\Big\|\Big(\sum_{j\in\mathbb{Z}}\big(\sup_{N\in\mathbb{N}}A_{N}|f_j|\big)^2\Big)^{1/2}\Big\|_{L_{\mu}^r(X)}\le C_{r,h}\Big\|\big(\sum_{j\in\mathbb{Z}}|f_j|^2\big)^{1/2}\Big\|_{L_{\mu}^r(X)}\text{ for any $(f_j)_{j\in\mathbb{Z}}\in L_{\mu}^r\big(\mathbb{Z};\ell^2(\mathbb{Z})\big)$}\text{.}
\end{equation} 
\end{theorem}
\subsection{Multiparameter oscillation inequalities and multiparameter pointwise ergodic theorems}The second main result is establishing uniform multiparameter oscillation inequalities for our ergodic averages and the corresponding multiparameter extension of Theorem~$\ref{MT}$.
\begin{theorem}\label{Multi}
 Let $k\in\mathbb{N}$, and for every $i\in[k]$ let $c_i\in(1,4/3)$ and $h_i\in\mathcal{R}_{c_i}$. Assume $(X,\mathcal{B},\mu)$ is a $\sigma$-finite measure space and $\{S_i:\,i\in[k]\}$ is a family of invertible $\mu$-invariant transformations. For any $\mathbf{N}=(N_1,\dotsc,N_k)\in\mathbb{N}^k$ let
\[
A_{\mathbf{N}} f(x)=\frac{1}{\prod_{i=1}^k |\mathbb{P}\cap[1,N_i]|}\sum_{\mathbf{p} \in \prod_{i=1}^k \mathbb{P}\cap[1,N_i]}f (S_1^{\lfloor h_1(p_1)\rfloor}\dots S_k^{\lfloor h_k(p_k)\rfloor}x)\text{,}\quad\text{where $\mathbf{p}=(p_1,\dotsc,p_k)$.}
\] 
Then for every $r\in (1,\infty)$ and every $f\in L^r_{\mu}(X)$ we have that
\begin{equation}\label{multpointwise}
\lim_{\min\{N_1,\dotsc,N_k\}\to\infty}A_{\mathbf{N}}f(x)\text{ exists for $\mu$-a.e. $x\in X$.} 
\end{equation}
Moreover, if $\{S_i:\,i\in[k]\}$ is additionally assumed to commute, then there exists $C_{r,h_1,\dotsc,h_k}$ such that
\begin{equation}\label{Unmultosc}
\sup_{J\in\mathbb{N}}\sup_{I\in\mathfrak{S}_J(\mathbb{N}^k)}\|O^2_{I,J}(A_{\mathbf{N}} f:\,\mathbf{N}\in\mathbb{N}^k)\|_{L_{\mu}^r(X)}\le C_{r,h_1,\dotsc,h_k}\|f\|_{L_{\mu}^r(X)}\quad\text{for any $f\in L_{\mu}^r(X)$}\text{.}
\end{equation}
\end{theorem}
Although we will not work directly with multiparameter oscillations we give the definition of $\mathfrak{S}_J(\mathbb{N}^k)$ and $O^2_{I,J}(A_{\mathbf{N}} f:\,\mathbf{N}\in\mathbb{N}^k)$ in the notation subsection in the end of the introduction.

The existence of the multiparameter limit can be deduced via an iterative argument using the one-dimensional case. More precisely, it can be obtained by exploiting the fact that the limit of the single averages exists, which is a corollary of the estimate $\ref{U2OE}$, together with the fact that we have $L^r_{\mu}\to L^r_{\mu}$ estimates for the corresponding maximal function,
\[
\text{i.e.: }\big\|\sup_{N\in\mathbb{N}}A_{N}|f|\big\|_{L_{\mu}^r(X)}\le C_{r,h}\|f\|_{L_{\mu}^r(X)}\text{,}
\]
which is immediate from $\ref{VVMFE}$.

We note that the estimate $\ref{Unmultosc}$ implies that the limit in $\ref{multpointwise}$ exists, see Remark~2.4 together with Proposition~2.8, page 15 in \cite{MOE}. One may deduce the uniform multiparameter 2-oscillation estimate $\ref{Unmultosc}$ as a corollary of Theorem~$\ref{OM}$ together with Proposition~4.1 in \cite{MOE}. To upgrade the simple existence of the multiparameter limit to the powerful quantitative estimate $\ref{Unmultosc}$, the assumption of commutativity is crucial. 

The technical heart of the present work is the following exponential sum estimate.

\begin{proposition}\label{MainExpSumEst}
Assume $c\in(1,4/3)$, $h\in\mathcal{R}_c$ with $\varphi$ its compositional inverse, and let $\varepsilon>0$. Then there exists a positive constant $C=C(h,\varepsilon)$ such that for all $\xi\in\mathbb{T}$ and $N\in\mathbb{N}$ we have 
\begin{equation}\label{ExpEstM}
\Big|\sum_{p\in\mathbb{P}\cap[1,N]}\log(p)e^{2\pi i \lfloor h(p) \rfloor \xi }-\sum_{n\le h(N)}\varphi'(n)e^{2\pi i n \xi }\Big|\le CNe^{-(\log N)^{1/3-\varepsilon}}\text{.}
\end{equation}
\end{proposition}
This estimate for the case $c\in(1,17/16)$ and $h(x)=x^c$ can be deduced from the work done in \cite{laporta}, the proof of which relies heavily on estimates already appearing in \cite{Tolev}. We manage to obtain the range $c\in(1,4/3)$, which is usually the limit of arguments of that nature involving van der Corput type estimates, for example, see Proposition~3.12 in \cite{LPE}, where the exponential sums are much simpler but the condition $c<4/3$ is present. We make further comments in subsection~$\ref{strategy}$.
\subsection{Waring-type problem for $\lfloor h_1(p_1)\rfloor+\lfloor h_2(p_2)\rfloor+\lfloor h_3(p_3)\rfloor$}
As a corollary, of our investigations of such exponential sum estimates we derive the following Waring-type result. For every $i\in[3]$ let $c_i\in(1,2)$ and $h_i\in\mathcal{R}_{c_i}$, with $\varphi_i$ its compositional inverse and $
\gamma_i=1/c_i$. For each $\lambda\in\mathbb{N}$, we define
\[
R_{h_1,h_2,h_3}(\lambda)=\sum_{\substack{p_1,p_2,p_3\in\mathbb{P},\\ \lfloor h_1(p_1)\rfloor+\lfloor h_2(p_2)\rfloor+\lfloor h_3(p_3)\rfloor=\lambda}}\log(p_1)\log(p_2)\log(p_3)
\]
and 
\[
r_{h_1,h_2,h_3}(\lambda)=\big|\big\{(m_1,m_2,m_3)\in\mathbb{N}^3:\,\lfloor h_1(m_1)\rfloor+\lfloor h_2(m_2)\rfloor+\lfloor h_3(m_3)\rfloor=\lambda\big\}\big|\text{.}
\]
\begin{theorem}\label{nondiag}For every $i\in[3]$ let $c_i\in(1,4/3)$ and $h_i\in\mathcal{R}_i$, and assume that
\begin{equation}\label{assumptiongamma}
4(1-\gamma_1)+\frac{45}{4}(1-\gamma_2)+\frac{45}{4}(1-\gamma_3)<1\text{.}
\end{equation}
Then we have that
\[
R_{h_1,h_2,h_3}(\lambda)=\frac{\Gamma(\gamma_1)\Gamma(\gamma_2)\Gamma(\gamma_3)}{\Gamma(\gamma_1+\gamma_2+\gamma_3)}\lambda^{2}\varphi_1'(\lambda)\varphi_2'(\lambda)\varphi_3'(\lambda)+o(\lambda^{2}\varphi_1'(\lambda)\varphi_2'(\lambda)\varphi_3'(\lambda))\text{,}
\]
and moreover for every $\varepsilon>0$ we have
\begin{equation}\label{SecondEstimatekey}
R_{h_1,h_2,h_3}(\lambda)=r_{h_1,h_2,h_3}(\lambda)+O_{\varepsilon}\Big(\lambda^{2}\varphi_1'(\lambda)\varphi_2'(\lambda)\varphi_3'(\lambda)e^{-(\log \lambda)^{\frac{1}{3}-\varepsilon}}\Big)\text{.}
\end{equation}
\end{theorem}
Of course the implicit constant in the latter assertion above may also depend on $h_1$, $h_2$, $h_3$. The reader is encouraged to compare this theorem with the main theorem in \cite{laporta}, as well as Theorem~1.2 from \cite{LPE}. Let us note that in \cite{laporta} the range of $c$ is $(1,17/16)$, and the assumption $\ref{assumptiongamma}$ for $h=h_1=h_2=h_3$ becomes $c\in(1,53/51)$. This comes from the fact that a different optimization took place for our proof of the exponential sum estimates. See also the discussion in subsection~$\ref{strategy}$. Apart from extending the work of \cite{laporta} to a wider class of functions for small exponents, our theorem also covers the ``non-diagonal'' case, namely the case where $h_1$, $h_2$, $h_3$ are different.

Finally, let us comment that the second estimate above shows that our argument does not lose quantitative control of the error for treating the restriction over the primes but rather for the estimate
\begin{equation}
r_{h_1,h_2,h_3}(
\lambda)=\frac{\Gamma(\gamma_1)\Gamma(\gamma_2)\Gamma(\gamma_3)}{\Gamma(\gamma_1+\gamma_2+\gamma_3)}\lambda^{2}\varphi_1'(\lambda)\varphi_2'(\lambda)\varphi_3'(\lambda)+o(\lambda^{2}\varphi_1'(\lambda)\varphi_2'(\lambda)\varphi_3'(\lambda))\text{.}
\end{equation}
This is due to the complications the class $\mathcal{R}_c$ brings, and is in some sense optimal, see Remark~3.6 in \cite{LPE}. The logarithmically weighted count of the solutions over the primes is the same as the unrestricted count with a quantifiable error.

\subsection{Strategy of our proofs}\label{strategy}
In this short subsection we give a brief description of the strategy of our proofs. As discussed earlier, to establish theorems $\ref{basic}\text{, }\ref{MT}\text{, }\ref{Multi}$ it suffices to establish Theorem~$\ref{OM}$, i.e.: to obtain uniform $2$-oscillation and vector-valued maximal estimates. By Calder\'on's transference principle it suffices to establish these on the integer shift system. We have adopted the strategy from \cite{WT11LD} and we encourage the reader to see the discussion after Theorem~1.12 from the aforementioned paper. Instead of directly repeating the arguments here, we have formulated two abstract lemmas allowing us:
\begin{enumerate}
\item[i)] on the one hand to compare families of averaging operators with sufficiently good $L^{\infty}(\mathbb{T})$-bounds on the difference of their corresponding multipliers, see Lemma~$\ref{Comparison}$ and 
\item[ii)] on the other hand to perform delicate summation by parts arguments to pass from one averaging operator to the other in the context of oscillations, see Lemma~$\ref{correctosccomp}$. 
\end{enumerate}
Apart from directly treating the so-called short oscillations for a suitable variant of our averaging operator, applying these two steps and taking into account Proposition~$\ref{MainExpSumEst}$ allows us to reduce our problem to the corresponding one for the following averaging operators
\[
H_Nf(x)=\frac{1}{N}\sum_{n\le h(N)}\varphi'(n)f(x-n)\text{,}
\]
which are much better-behaving and have been studied in the context of these problems in \cite{WT11LD}. 

For the exponential sum estimates certain inputs from analytic number theory are needed and the strategy is the following: We perform the circle method with one major arc around $0$. Our departure point is the work from \cite{laporta} and \cite{Tolev}, where an analogous exponential sum estimate for the case $h(x)=x^c$ is established. Notably, our treatment deviates from the one in \cite{Tolev} in the minor arc analysis and we obtain the range $c\in(1,4/3)$.

For the minor arc we use certain truncated fourier series expansions with uniform bounds on the tail to remove the floor function in a similar manner to \cite{LPE}. The main tool for estimating exponential sums is van der Corput. The argument also utilizes Vaughan's identity in a similar manner to \cite{MMR} and \cite{RHLMLD}. Our arguments for the minor arc yield polynomial saving in $N$, see Lemma~$\ref{Minor}$, which although does not play any role in the proof of Theorem~$\ref{OM}$, it is indispensable for the proof of Theorem~$\ref{nondiag}$. In fact, if one wishes to obtain the best range for exponents in Theorem~$\ref{nondiag}$, certain choices made in \cite{laporta} are more suitable, since one would have to balance the size of the polynomial saving in the minor arc analysis with the allowed largeness of $c$. Our choices are made to optimize the range of our exponents for our ergodic theorems, and indeed we manage to obtain the rather large range $(1,4/3)$. 

For the major arc, after passing to a dyadic variant of our exponential sum, see Lemma~$\ref{majarckey}$, we apply the (truncated) explicit formula, which after summation and integration by parts immediately gives the main term of our approximation. To conclude, we are left with the task of estimating certain sums over roots of the zeta function, see Lemma~$\ref{Keyroots}$. This is the only part of the argument yielding a subpolynomial saving. For the application of our exponential sum estimates to obtain Theorem~$\ref{nondiag}$, it is important that we have allowed for large major arcs, see Lemma~$\ref{Major}$.

For the proof of Theorem~$\ref{nondiag}$, we begin by expressing $R_{h_1,h_2,h_3}$ as an integral over the corresponding multipliers. We cannot compare exponential sums directly, since our errors do not have polynomial decay for all the frequencies, so the analysis is split into major and minor arcs. A crucial observation is that the major and minor arcs do not have to be adapted precisely for all three functions $h_1,h_2,h_3$. The flexibility for choosing the parameters specifying admissible major arcs in Lemma~$\ref{Major}$ and the precise quantification of $\chi$ in Lemma~$\ref{Minor}$ allow us to conclude by adapting certain techniques from \cite{Tolev}, to accommodate for the sophisticated nature of $c$-regularly varying functions and the complications they introduce.

We have organized the paper as follows. In Section~$\ref{PET}$, we establish Theorem~$\ref{OM}$ assuming our exponential sum estimates hold. In Section~$\ref{ExpEstSec}$ we establish the exponential sum estimates, namely, we prove Proposition~$\ref{MainExpSumEst}$ and in the final section, we prove Theorem~$\ref{nondiag}$.

\subsection{Some remarks on the convergence on $L^1$} We wish to finish our introduction with certain comments regarding the natural question about pointwise convergence of our ergodic averages on $L^1$. Understanding the behavior of ergodic averages modeled along sparse sequences for functions in $L^1$ in the context of pointwise convergence is a rather challenging task. For more than fifteen years after Bourgain's seminal work \cite{bg2},\cite{bg3},\cite{BET} on pointwise ergodic theorems along polynomial and prime orbits there were no examples of sparse sequences for which the corresponding ergodic averaging operators converged pointwise for every function in $L^1$ and every dynamical system. This led Rosenblatt and Wierdl \cite[Conjecture~4.1]{RWC} to formulate a famous conjecture which, loosely speaking, asserts that such a result is impossible for sparse orbits. Several counterexamples have been given for this conjecture, for brief historical remarks we refer the reader to the introduction of \cite{WT11LD}. 

Notably, no example involving prime numbers has been given, and for the prime orbits pointwise convergence is known to fail on $L^1$, see \cite{N^kfail}. The orbits $\lfloor h(p)\rfloor$ are much better-behaving than the primes themselves and we believe that the ergodic averages in $\ref{ErgAv}$ do converge pointwise even for $f\in L_{\mu}^1(X)$ for small exponents.

Calder\'on's transference principle together with the results of the present work suggest that it would suffice to establish the following estimate
\[
\Big|\Big\{x\in\mathbb{Z}\,:\,\sup_{N\in\mathbb{N}}\frac{1}{N}\sum_{p\in\mathbb{P}\cap[1,N]}\log(p)|f|(x-\lfloor h(p) \rfloor)>\lambda\Big\}\Big|\le C \lambda^{-1}\|f\|_{\ell^1(\mathbb{Z})}\text{.}
\]
According to Proposition~$\ref{MainExpSumEst}$, the appropriate approximant for the corresponding multiplier is 
\[
\frac{1}{N}\sum_{n\le h(N)}\varphi'(n)e^{2\pi i n \xi }\text{,}
\]
which is precisely the same approximant used for establishing pointwise convergence for the following ergodic averages $
\frac{1}{N}\sum_{n\le N}f(T^{\lfloor h(n)\rfloor}x)$, which do converge even on $L^1$, see \cite{W11}. Moreover, if one assumes that the Riemann zeta function has a zero-free strip, namely, that there exists $\delta\in(0,1/2)$ such that $\zeta(s)\neq 0$ for every $s$ with $\Ree(s)\in[1-\delta,1]$, then our exponential sum estimates can be upgraded to give polynomial saving, which is the same type of saving as the one achieved in \cite{W11}.   

Establishing weak type (1,1) bounds for the maximal function corresponding to such averages goes beyond exponential sum estimates of the form $\ref{ExpEstM}$, see for example Section~5 and Section~6 from \cite{W11}. Nevertheless, we do believe that the weak type $(1,1)$ bound holds here and moreover we believe that a conditional result for small exponents may not be entirely out of reach.

\subsection{Notation} We denote by $C$ a positive constant that may change from occurrence to occurrence. If $A,B$ are two non-negative quantities, we write $A\lesssim B$ or $B \gtrsim A$  to denote that there exists a positive constant $C$ such that $A\le C B$. We write $A\simeq B$ to indicate that $A\lesssim B$ and $A\gtrsim B$. For two complex-valued functions $f,g$ we write $f\sim g$ to denote that $\lim_{x\to\infty}\frac{f(x)}{g(x)}=1$. In the sequel, $e(x)$ stands for $e^{2\pi i x}$. We also make use of the standard convention that if we write ``$p\le N$'' under the summation symbol, the summation is taken over primes. We adopt the notation from \cite{Tolev} and whenever we write 
\[
\sum_{0<\gamma <T} 
\]
 we mean that the summation takes place over all the nontrivial roots of the zeta function with imaginary part $\gamma\in(0,T)$. In general $\rho=\beta+i\gamma$, $\beta\in(0,1)$ and $\gamma\in\mathbb{R}$ denotes a nontrivial zero of the Riemann zeta function and any summation with restrictions written on those parameters means that we sum over all such roots. For any $t\in[1,\infty)$ we denote by $\mathbb{N}_{\ge t}$ the set $\{n\in\mathbb{N}:\,n\ge t\}$, and we also let $[t]\coloneqq\{1,2,\dotsc,\lfloor t\rfloor\}$. We note that $h(x)$ is not defined for $x\le x_0$ but we abuse notation, see for example $\ref{ExpEstM}$. We can let $h(x)$ take arbitrary values for $x\le x_0$ and all our results and estimates will remain the same. In a similar spirit, certain averages make sense for large values of $N$, we again choose to abuse notation to keep the exposition more reasonable.
 
 We end this subsection by giving the definition of multiparameter oscillations. For any $x$, $y\in\mathbb{R}^k$ we write $x\prec_s y$ if $x_i<y_i$ for all $i\in[k]$. For every $J\in\mathbb{N}$ we define the set
 \[
 \mathfrak{S}_J(\mathbb{N}^k)\coloneqq \big\{\{t_0,\dotsc,t_J\}\subseteq \mathbb{N}^k:\,t_0\prec_s\dots\prec_s t_J\big\}\text{,}
 \]
 i.e.: $\mathfrak{S}_J(\mathbb{N}^k)$ contains all the strictly increasing sequences (with respect to $\prec_s$) of length $J+1$ taking values in $\mathbb{N}^k$. For every $J\in\mathbb{N}$, $I\in \mathfrak{S}_J(\mathbb{N}^k)$ and every $k$-parameter family of complex-valued functions $\big(a_t(x):\,t\in\mathbb{N}^k\big)$ defined an a set $X$, we define 
\[
 O^2_{I,J}(a_{\mathbf{N}}(x):\,\mathbf{N}\in\mathbb{N}^k)=\Big(\sum_{j=0}^{J-1}\sup_{t\in\mathbb{B}[I_j]}|a_t(x)-a_{I_j}(x)|^2\Big)^{1/2}\text{,}
\] 
where $\mathbb{B}[I_i]\coloneqq [I_{i1},I_{(i+1)1})\times\dots\times[I_{ik},I_{(i+1)k})\cap\mathbb{N}^k$ can be thought of as a ``box'' in $\mathbb{N}^k$ with sides parallel to the axes defined by its ``diagonal'' points $I_i$ and $I_{i+1}$.  
 \section*{Acknowledgment}We would like to thank Mariusz Mirek for proposing the problem and for his constant support, as well as Alex Kontorovich for helpful discussions.

Part of this research was conducted during the 2024 DIMACS REU program at Rutgers University, CNS-2150186. All authors were partly supported by the Rutgers Department of Mathematics and the last three authors were partly supported by the National Science Foundation under Grant DMS-2154712.
 \section{Oscillation and Vector-valued maximal estimates}\label{PET}
 Here we give a proof of Theorem~$\ref{OM}$ assuming Proposition~$\ref{MainExpSumEst}$. Our proof relies on comparison with a better-behaving averaging operator corresponding to the approximant multiplier. The treatment is very similar to the one of sections~3 and 4 from \cite{WT11LD}. We begin by formulating an abstract principle which carries out this comparison.

\begin{lemma}\label{Comparison}
Let $\big(A_N\colon\ell^{\infty}(\mathbb{Z})\to\ell^{\infty}(\mathbb{Z})\big)_{N\in\mathbb{N}}$ and $\big(B_N\colon\ell^{\infty}(\mathbb{Z})\to\ell^{\infty}(\mathbb{Z})\big)_{N\in\mathbb{N}}$ be two families of linear operators such that for all $r\in(1,\infty)$ we have
\[
\sup_{N\in\mathbb{N}}\|A_N\|_{\ell^r(\mathbb{Z})\to\ell^r(\mathbb{Z})}<\infty\quad\text{and}\quad\sup_{N\in\mathbb{N}}\|B_N\|_{\ell^r(\mathbb{Z})\to\ell^r(\mathbb{Z})}<\infty\text{.}
\]  
Assume that for any $\rho\in(0,2]$ we have $\sum_{N\in\mathbb{N}}
\|A_N-B_N\|^{\rho}_{\ell^2(\mathbb{Z})\to\ell^2(\mathbb{Z})}<\infty$. 
Then for any $r\in(1,\infty)$ we have that 
\begin{equation}\label{OSCcomp}
\sup_{J\in\mathbb{N}}\sup_{I\in\mathfrak{S}_J(\mathbb{N})}\|O^2_{I,J}(A_N f:\,N\in\mathbb{N})\|_{\ell^r(\mathbb{Z})}\lesssim_r\sup_{J\in\mathbb{N}}\sup_{I\in\mathfrak{S}_J(\mathbb{N})}\|O^2_{I,J}(B_N f:\,N\in\mathbb{N})\|_{\ell^r(\mathbb{Z})}
+\|f\|_{\ell^r(\mathbb{Z})}\end{equation}
and
\begin{equation}\label{VVcomp}
\Big\|\Big(\sum_{j\in\mathbb{Z}}\big(\sup_{N\in\mathbb{N}}A_{N}|f_j|\big)^2\Big)^{1/2}\Big\|_{\ell^r(\mathbb{Z})}\lesssim_{r}\Big\|\Big(\sum_{j\in\mathbb{Z}}\big(\sup_{N\in\mathbb{N}}B_{N}|f_j|\big)^2\Big)^{1/2}\Big\|_{\ell^r(\mathbb{Z})}+\Big\|\big(\sum_{j\in\mathbb{Z}}|f_j|^2\big)^{1/2}\Big\|_{\ell^r(\mathbb{Z})}\text{.}
\end{equation}
\end{lemma}
\begin{proof}
For convenience let $C_N\coloneq \|A_N-B_N\|_{\ell^2(\mathbb{Z})\to\ell^2(\mathbb{Z})}$. We begin with the vector-valued maximal estimates. Note that
\[
\bigg\|\Big(\sum_{j\in\mathbb{Z}}\big(\sup_{N\in\mathbb{N}}A_{N}|f_j|\big)^2\Big)^{1/2}\bigg\|_{\ell^r(\mathbb{Z})}\]
\[
\le\bigg\|\Big(\sum_{j\in\mathbb{Z}}\big(\sup_{N\in\mathbb{N}}B_{N}|f_j|\big)^2\Big)^{1/2}\bigg\|_{\ell^r(\mathbb{Z})}+\bigg\|\Big(\sum_{j\in\mathbb{Z}}\sum_{N\in\mathbb{N}}\big|A_{N}|f_j|-B_{N}|f_j|\big|^2\Big)^{1/2}\bigg\|_{\ell^r(\mathbb{Z})}\text{,}
\]
so it only remains to show that 
\begin{equation}\label{KeyEstimate}
\bigg\|\Big(\sum_{j\in\mathbb{Z}}\sum_{N\in\mathbb{N}}\big|A_{N}|f_j|-B_{N}|f_j|\big|^2\Big)^{1/2}\bigg\|_{\ell^r(\mathbb{Z})}\lesssim_r \Big\|\big(\sum_{j\in\mathbb{Z}}|f_j|^2\big)^{1/2}\Big\|_{\ell^r(\mathbb{Z})}\text{.}
\end{equation}
For $r=2$ we note that 
\[
\bigg\|\Big(\sum_{j\in\mathbb{Z}}\sum_{N\in\mathbb{N}}\big|A_{N}|f_j|-B_N|f_j|\big|^2\Big)^{1/2}\bigg\|_{\ell^2(\mathbb{Z})}\le\Big(\sum_{j\in\mathbb{Z}}\sum_{N\in\mathbb{N}}C_N^2\|f_j\|^2_{\ell^2(\mathbb{Z})}\Big)^{1/2}\lesssim\Big\|\Big(\sum_{j\in\mathbb{Z}}|f_j|^2\Big)^{1/2}\Big\|_{\ell^2(\mathbb{Z})} \text{,}
\]
where we have used our assumption for $\rho=2$.

For $r\in(2,\infty)$, let us fix a $r_0=2r$ and by our assumptions there exists a positive constant $K_{r}$ such that $\|A_Nf-B_Nf\|_{\ell^{r_0}(\mathbb{Z})}\le K_{r}\|f\|_{\ell^{r_0}(\mathbb{Z})}$. Choose $\theta\in (0,1)$ such that $\frac{1}{r}=\frac{\theta}{2}+\frac{1-\theta}{r_0}$, and note that Riesz--Thorin interpolation yields
\begin{equation}\label{ERES}
\|A_Nf-B_Nf\|_{\ell^{r}(\mathbb{Z})}\le C_N^{\theta}K_{r}^{1-\theta}\|f\|_{\ell^{r}(\mathbb{Z})}\lesssim_r C_N^{\theta}\|f\|_{\ell^{r}(\mathbb{Z})}\text{.}
\end{equation}
Thus if we let $T_N\colon \ell^{r}(\mathbb{Z})\to \ell^r(\mathbb{Z})$ such that $T_Nf=A_Nf-B_Nf$, then we know that $T_N$ is a bounded linear operator with $\|T_N\|_{\ell^r(\mathbb{Z})\to\ell^r(\mathbb{Z})}\lesssim_r C_N^\theta $. Thus $T_N$ has an $\ell^2$-valued extension with the same norm, see Theorem~5.5.1 in \cite{Graf}, namely
\begin{equation}\label{VVuse}
\Big\|\big(\sum_{j\in\mathbb{Z}}|T_N(f_j)|^2\big)^{1/2}\Big\|_{\ell^r(\mathbb{Z})}\lesssim_r C_N^\theta \Big\|\big(\sum_{j\in\mathbb{Z}}|f_j|^2\big)^{1/2}\Big\|_{\ell^r(\mathbb{Z})}\text{.}
\end{equation}
Therefore we may estimate as follows
\[
\bigg\|\Big(\sum_{j\in\mathbb{Z}}\sum_{N\in\mathbb{N}}\big|A_N|f_j|-B_N|f_j|\big|^2\Big)^{1/2}\bigg\|_{\ell^r(\mathbb{Z})}\le\Big(\sum_{N\in\mathbb{N}}\Big\|\Big(\sum_{j\in\mathbb{Z}}\big|T_N|f_j|\big|^2\Big)^{1/2}\Big\|_{\ell^r(\mathbb{Z})}^2\Big)^{1/2}
\]
\[
\lesssim_r\Big(\sum_{N\in\mathbb{N}}C_N^{2\theta}\Big)^{1/2}\Big\|\big(\sum_{j\in\mathbb{Z}}|f_j|^2\big)^{1/2}\Big\|_{\ell^r(\mathbb{Z})}\lesssim_{r} \Big\|\big(\sum_{j\in\mathbb{Z}}|f_j|^2\big)^{1/2}\Big\|_{\ell^r(\mathbb{Z})} \text{,}
\] 
where we used the fact that $r/2>1$ and our assumption with $\rho=2\theta\in(0,2)$. 

For $r\in (1,2)$, we choose $r_0=r-\frac{r-1}{2}\in (1,r)$, and $\theta\in (0,1)$ such that $\frac{1}{r}=\frac{\theta}{2}+\frac{1-\theta}{r_0}$ and the same reasoning as before yields $\ref{VVuse}$. We proceed as follows
\[
\bigg\|\Big(\sum_{j\in\mathbb{Z}}\sum_{N\in\mathbb{N}}\big|A_N|f_j|-B_N|f_j|\big|^2\Big)^{1/2}\bigg\|_{\ell^r(\mathbb{Z})}\le\Big(\sum_{N\in\mathbb{N}}\Big\|\Big(\sum_{j\in\mathbb{Z}}\big|T_N|f_j|\big|^2\big)^{1/2}\Big\|_{\ell^r(\mathbb{Z})}^r\Big)^{1/r}
\]
\[
\lesssim_r\Big( \sum_{N\in\mathbb{N}} C_N^{\theta r}\Big)^{1/r}\big\|\big(\sum_{j\in\mathbb{Z}}|f_j|^2\big)^{1/2}\big\|_{\ell^r(\mathbb{Z})}\lesssim_r \big\|\big(\sum_{j\in\mathbb{Z}}|f_j|^2\big)^{1/2}\big\|_{\ell^r(\mathbb{Z})}\text{,}
\]
where for the first estimate we used the fact that $r<2$, and for the last one we used our assumption for $\rho=\theta r\in(0,2)$.

We have established the estimate $\ref{KeyEstimate}$ for all $r\in(1,\infty)$, and thus we have obtained the vector-valued maximal estimate $\ref{VVcomp}$. 

We now turn our attention to $\ref{OSCcomp}$. Firstly, we note that the following pointwise estimate holds
\[
O^2_{I,J}(A_N f(x):\,N\in\mathbb{N})\lesssim O^2_{I,J}(B_N f(x):\,N\in\mathbb{N})+\Big(\sum_{N\in\mathbb{N}}|A_N f(x)-B_N(x)|^2\Big)^{1/2}\text{.} 
\]
Therefore we immediately obtain
\[
\sup_{J\in\mathbb{N}}\sup_{I\in\mathfrak{S}_J(\mathbb{N})}\|O^2_{I,J}(A_N f:\,N\in\mathbb{N})\|_{\ell^r(\mathbb{Z})}
\]
\[
\lesssim\sup_{J\in\mathbb{N}}\sup_{I\in\mathfrak{S}_J(\mathbb{N})}\|O^2_{I,J}(B_N f:\,N\in\mathbb{N})\|_{\ell^r(\mathbb{Z})}+\Big\|\Big(\sum_{N\in\mathbb{N}}|A_N f(x)-B_Nf(x)|^2\Big)^{1/2}\Big\|_{\ell^r(\mathbb{Z})}\text{,}
\]
and the second summand can be appropriately bounded using $\ref{KeyEstimate}$, yielding the desired estimate.
\end{proof}

\begin{remark} Let us note that even if the condition $\sum_{N\in\mathbb{N}}
\|A_N-B_N\|^{\rho}_{\ell^2(\mathbb{Z})\to\ell^2(\mathbb{Z})}<\infty$ is available for a smaller range for $\rho$, an argument as the one presented above may still yield the estimates $\ref{OSCcomp}$, $\ref{VVcomp}$ but for a smaller range for $r$.
\end{remark}
We now formulate a lemma allowing us to perform summation by parts arguments in the context of oscillations.
\begin{lemma}\label{correctosccomp}Let $\big(A_N\colon\ell^{\infty}(\mathbb{Z})\to\ell^{\infty}(\mathbb{Z})\big)_{N\in\mathbb{N}}$ be operators such that $\sup_{N\in\mathbb{N}}\|A_N\|_{\ell^{\infty}(\mathbb{Z})\to\ell^{\infty}(\mathbb{Z})}<\infty$ and let $(\lambda_s^k)_{s,k\in\mathbb{N}}$ be non-negative real numbers satisfying the following: 
\begin{itemize}
\item[i)] There exists $\Lambda\in(0,\infty)$ such that for every $k$ we have that $\sum_{s=1}^{\infty}\lambda^k_s=\Lambda$.
\item[ii)] For any fixed $N\in\mathbb{N}$ we have that $\sum_{s=1}^{N}\lambda_{s}^k$ is decreasing in $k$.
\end{itemize}
Then for any $r\in[1,\infty)$, $q\in[2,\infty)$, $J\in\mathbb{N}$ and $f\in\ell^{r}(\mathbb{Z})$ we have
\[
\sup_{I\in\frak{S}_J(\mathbb{N})}\Big\|O^q_{I,J}\Big(\sum_{s=1}^\infty\lambda_{s}^k A_sf(x):\,k\in\mathbb{N}\Big)\Big\|_{\ell^r(\mathbb{Z})}\le \Lambda\sup_{I\in\frak{S}_J(\mathbb{N})}\big\|O^q_{I,J}\big(A_kf(x):\,k\in\mathbb{N}\big)\big\|_{\ell^r(\mathbb{Z})}\text{.}
\] 
\end{lemma}
\begin{proof}
An analogous argument to the one of Section~3, pages 21-23 in \cite{WT11LD}, which in turn is a variant of the  argument given in Lemma~2 from \cite{UA}, gives the result.
\end{proof}
We are now ready to give a proof of Theorem~\ref{OM}. 
\begin{proof}[Proof of Theorem~$\ref{OM}$]
We begin with the estimate $\ref{VVMFE}$. By Calder\'on's transference principle it suffices to prove that for any $r\in(1,\infty)$, there exists $C=C(r,h)$ such that
\begin{equation}\label{VVMCT}
\Big\|\Big(\sum_{j\in\mathbb{Z}}\big(\sup_{N\in\mathbb{N}_{\ge 2}}A_{N}|f_j|\big)^2\Big)^{1/2}\Big\|_{\ell^r(\mathbb{Z})}\le C_{r,h}\Big\|\big(\sum_{j\in\mathbb{Z}}|f_j|^2\big)^{1/2}\Big\|_{\ell^r(\mathbb{Z})}\text{ for any $(f_j)_{j\in\mathbb{Z}}\in \ell^r\big(\mathbb{Z};\ell^2(\mathbb{Z})\big)$}\text{,}
\end{equation}
where $A_tf(x)=\frac{1}{|\mathbb{P}\cap[t]|}
\sum_{p\le t}f(x-\lfloor h(p)\rfloor)$, $t\in[2,\infty)$. Let us define two auxiliary operators 
\[
M_tf(x)=\frac{1}{t}\sum_{p\le t}\log(p)f(x-\lfloor h(p)\rfloor)\quad \text{and}\quad H_tf(x)=\frac{1}{t}\sum_{n\le h(t)}\varphi'(n)f(x-n)\text{, }t\in[1,\infty)\text{.}
\]
A standard summation by parts argument together with a passage to dyadic scales yield the pointwise estimate 
\begin{equation}\label{easyptwise}
\sup_{N\in\mathbb{N}_{\ge 2}}A_N|f|(x)\lesssim \sup_{n\in\mathbb{N}}M_{2^n}|f|(x)\text{.} 
\end{equation}
For the convenience of the reader we provide some details here. Summation by part yields
\[
A_N|f|(x)\lesssim\frac{N\log(N)^{-1}}{|\mathbb{P}\cap[N]|}M_N|f|(x)+\frac{1}{|\mathbb{P}\cap[N]|}\sum_{m=2}^{N-1}M_m|f|(x)m(\log(m)^{-1}-\log(m+1)^{-1})
\]
\[
\lesssim M_N|f|(x)+\frac{\log N}{N}\Big(\sum_{m=2}^{N-1}m(\log(m)^{-1}-\log(m+1)^{-1})\Big)\sup_{m\in\mathbb{N}_{\ge 2}}M_m|f|(x)\lesssim\sup_{m\in\mathbb{N}_{\ge 2}}M_m|f|(x)\text{.}
\]
Thus $\sup_{N\in\mathbb{N}_{\ge 2}}A_N|f|(x)\lesssim \sup_{N\in\mathbb{N}_{\ge 2}}M_{N}|f|(x)\le2\sup_{m\in\mathbb{N}}M_{2^m}|f|(x)$, this establishes the estimate $\ref{easyptwise}$.

We will apply Lemma~$\ref{Comparison}$. For any $t\in\mathbb{N}$ we have by Plancherel together with an application of Proposition~$\ref{MainExpSumEst}$ with $\varepsilon=1/12$ the following estimate
\[
\|M_{t}f-H_{t}f\|_{\ell^2(\mathbb{Z})}=\Big\|\Big(\frac{1}{t}\sum_{p\le t}\log(p)e(\lfloor h(p)\xi\rfloor)-\frac{1}{t}\sum_{n\le h(t)}\varphi'(n)e(n\xi)\Big)\hat{f}(\xi)\Big\|_{L_{d\xi}^2(\mathbb{T})}\lesssim e^{-(\log t)^{1/4}}
\|f\|_{\ell^2(\mathbb{Z})}\text{.}
\]
Note that for any $\rho>0$ we have that
\[
\sum_{N\in\mathbb{N}}\|M_{2^N}f-H_{2^N}f\|^{\rho}_{\ell^2(\mathbb{Z})\to\ell^2(\mathbb{Z})}\lesssim \sum_{N\in\mathbb{N}}e^{-\rho(\log^{1/4} 2) N^{1/4}}<\infty\text{,}
\]
and $\|M_tf\|_{\ell^r(\mathbb{Z})}\lesssim \|f\|_{\ell^r(\mathbb{Z})} $. Using the basic properties of $\varphi$ is it easy to see that 
\[
\|H_tf\|_{\ell^r(\mathbb{Z})}\le \frac{1}{t}\sum_{n\le h(t)}\varphi'(n)\|f\|_{\ell^r(\mathbb{Z})}\lesssim\|f\|_{\ell^r(\mathbb{Z})}\text{,}
\]
see for example Step 4 from the proof~\ref{Steps} (we apply it here for $\xi=0$). Thus Lemma~$\ref{Comparison}$ is applicable for $M_{2^N}$ and $H_{2^N}$ and we obtain
\[
\Big\|\Big(\sum_{j\in\mathbb{Z}}\big(\sup_{N\in\mathbb{N}_{\ge 2}}A_{N}|f_j|\big)^2\Big)^{1/2}\Big\|_{\ell^r(\mathbb{Z})}\lesssim \Big\|\Big(\sum_{j\in\mathbb{Z}}\big(\sup_{N\in\mathbb{N}}M_{2^N}|f_j|\big)^2\Big)^{1/2}\Big\|_{\ell^r(\mathbb{Z})}
\]
\begin{equation}\label{finalstepinVV}
\lesssim_r\Big\|\Big(\sum_{j\in\mathbb{Z}}\big(\sup_{N\in\mathbb{N}}H_{2^N}|f_j|\big)^2\Big)^{1/2}\Big\|_{\ell^r(\mathbb{Z})}+\Big\|\big(\sum_{j\in\mathbb{Z}}|f_j|^2\big)^{1/2}\Big\|_{\ell^r(\mathbb{Z})}\lesssim_r \Big\|\big(\sum_{j\in\mathbb{Z}}|f_j|^2\big)^{1/2}\Big\|_{\ell^r(\mathbb{Z})} 
\end{equation}
where for the last estimate we note that for a variant of the operator $H_t$ these vector-valued maximal estimates have already been established in \cite{WT11LD}, see page 27. To keep the exposition reasonably self-contained we provide a short proof. For any $n\in\mathbb{N}_0$ let $l=l(n)\in\mathbb{N}_0$ be such that $2^{l-1}\le h(2^n)<2^l$ and note
\[
H_{2^n}|f|(x)\le\frac{1}{2^n}\sum_{k=0}^{l}\sum_{2^k \le m<2^{k+1}}\varphi'(m)|f|(x-m)\lesssim\frac{1}{2^n}\sum_{k=0}^l 2^k\varphi'(2^k)\Big(\frac{1}{2^k}\sum_{2^k \le m<2^{k+1}}|f|(x-m)\Big)
\]
\begin{equation}\label{geom}
\lesssim\Big(\frac{1}{2^n}\sum_{k=0}^l \varphi(2^k)\Big)\sup_{k\in\mathbb{N}_0}\frac{1}{2^k}\sum_{1\le m\le  2^{k}}|f(x-m)|\lesssim \sup_{k\in\mathbb{N}_0}\frac{1}{2^k}\sum_{1\le m\le  2^{k}}|f(x-m)|\text{,}
\end{equation}
where for the last estimate we note that the behavior of this sum is of geometric type in the sense that the last term dictates its growth. We provide the details here. There exists $d=d(h)\in(0,1)$ and $x_0=x_0(h)$ such that $\varphi(x)\le d\varphi(2x)$ for $x\ge x_0$. To see this note that by Lemma~2.6 \cite{MMR} we get that
\[
\frac{\varphi(x)}{\varphi(2x)}=2^{-\gamma}e^{-\int_{x}^{2x}\frac{\theta(t)}{t}dt}\text{, where $\theta(t)\to 0$ as $t\to\infty$.}
\]
Thus there exists $x_0=x_0(h)$ such
that $|\theta(t)|\le \gamma/2$ for all $t\ge x_0$. For all $x\ge x_0$ we get
\[
\frac{\varphi(x)}{\varphi(2x)}=2^{-\gamma}e^{\int_{x}^{2x}\frac{-\theta(t)}{t}dt}\le 2^{-\gamma}e^{\int_{x}^{2x}\frac{\gamma}{2t}dt}=2^{-\gamma}e^{\frac{\gamma}{2}\log 2}=2^{-\gamma/2}\eqqcolon d\in(0,1)\text{.}
\] 
We have $\varphi(x)\le d\varphi(2x)\iff\varphi(x)\le \frac{d}{1-d}(\varphi(2x)-\varphi(x))$, and thus
\[
\frac{1}{2^n}\sum_{k=0}^l \varphi(2^k)\lesssim O(2^{-n})+\frac{1}{2^n} \sum_{k=0}^l\big( \varphi(2^{k+1})- \varphi(2^{k})\big)=O(1)+\frac{\varphi(2^{l+1})}{2^n}\lesssim 1\text{, justifying the estimate in $\ref{geom}$.}
\]
Thus $\sup_{n\in\mathbb{N}_0}H_{2^n}|f|(x)\lesssim \sup_{n\in\mathbb{N}_0}S_{2^n}|f|(x)$, where $S_t|f|(x)=\frac{1}{t}\sum_{1\le m\le  t}|f(x-m)|$ and we get 
\[
\Big\|\Big(\sum_{j\in\mathbb{Z}}\big(\sup_{N\in\mathbb{N}}H_{2^N}|f_j|\big)^2\Big)^{1/2}\Big\|_{\ell^r(\mathbb{Z})}\lesssim \Big\|\Big(\sum_{j\in\mathbb{Z}}\big(\sup_{N\in\mathbb{N}_0}S_{2^N}|f_j|\big)^2\Big)^{1/2}\Big\|_{\ell^r(\mathbb{Z})}\lesssim_r \Big\|\big(\sum_{j\in\mathbb{Z}}|f_j|^2\big)^{1/2}\Big\|_{\ell^r(\mathbb{Z})}\text{,}
\]
since the vector-valued maximal estimates hold for the discrete Hardy--Littlewood maximal operator, see Theorem~1 in \cite{SteinFef} or Theorem~C  \cite{VectorValuedMax}. This justifies $\ref{finalstepinVV}$ and the proof of $\ref{VVMCT}$ is complete..

We now turn our attention to establishing the uniform $2$-oscillation estimates $\ref{U2OE}$. Again, by Calder\'on's transference principle it suffices to establish the following estimate
\[
\sup_{J\in\mathbb{N}}\sup_{I\in\mathfrak{S}_J(\mathbb{N})}\|O^2_{I,J}(A_n f:\,n\in\mathbb{N})\|_{\ell^r(\mathbb{Z})}\le C_{r,h}\|f\|_{\ell^r(\mathbb{Z})}\text{  for any $f\in \ell^r(\mathbb{Z})$}
\text{.}
\]
We begin by using Lemma~$\ref{correctosccomp}$ to reduce this task to the analogous one for the operator
\[
D_tf(x)=\frac{1}{\big(\sum_{p\le t}\log p\big)}\sum_{p\le t}\log(p) f(x-\lfloor h(p)\rfloor)\text{.}
\]
For convenience, let $\Psi(k)=\sum_{p\le k}\log(p)$ and $\pi(k)=|\mathbb{P}\cap[k]|$. Summation by parts yields 
\begin{multline*}
{A}_kf(x)=\frac{1}{\pi(k)}\sum_{p \le k}f(x-\lfloor h(p)\rfloor)=
\\
=\frac{\Psi(k)}{\pi(k)\log(k)}D_kf(x)-\sum_{2\le s\le k}\frac{\Psi(s)(\log^{-1}(s+1)-\log^{-1}(s))}{\pi(k)}
D_sf(x)=\sum_{s=2}^{\infty}\lambda_{s}^kD_sf(x)\text{,}
\end{multline*}
where
\[ 
\lambda_{s}^k= \left\{
\begin{array}{ll}
     \frac{\Psi(s)(\log^{-1}(s)-\log^{-1}(s+1))}{\pi(k)}&
     \text{if }2\le s\le k-1\\
      \frac{\Psi(k)}{\pi(k)\log(k)}&\text{if }s=k\\
      0&\text{if }s>k.\\
\end{array} 
\right. 
\]
Note that
\[
\sum_{s=2}^{\infty}\lambda_s^k=\sum_{s=2}^{k-1}\frac{\Psi(s)(\log^{-1}(s)-\log^{-1}(s+1))}{\pi(k)}+\frac{\Psi(k)}{\pi(k)\log(k)}=1\text{,}
\] 
and for any fixed $N\in\mathbb{N}_{\ge 2}$ the sequence $\sum_{s=2}^N\lambda_s^k$ is decreasing in $k$, since
\[
\sum_{s=2}^N\lambda_s^k=\left\{
\begin{array}{ll}
\frac{1}{\pi(k)}\sum_{s=2}^N \Psi(s)(\log^{-1}(s)-\log^{-1}(s+1))&\text{if }1\le N\le k-1\\
1&\text{if }N\ge k\\
\end{array}
\right.
\]
and for any $1\le N\le k-1$ we have by summation by parts
\[
\frac{1}{\pi(k)}\sum_{s=2}^N \Psi(s)(\log^{-1}(s)-\log^{-1}(s+1))=\frac{\pi(N+1)}{\pi(k)}-\frac{\Psi(N+1)}{\pi(N+1)\log(N+1)}\le 1\text{.}
\]
By Lemma~$\ref{correctosccomp}$ we get that
\[
\sup_{J\in\mathbb{N}}\sup_{I\in\mathfrak{S}_J(\mathbb{N})}\|O^2_{I,J}(A_n f:\,n\in\mathbb{N})\|_{\ell^r(\mathbb{Z})}\le \sup_{J\in\mathbb{N}}\sup_{I\in\mathfrak{S}_J(\mathbb{N})}\|O^2_{I,J}(D_n f:\,n\in\mathbb{N})\|_{\ell^r(\mathbb{Z})} \text{,}
\]
so we may focus on the new operator. 

We break 2-oscillations into short and long ones. Let $r_0\in(1,\infty)$ be such that $r\in(r_0,r_0')$, let $\tau\in(0,\min\big\{\frac{r_0-1}{2},\frac{1}{2}\big\})$ and let $\mathbb{D}_{\tau}=\{2^{n^{\tau}}:\,n\in\mathbb{N}_0\}$. We have that
 \begin{multline}\label{SL}
 \sup_{J\in\mathbb{N}}\sup_{I\in\mathfrak{S}_J(\mathbb{N})}\|O^2_{I,J}(D_{t}f:t\in\mathbb{N})\|_{\ell^r(\mathbb{Z})}\lesssim
 \sup_{J\in\mathbb{N}}\sup_{I\in\mathfrak{S}_J(\mathbb{D}_{\tau})}\|O^2_{I,J}(D_{t}f:t\in\mathbb{D}_{\tau})\|_{\ell^r(\mathbb{Z})}+
 \\
 +\bigg\|\bigg(\sum_{n=0}^{\infty}V^2\Big(D_nf:\,n\in\big[2^{n^{\tau}},2^{(n+1)^{\tau}}\big)\Big)^2\bigg)^{1/2}\bigg\|_{\ell^r(\mathbb{Z})}
 \end{multline}
 where 
 \[
V^2\Big(D_tf(x):\,t\in\big[2^{n^{\tau}},2^{(n+1)^{\tau}}\big)\Big)=\sup_{J\in\mathbb{N}}\sup_{\substack{t_0<\dots<t_J\\t_j\in\big[2^{n^{\tau}},2^{(n+1)^{\tau}}\big)}}\Big(\sum_{j=0}^{J-1}|D_{t_{j+1}}f(x)-D_{t_j}f(x)|^2 \Big)^{1/2}\text{, see page 17 in \cite{MOE}.}
\]
We begin with the short oscillations. One useful observation is that for large $n$ we have
\begin{equation}\label{passtoprimes}
\sup_{\substack{t_0<\dots<t_J\\t_j\in\big[2^{n^{\tau}},2^{(n+1)^{\tau}}\big)}}\Big(\sum_{j=0}^{J-1}|D_{t_{j+1}}f(x)-D_{t_j}f(x)|^2 \Big)^{1/2}\lesssim\sup_{\substack{t_0<\dots<t_J\\t_j\in\big[2^{{(n-1)}^{\tau}},2^{(n+2)^{\tau}}\big)\cap\mathbb{P}}}\Big(\sum_{j=0}^{J-1}|D_{t_{j+1}}f(x)-D_{t_j}f(x)|^2 \Big)^{1/2}\text{.}
\end{equation}
To see why this is true firstly note that $D_t=D_{t'}$ whenever $\mathbb{P}\cap[t]=\mathbb{P}\cap[t']$. Also, for any $A>0$, we have that there exists $C=C(A)$ such that for any $n\ge C$, we get
\begin{equation}\label{primesinpsdyad}
|\mathbb{P}\cap\big[2^{n^{\tau}},2^{(n+1)^{\tau}}\big)|=\int_{2^{n^{\tau}}}^{2^{(n+1)^\tau}}\frac{dt}{\log t}+O_A\bigg(\frac{2^{(n+1)^{\tau}}}{\log^A\big(2^{(n+1)^\tau}\big)}\bigg)\text{,}
\end{equation}
where we have used Corollary 5.26 in \cite{IWKO}. The error term can be simplified to $O_A\big(n^{-A\tau}2^{n^{\tau}}\big)$,
and it is not difficult to see using the Mean Value Theorem that
\begin{equation}\label{countprimes}
\int_{2^{n^{\tau}}}^{2^{(n+1)^\tau}}\frac{dt}{\log t}\simeq_{\tau}n^{-\tau}\big(2^{(n+1)^\tau}-2^{n^\tau}\big)\simeq_{\tau} \frac{2^{n^{\tau}}}{n}\text{.}
\end{equation}
Returning to $\ref{primesinpsdyad}$, we see that our estimate applied for $A=2\tau^{-1}$ implies that 
\begin{equation}\label{upperboundforshort}
|\mathbb{P}\cap\big[2^{n^{\tau}},2^{(n+1)^{\tau}}\big)|\simeq_{\tau}n^{-1}2^{n^{\tau}} \text{.}
\end{equation}
In particular this means that for $n\gtrsim_{\tau}1$ we get that $\mathbb{P}\cap\big[2^{n^{\tau}},2^{(n+1)^{\tau}}\big)\neq \emptyset$ and the reasoning for $\ref{passtoprimes}$ is justified. Now we see that
\[
\bigg\|\bigg(\sum_{n=0}^{\infty}V^2\Big(D_tf:\,t\in\big[2^{n^{\tau}},2^{(n+1)^{\tau}}\big)\Big)^2\bigg)^{1/2}\bigg\|_{\ell^r(\mathbb{Z})}
\]
\[
\lesssim\bigg\|\bigg(\sum_{n=1}^{\infty}\bigg(\sup_{J\in\mathbb{N}}\sup_{\substack{t_0<\dots< t_J\\t_j\in\big[2^{(n-1)^{\tau}},2^{(n+2)^{\tau}}\big)\cap \mathbb{P}}}\Big(\sum_{j=0}^{J-1}|D_{t_{j+1}}f-D_{t_{j}}f|\Big)\bigg)^{2}\bigg)^{1/2}\bigg\|_{\ell^r(\mathbb{Z})}\text{,}
\] 
where we also used that $\|\cdot\|_{\ell^2}\le\|\cdot\|_{\ell^1}$. For any $n
\in\mathbb{N}$, let $\big\{p_{m}^{(n)}:\,m\in\{0,\dotsc,l_n\}\big\}$ be an increasing enumeration of $\big[2^{(n-1)^{\tau}},2^{(n+2)^{\tau}}\big)\cap \mathbb{P}$. We may use the triangle inequality to bound the last expression by 
\[
\bigg\|\bigg(\sum_{n=1}^{\infty}\bigg(\Big(\sum_{m=1}^{l_n}|D_{p_{m}^{(n)}}f-D_{p_{m-1}^{(n)}}f|\Big)\bigg)^{2}\bigg)^{1/2}\bigg\|_{\ell^r(\mathbb{Z})}\text{.}
\]
Let $K_t(x)=\frac{1}{\Psi(t)}\sum_{p\le t}\log(p)\delta_{\lfloor h(p)\rfloor}(x)$, where $\delta_n(x)=1_{\{n\}}(x)$, and note that the expression above becomes
\[
\bigg\|\bigg(\sum_{n=1}^{\infty}\bigg(\Big(\sum_{m=1}^{l_n}|(K_{p_{m}^{(n)}}-K_{p_{m-1}^{(n)}})*f|\Big)\bigg)^{2}\bigg)^{1/2}\bigg\|_{\ell^r(\mathbb{Z})}\text{.}
\] 
With an argument identical to the one in page 15 in \cite{WT11LD} or in page 22 in \cite{shortvar}, it is easy to see that
\begin{equation}\label{SFEst}
\bigg\|\bigg(\sum_{n=0}^{\infty}V^2\Big(D_tf:\,t\in\big[2^{n^{\tau}},2^{(n+1)^{\tau}}\big)\Big)^2\bigg)^{1/2}\bigg\|_{\ell^r(\mathbb{Z})}\lesssim \bigg(\sum_{n=1}^{\infty}\bigg(\sum_{m=1}^{l_n}\big\|K_{p_{m}^{(n)}}-K_{p_{m-1}^{(n)}}\big\|_{\ell^1(\mathbb{Z})}\bigg)^{q}\bigg)^{1/q}\|f\|_{\ell^r(\mathbb{Z})}\text{,}
\end{equation}
where $q=\min\{2,r\}$. Note that
\[
\big\|K_{p_{m}^{(n)}}-K_{p_{m-1}^{(n)}}\big\|_{\ell^1(\mathbb{Z})}=\sum_{p\le p^{(n)}_{m-1}}\log(p)\bigg(\frac{1}{\Psi(p^{(n)}_{m-1})}-\frac{1}{\Psi(p^{(n)}_{m})}\bigg)+\frac{\log(p_m^{(n)})}{\Psi(p_m^{(n)})}
\]
\[
=1-\frac{\Psi(p_{m-1}^{(n)})}{\Psi(p_m^{(n)})}+\frac{\log(p_m^{(n)})}{\Psi(p_m^{(n)})}=\bigg(1-\frac{\Psi(p_{m-1}^{(n)})+\log(p_m^{(n)})}{\Psi(p_m^{(n)})}\bigg)+\frac{2\log(p_m^{(n)})}{\Psi(p_m^{(n)})}=\frac{2\log(p_m^{(n)})}{\Psi(p_m^{(n)})}\text{.}
\]
Thus
\begin{equation}\label{shortalm}
\sum_{m=1}^{l_n}\big\|K_{p_{m}^{(n)}}-K_{p_{m-1}^{(n)}}\big\|_{\ell^1(\mathbb{Z})}\lesssim\frac{l_n\log(2^{(n+2)^\tau})}{\Psi(p_{0}^{(n)})}\lesssim\frac{\big|\mathbb{P}\cap\big[2^{(n-1)^{\tau}},2^{(n+2)^{\tau}}\big)\big|n^{\tau}}{\Psi(2^{(n-1)^\tau})}\lesssim \frac{n^{\tau-1}2^{n^{\tau}}}{2^{(n-1)^\tau}}\lesssim n^{\tau-1}\text{,}
\end{equation}
where we used that
$\Psi(x)\simeq x$ as well as that $\big|\mathbb{P}\cap\big[2^{(n-1)^{\tau}},2^{(n+2)^{\tau}}\big)\big|\lesssim n^{-1}2^{n^{\tau}}$, see $\ref{upperboundforshort}$.

Using $\ref{shortalm}$ we can conclude by noting that
\[
\bigg(\sum_{n=1}^{\infty}\bigg(\sum_{m=1}^{l_n}\big\|K_{p_{m}^{(n)}}-K_{p_{m-1}^{(n)}}\big\|_{\ell^1(\mathbb{Z})}\bigg)^{q}\bigg)^{1/q}\lesssim\bigg(\sum_{n=1}^{\infty}n^{q(\tau-1)}\bigg)^{1/q}\text{.}
\]
Note that if $r>2$ then $q=2$ and then  $q(1-\tau)>2(1-1/2)=1$ thus the sum converges. Similarly, if $r\le 2$, then $q=r$. Note that $\tau<(r_0-1)/2\le(r_0-1)/r_0=1/r_0'$, but then $q(1-\tau)>r(1-1/r_0')=r/r_0>1$, as desired. In either case, the series is summable and $\ref{SFEst}$ becomes
\[
\bigg\|\bigg(\sum_{n=0}^{\infty}V^2\Big(D_tf:\,t\in\big[2^{n^{\tau}},2^{(n+1)^{\tau}}\big)\Big)^2\bigg)^{1/2}\bigg\|_{\ell^r(\mathbb{Z})}\lesssim_{r,\tau}\|f\|_{\ell^r(\mathbb{Z})}\text{,}
\]
completing the bounds for the short oscillations.

It remains to bound the first term of $\ref{SL}$. Note that
\[ \sup_{J\in\mathbb{N}}\sup_{I\in\mathfrak{S}_J(\mathbb{D}_{\tau})}\|O^2_{I,J}(D_{t}f:t\in\mathbb{D}_{\tau})\|_{\ell^r(\mathbb{Z})}=
\sup_{J\in\mathbb{N}_0}\sup_{I\in \mathfrak{S}_J(\mathbb{N}_0)}\|O^2_{I,J}\big(D_{2^{n^{\tau}}}f:\,n\in\mathbb{N}_0\big)\|_{\ell^r(\mathbb{Z})}\text{,}
\]
and we now establish the following estimate.
\begin{claim} There exists a constant $C_{r,\tau}>0$ such that
\[
\sup_{J\in\mathbb{N}}\sup_{I\in \mathfrak{S}_J(\mathbb{N}_0)}\|O^2_{I,J}\big(D_{2^{n^{\tau}}}f:\,n\in\mathbb{N}_0\big)\|_{\ell^r(\mathbb{Z})}
\le C_{r,\tau}\|f\|_{\ell^r(\mathbb{Z})}\text{.}
\] 
\end{claim}\,\\
\textit{Proof of the Claim.} Let
\[
\widetilde{H}_tf(x)=\frac{1}{\big(\sum_{n\le h(t)}\varphi'(n)\big)}\sum_{n\le h(t)}\varphi'(n)f(x-n)\text{.}\]
We have that $\|D_t\|_{\ell^r(\mathbb{Z})\to\ell^r(\mathbb{Z})},\|M_t\|_{\ell^r(\mathbb{Z})\to\ell^r(\mathbb{Z})},\|H_t\|_{\ell^r(\mathbb{Z})\to\ell^r(\mathbb{Z})},\|\widetilde{H}_t\|_{\ell^r(\mathbb{Z})\to\ell^r(\mathbb{Z})}<\infty$. We show that for any $\rho\in(0,2]$ we have
\[
\sum_{N\in\mathbb{N}_0}\|D_{2^{N^{\tau}}}-M_{2^{N^{\tau}}}\|^{\rho}_{\ell^2(\mathbb{Z})\to\ell^2(\mathbb{Z})},\sum_{N\in\mathbb{N}_0}\|M_{2^{N^{\tau}}}-H_{2^{N^{\tau}}}\|^{\rho}_{\ell^2(\mathbb{Z})\to\ell^2(\mathbb{Z})},\sum_{N\in\mathbb{N}_0}\|H_{2^{N^{\tau}}}-\widetilde{H}_{2^{N^{\tau}}}\|^{\rho}_{\ell^2(\mathbb{Z})\to\ell^2(\mathbb{Z})}<\infty\text{.}
\] 
For the first estimate note that 
\[
\|D_{t}f-M_{t}f\|_{\ell^2(\mathbb{Z})}
\le \bigg|1-\frac{\sum_{p\le t}\log(p)}{t}\bigg|\|f\|_{\ell^2(\mathbb{Z})}\text{,}
\]
and applying Corollary~5.29 with $A=\frac{2}{\tau \rho}$ yields 
\[
\|D_{2^{N^{\tau}}}f-M_{2^{N^{\tau}}}f\|_{\ell^2(\mathbb{Z})}=\bigg|1-\frac{2^{N^{\tau}}+O_{A}(2^{N^{\tau}}\log^{-A}(2^{N^{\tau}}))}{2^{N^{\tau}}}\bigg|\|f\|_{\ell^2(\mathbb{Z})}\lesssim_AN^{-A\tau}\|f\|_{\ell^2(\mathbb{Z})}\text{,}
\]
and thus
\[`
\sum_{N\in\mathbb{N}_0}\|D_{2^{N^{\tau}}}-M_{2^{N^{\tau}}}\|^{\rho}_{\ell^2(\mathbb{Z})\to\ell^2(\mathbb{Z})}\lesssim_{\rho,\tau}\sum_{N\in\mathbb{N}}N^{-A\tau\rho}<\infty\text{.}
\]
For the second estimate, we have already established that $\|M_{t}f-H_{t}f\|_{\ell^2(\mathbb{Z})}\lesssim e^{-(\log t)^{1/4}}
\|f\|_{\ell^2(\mathbb{Z})}$, and thus
\[
\sum_{N\in\mathbb{N}_0}\|M_{2^{N^{\tau}}}-H_{2^{N^{\tau}}}\|^{\rho}_{\ell^2(\mathbb{Z})\to\ell^2(\mathbb{Z})}\lesssim \sum_{N\in\mathbb{N}_0}\Big(e^{-(\log 2^{N^\tau})^{1/4}}\Big)^{\rho}=\sum_{N\in\mathbb{N}_0} e^{-\rho (\log(2))^{1/4}N^{\tau/4}}<\infty\text{.}
\]
For the third estimate note that 
\[
\|H_{t}f-\widetilde{H}_{t}f\|_{\ell^2(\mathbb{Z})}\lesssim \bigg|\frac{\sum_{n\le h(t)}\varphi'(n)}{t}-1\bigg|\|f\|_{\ell^2(\mathbb{Z})}\lesssim \frac{1}{t}\|f\|_{\ell^2(\mathbb{Z})}\text{,}
\]
since from the basic properties of $\varphi$ one gets $\sum_{n\le h(t)}\varphi'(n)=t+O(1)$, see for example Step~4 of Lemma~$\ref{majarckey}$. Thus  
\[
\sum_{N\in\mathbb{N}_0}\|H_{2^{N^{\tau}}}-\widetilde H_{2^{N^{\tau}}}\|^{\rho}_{\ell^2(\mathbb{Z})\to\ell^2(\mathbb{Z})}\lesssim \sum_{N\in\mathbb{N}_0}2^{-\rho N^{\tau}}<\infty\text{.}
\]
Using Lemma~$\ref{Comparison}$ thrice yields
\[
\sup_{J\in\mathbb{N}}\sup_{I\in \mathfrak{S}_J(\mathbb{N}_0)}\|O^2_{I,J}\big(D_{2^{n^{\tau}}}f:\,n\in\mathbb{N}_0\big)\|_{\ell^r(\mathbb{Z})}
\lesssim \sup_{J\in\mathbb{N}}\sup_{I\in \mathfrak{S}_J(\mathbb{N}_0)}\|O^2_{I,J}\big(\widetilde{H}_{2^{n^{\tau}}}f:\,n\in\mathbb{N}_0\big)\|_{\ell^r(\mathbb{Z})}+\|f\|_{\ell^r(\mathbb{Z})}\text{.}
\]
For the operator $\widetilde{H}_t$ full $2$-oscillations can be established straightforwardly using our comparison Lemma~$\ref{correctosccomp}$ and comparing with the standard averaging operator $S_tf(x)=\frac{1}{t}\sum_{m\le t}f(x-m)$, in fact this has already been done in \cite{WT11LD}, see page 21-23. We get that
\[
\sup_{J\in\mathbb{N}}\sup_{I\in\mathfrak{S}_J(\mathbb{N})}\|O^2_{I,J}(\widetilde{H}_n f:\,n\in\mathbb{N})\|_{\ell^r(\mathbb{Z})}\le \sup_{J\in\mathbb{N}}\sup_{I\in\mathfrak{S}_J(\mathbb{N})}\|O^2_{I,J}(S_n f:\,n\in\mathbb{N})\|_{\ell^r(\mathbb{Z})}\lesssim \|f\|_{\ell^r(\mathbb{Z})} \text{,}
\]
since the $2$-oscillation estimates hold for the standard averaging operator, see for example Theorem~1.4 in \cite{UOIAO}. The proof is complete.
\end{proof}
\section{Exponential sum estimates}\label{ExpEstSec}
This section is devoted to the proof of Proposition~$\ref{MainExpSumEst}$. We use the basic properties of the functions $h\in\mathcal{R}_c$ and their compositional inverses $\varphi$, specifically Lemma~2.1, 2.6 and 2.14 from \cite{MMR} throughout this section. We break the analysis into major and minor arc regimes in the following two lemmas.

\begin{lemma}[Minor arc estimate]\label{Minor}
Assume $c\in (1,4/3)$, $ h\in\mathcal{R}_c$ and $\theta_1\in\big(0,\frac{6c-2}{9}\big)$. Then there exists positive constants $\chi=\chi(h,\theta_1)$ and $C=C(h,\theta_1)$, such that for all $N\in\mathbb{N}$ and $\xi\in\mathbb{T}$ with $\|\xi\|>N^{-\theta_1}$, we have
\[
\Big|\sum_{p\le N}\log(p)e(\lfloor h(p)\rfloor\xi)\Big|\le C N^{1-\chi} 
\quad
\text{and}
\quad
\Big|\sum_{n\le h(N)}\varphi'(n)e(n\xi)\Big|\le CN^{1-\chi}\text{.} 
\]
Moreover, any $\chi\in\Big(0,\min\Big\{\frac{8-6c}{45},\frac{6c-2-9\theta_1}{36}\Big\}\Big)$ is admissible and thus for the choice $\theta_1=\frac{6c}{5}-\frac{14}{15}<\frac{6c-2}{9}$ we get that any $\chi\in\big(0,\frac{8-6c}{45}\big)$ is admissible.
\end{lemma}

\begin{lemma}[Major arc estimate]\label{Major}
Assume $c\in (1,2)$, $ h\in\mathcal{R}_c$, $\theta_1\in(c-1,1)$ and $\varepsilon>0$. Then there exists a positive constant $C=C(\varepsilon,h,\theta_1)$ such that for all $N\in\mathbb{N}$ and $\xi\in\mathbb{T}$ with $\|\xi\|\le N^{-\theta_1}$, we have
\begin{equation}\label{mainmajarc}
\Big|\sum_{p\le N}\log(p)e(\lfloor h(p)\rfloor\xi)-\sum_{n\le h(N)}\varphi'(n)e(n\xi)\Big|\le CNe^{-(\log N)^{1/3-\varepsilon}}\text{.}
\end{equation}
\end{lemma}
Notice that we can apply both lemmas simultaneously for $c\in(1,4/3)$ by choosing $\theta_1=\frac{6c}{5}-\frac{14}{15}>c-1$.
\subsection{Minor arc estimates}
\begin{proof}[Proof of Lemma~\ref{Minor}]Fix $c\in(1,4/3)$, $h\in\mathcal{R}_c$, $\theta_1\in(0,(6c-2)/9)$, $\chi<\min\Big\{\frac{8-6c}{45},\frac{6c-2-9\theta_1}{36}\Big\}$ and we note that all implicit constants appearing may depend on them. We also identify $\mathbb{T}$ with  $[-1/2,1/2)$ and then we get that $\|\xi\|=|\xi|$ for all $\xi\in\mathbb{T}$.

We begin with the rather straightforward second assertion of our lemma. Specifically, we show that for all $N\in\mathbb{N}$ and $\xi\in[-1/2,N^{-\theta_1})\cup (N^{-\theta_1},1/2)$, we get 
\begin{equation}\label{minorarcapprox}
\bigg|\sum_{n\le h(N)}\varphi'(n)e(n\xi)\bigg|\lesssim N^{\theta_1}\text{.}
\end{equation}
To see this, note that for all $L\in\mathbb{N}$, summation by part yields 
\[
\bigg|\sum_{1\le n\le L}\varphi'(n)e(n\xi)\bigg|=\bigg|\Big(\sum_{1\le n\le L}e(n\xi)\Big)\varphi'(L)-\sum_{1\le n\le L-1}\Big(\sum_{1\le m\le n}e(m\xi)\Big)\big(\varphi'(n+1)-\varphi'(n)\big)\bigg|\text{.}
\]
We use the fact that $\Big|\sum_{1\le n\le L}e(n\xi)\Big|\lesssim \|\xi\|^{-1}$ to bound the previous expression by a constant multiple of
\[
\|\xi\|^{-1}\varphi'(L)+\|\xi\|^{-1}\sum_{1\le n\le N-1}|\varphi'(n)-\varphi'(n+1)|\lesssim\|\xi\|^{-1}
\]
since the sum eventually telescopes, $\varphi'$ is eventually decreasing and $\varphi'(x)\to 0$ as $x\to\infty$. For all $\xi$ with $\|\xi\|>N^{-\theta_1}$, we have 
\begin{equation}\label{firstres}
\bigg|\sum_{n\in[1,\lfloor h(N)\rfloor ]}\varphi'(n)e(n\xi)\bigg|\lesssim \|\xi\|^{-1}\lesssim N^{\theta_1}\text{,}
\end{equation}
and the proof of the estimate $\ref{minorarcapprox}$ is complete since $\theta_1<1-\chi$, because $
1-\theta_1>1-\frac{6c-2}{9}=\frac{11-6c}{9}>\frac{8-6c}{45}>\chi$, where one can easily see that the second to last inequality holds because it is equivalent to $c<47/24$. We now focus on the first estimate of Lemma~$\ref{Minor}$, namely we establish the estimate
\begin{equation}\label{MminAest}
\bigg|\sum_{p\le N}\log(p)e(\lfloor h(p)\xi\rfloor)\bigg|\lesssim N^{1-\chi}\text{.}
\end{equation}
We can replace $\log$ with $\Lambda$ at the expense of an acceptable error. More precisely, we have
\begin{equation}\label{ReplwL}
\sum_{p\le N}\log(p)e(\lfloor h(p)\rfloor\xi)=\sum_{n\le N}\Lambda(n)e(\lfloor h(n)\rfloor\xi)+O(N^{1/2})\text{.}
\end{equation}
To see this, note that
\[
\Big|\sum_{p\le N}\log(p)e(\lfloor h(p)\rfloor\xi)-\sum_{n\le N}\Lambda(n)e( \lfloor h(n)\rfloor\xi)\Big|\le \sum_{\substack{1\le p^s\le N,\\p\in\mathbb{P},s\ge2}}\log p\text{.}
\]
 Each prime $p$ will contribute $\log p $ to the sum exactly $s_p-1$ times where $s_p=s_p(N)$ is the unique integer satisfying $p^{s_p}\le N< p^{s_p+1}$ or equivalently $s_p=\lfloor \log(N)/\log(p)\rfloor$. Thus
\begin{equation}\label{Ltolog}
\sum_{\substack{1\le p^s\le N,\\p\in\mathbb{P},s\ge2}}\log p \le \sum_{ p\le N^{1/2}}\bigg\lfloor \frac{\log N}{\log p}\bigg\rfloor \log p \le \log(N)\sum_{p\le N^{1/2}}1\lesssim \frac{\log(N)N^{1/2}}{\log(N^{1/2})}\lesssim N^{1/2}\text{,}
 \end{equation}
where we have used the fact that $|\mathbb{P}\cap [1,x]|\lesssim x(\log x)^{-1}$. Let us note that for our choice of $\chi$ we do have that $1/2<1-\chi$, since $\chi<\frac{8-6c}{45}<\frac{8}{45}<\frac{1}{2}$. This establishes $\ref{ReplwL}$ and now we may focus on the approximant.

We have that
\[
\sum_{n\le N}\Lambda(n)e(\lfloor h(n)\rfloor\xi)=\sum_{n\le N}\Lambda(n)e(h(n)\xi)e(-\{h(n)\}\xi)\text{,}
\]
and we now use a standard Fourier series expansion. For every $M\ge 2$, $x\in[-1/2,1/2)\setminus\{0\}$ and $y\in\mathbb{R}$ we have that
\[
e(-x\{y\})=\sum_{|m|\le M}\frac{1-e(-x)}{2\pi i (m+x)}e(my)+O\bigg(\min\bigg(1,\frac{1}{M\|y\|}\bigg)\bigg)\text{,}
\]
where the implied constant is absolute, see page 29 in \cite{LPE}. We apply this for $M=\lfloor N^{\theta_2}\rfloor$, for $\theta_2=(8-6c)/45>0$. Thus we get
\begin{equation}\label{1stapprox}
\sum_{n\le N}\Lambda(n)e(\lfloor h(n)\rfloor\xi)=\sum_{n\le N}\Lambda(n)e(h(n)\xi)\sum_{|m|\le M}\frac{1-e(-\xi)}{2\pi i (m+\xi)}e(mh(n))+\sum_{n\le N}\Lambda(n)e(h(n)\xi)g_M(\xi,n)
\end{equation}
where $|g_M(\xi,n)|\lesssim \min\Big(1,\frac{1}{M\|h(n)\|}\Big)$.

We bound the second term. We have that 
\begin{equation}\label{ESTerrmin}
\Big|\sum_{n\le N}\Lambda(n)e(h(n)\xi)g_M(\xi,n)\Big|\lesssim  \log(N)\sum_{n\le N} \min\bigg(1,\frac{1}{M\|h(n)\|}\bigg)\text{.}
\end{equation}
We use the following lemma to conclude. It is an adaptation of Lemma~4.4 from \cite{LPE}.
\begin{lemma}\label{Easyerror}
For every $c\in(1,2)$ and $h\in\mathcal{R}_c$, there exists a positive real number $C_h$ such that for all $N,M\in\mathbb{N}_{\ge 2}$ we have
\begin{equation}\label{Lemerroroferror}
\sum_{n\le N} \min\bigg(1,\frac{1}{M\|h(n)\|}\bigg)\le C_h \log(N)\log(M)\big(NM^{-1}+h(N)^{1/2}M^{1/2}\big)\text{.}
\end{equation}
\end{lemma}
\begin{proof}
We start by estimating a dyadic piece. Let $2\leq P\leq P_1\leq2P$. We now use another standard Fourier series expansion, namely, for every $M\in\mathbb{N}$ and $x\in\mathbb{R}$ we have that
\[
\min\bigg(1,\frac{1}{M\|x\|}\bigg)=\sum_{m\in\mathbb{Z}}a_me(mx)
\]
where $|a_m|\lesssim \min\Big(\frac{\log M}{M},\frac{M}{|m|^2}\Big)$, see page 23 in \cite{LPE}. We note that $a_m$ also depends on $M$, but we decide to suppress that dependence to make the exposition more reasonable. We have that
\[
\sum_{P\leq n\leq P'}\min\left(1,\frac{1}{M\|h(n)\|}\right)=\sum_{P\leq n\leq P'}\sum_{m\in\mathbb{Z}}a_me(mh(n))\le \sum_{m\in\mathbb{Z}}|a_m|\Big|\sum_{P\leq n\leq P'}e(mh(n))\Big|\text{.}
\]
We now estimate $\Big|\sum_{P\leq n\leq P'}e(mh(n))\Big|$ using van der Corput's lemma. Note that for all $m\neq 0$ and $F(x)=mh(x)$, we have that for all $x\in [P,P_1]$ 
\[
|F''(x)|=|m|h''(x)\simeq |m|\frac{h(x)}{x^2}\simeq |m|\frac{h(P)}{P^2}
\]
for sufficiently large $P$, where we have used the basic properties of $h$. More precisely, according to Lemma~2.1 in \cite{MMR}, we have that there exist real-valued functions $\vartheta_1,\vartheta_2$ with
\[
x^2h''(x)=xh'(x)(c-1+\vartheta_2(x))=h(x)(c-1+\vartheta_2(x))(c+\vartheta_1(x))
\]
and $\lim_{x\to\infty}\vartheta_i(x)=0$ for $i=1,2$. This means that there exist a constant $x_0$ and constants $0<\tilde{c}<\widetilde{C}$ such that $\tilde{c}h(x)\le x^2h''(x)\le \widetilde{C}h(x)$ for all $x\ge x_0$. It is not difficult to see that $h(x)\le h(2x)\lesssim h(x)$ (the first inequality is immediate and the second estimate is easily derivable from the fact that $\vartheta$ in the definition of $h$ is bounded). Applying van der Corput's lemma, see for example Lemma~3.1 in \cite{MMR}, gives
\[
 \Big|\sum_{P\leq n\leq P'}e(mh(n))\Big|\lesssim  P(|m|^{1/2}P^{-1}h(P)^{1/2})+|m|^{-1/2}Ph(P)^{-1/2}
\]
\[
\lesssim|m|^{1/2}h(P)^{1/2}+|m|^{-1/2}Ph(P)^{-1/2}\lesssim |m|^{1/2}h(P)^{1/2}\text{,}
\]
where for the last estimate we used the fact that $P\lesssim h(P)$ (see Lemma~2.6 in \cite{MMR}). Thus we get
\[
\sum_{m\in\mathbb{Z}}|a_m|\Big|\sum_{P\leq n\leq P'}e(mh(n))\Big|=\sum_{|m|\leq M}|a_m|\Big|\sum_{P\leq n\leq P'}e(mh(n))\Big|+\sum_{|m|\geq M+1}|a_m|\Big|\sum_{P\leq n\leq P'}e(mh(n))\Big|
        \]
        \[
        \lesssim\frac{\log M}{M}\Big(P+\sum_{1\leq|m|\leq M}|m|^{1/2}h(P)^{1/2}\Big)+\sum_{|m|\geq M+1}\frac{M}{|m|^2}|m|^{1/2}h(P)^{1/2}
\]
\[
\lesssim\frac{\log M}{M}\Big(P+h(P)^{1/2}M^{3/2}\Big)+h(P)^{1/2}M^{1/2}
\]
\[
\lesssim\log(M)PM^{-1}+\log(M)h(P)^{1/2}M^{1/2}+h(P)^{1/2}M^{1/2}\lesssim \log(M)\big(PM^{-1}+h(P)^{1/2}M^{1/2}\big)\text{.}        
\]
Notice that the bound is increasing in $P$. Thus for all $M\ge 2$ we have that
\[
\sum_{n\le N} \min\bigg(1,\frac{1}{M\|h(n)\|}\bigg)\lesssim \log(N)\sup_{1\le P\le P'\le 2P\le N}\sum_{P\leq n\leq P'}\min\left(1,\frac{1}{M\|h(n)\|}\right)
\]
\[
\lesssim\log(N)\log(M)\big(NM^{-1}+h(N)^{1/2}M^{1/2}\big)\text{,}\quad\text{and the proof of the lemma is complete.}
\]
\end{proof}
Combining the estimate \ref{Lemerroroferror} of our lemma with $\ref{ESTerrmin}$ we obtain
\begin{equation}\label{finishedestforgm}
\Big|\sum_{n\le N}\Lambda(n)e(h(n)\xi)g_M(\xi,n)\Big|\lesssim \log^2(N)\log(M)\big(NM^{-1}+h(N)^{1/2}M^{1/2}\big)\text{.}
\end{equation}
Remembering that $\theta_2=(8-6c)/45$ and $M=\lfloor N^{\theta_2}\rfloor$ we get that
\begin{equation}\label{thirdres}
\Big|\sum_{n\le N}\Lambda(n)e(h(n)\xi)g_M(\xi,n)\Big|\lesssim \log^3(N)\big(N^{1-\theta_2}+h(N)^{1/2}N^{\theta_2/2}\big)\lesssim\log^3(N)\big(N^{1-\theta_2}+N^{c/2+1/100+\theta_2/2}\big)\text{,}
\end{equation}
where we have used that $h(x)\lesssim_{\varepsilon} x^{c+\varepsilon}$, see Lemma~2.6 in \cite{MMR}. We are done since $1-\theta_2<1-\chi$ and $c/2+1/100+\theta_2/2=c/2+(4-3c)/45+1/100<1-\chi$. Thus
\[
\Big|\sum_{n\le N}\Lambda(n)e(h(n)\xi)g_M(\xi,n)\Big|\lesssim N^{1-\chi}\text{,}
\]
and this concludes the estimation of the second summand in $\ref{1stapprox}$.

We now focus on the first term; it is not difficult to see that 
\[
\Big|\sum_{n\le N}\Lambda(n)e(h(n)\xi)\sum_{|m|\le M}\frac{1-e(-\xi)}{2\pi i (m+\xi)}e(h(n)m)\Big|=\Big|\sum_{|m|\le M}\frac{1-e(-\xi)}{2\pi i (m+\xi)}\sum_{n\le N}\Lambda(n)e(h(n)(\xi+m))\Big|
\]
\begin{equation}\label{LastMinorArc}
\lesssim\log(N)\bigg(\sup_{1\le P<P_1\le 2P\le N}\Big|\sum_{P<n\le P_1}\Lambda(n)e(h(n)\xi)\Big|+\sum_{1\le|m|\le M}\frac{1}{|m|}\sup_{1\le P<P_1\le 2P\le N}\Big|\sum_{P<n\le P_1}\Lambda(n)e(h(n)(\xi+m))\Big|\bigg)\text{.}
\end{equation}
Again, we focus on a dyadic piece and we begin by replacing the von Mangoldt function using Vaughan's identity. We state it here for the sake of clarity.
\begin{lemma}Let $v,w\ge 1$ be real numbers and let $n\in\mathbb{N}$ be such that $n>v$, then
 		\[
 		\Lambda(n)=\sum_{\substack{b|n\\b\le w}}
 		\mu(b)\log(n/b)-\mathop{\sum\sum}\limits_{\substack{bc|n\\b\le w,c\le v}}\mu(b)\Lambda(c)+\mathop{\sum\sum}\limits_{\substack{bc|n\\b>w,c>v}}\mu(b)\Lambda(c)
 		\] 
 		Or equivalently, for every $n>v$ we have
 		\[
 		\Lambda(n)=\sum_{kl=n\text{, }l\le w}\log(k)\mu(l)-\sum_{l\le vw}\sum_{kl=n}\Pi_{v,w}(l)+\sum_{kl=n\text{, }k>v\text{, }l> w}\Lambda(k)\Xi_w(l)
 		\]
 		where \[\Pi_{v,w}(l)=\sum_{\substack{rs=l\\r\le v\text{, }s\le w}}\Lambda(r)\mu(s)\quad\text{and}\quad\Xi_w(l)=\sum_{\substack{d|l\\d>w}}\mu(d)\text{.}
 		\]
 	\end{lemma}
\begin{proof}
The proof can be found in \cite{IWKO}, see section~13.4, page 344.
\end{proof}

Let $m\in\mathbb{Z}$ and $1\le P < P_1\le 2P\le N$. We apply Vaughan's identity for $v=w=P_1^{1/3}/2$ and we get that
\[
\sum_{P<n\le P_1}\Lambda(n)e((\xi+m)h(n))
\]
\[
=\bigg(\sum_{P<n\le P_1}\sum_{kl=n\text{, }l\le v}\log(k)\mu(l)e((\xi+m)h(n))-\sum_{P<n\le P_1}\sum_{l\le v^2}\sum_{kl=n}\Pi_{v,v}(l)e((\xi+m)h(n))
\]
\[
+\sum_{P<n\le P_1}\sum_{kl=n\text{, }k>v\text{, }l> v}\Lambda(k)\Xi_v(l)e((\xi+m)h(n))\bigg)
\]
\[
=\bigg(\sum_{l\le v}\sum_{P/l<k\le P_1/l}\log(k)\mu(l)e((\xi+m)h(kl))
\]
\[
-\sum_{l\le v}\sum_{P/l<k\le P_1/l}\Pi_{v,v}(l)e((\xi+m)h(kl))-\sum_{v<l\le v^2}\sum_{P/l<k\le P_1/l}\Pi_{v,v}(l)e((\xi+m)h(kl))
\]
\[
+\sum_{v<l\le P_1/v}\sum_{\substack{P/l<k\le P_1/l\\k>v}}\Lambda(k)\Xi_v(l)e((\xi+m)h(kl))\bigg)\eqqcolon S_1-S_{2,1}-S_{2,2}+S_3\text{.}
\]
The proof is now reduced to estimating these four terms. For $S_1$, we use summation by parts. Let us denote by $U_l(t)=\sum_{P/l<k\le t}e((\xi+m)h(kl))$, then
 	\[
 	|S_1|\le \sum_{l\le v} |\mu(l)| \bigg|\sum_{P/l<k\le P_1/l}\log(k)e((\xi+m)h(kl))\bigg|\le\sum_{l\le v}  \bigg|U_l(P_1/l)\log(P_1/l)-\int_{P/l}^{P_1/l}U_l(t)/tdt\bigg|
 	\]
 	\[
 	\le\sum_{l\le v}\Big(|U_l(P_1/l)|\log(P_1/l)+\sup_{P/l<t\le P_1/l}|U_l(t)|\big(\log(P_1/l)-\log(P/l)\big)\Big)\le 2\log(P_1)\sum_{l\le v}\sup_{P/l<t\le P_1/l}|U_l(t)|\text{.}
 	\] We estimate the dyadic pieces appearing using van der Corput. For $l\le v$, $x\in[P/l, P_1/l]$ and for $F(x)=(\xi+m)h(lx)$, let us notice that 
\[
|F''(x)|=|l^2(\xi+m)h''(lx)| \simeq |\xi+m| l^2\frac{h(xl)}{(xl)^2}\simeq|\xi+m|l^{2}h(P)P^{-2}\text{,}
\]
where we have used the basic properties of $h$ in an identical manner to the proof of Lemma~$\ref{Easyerror}$. Notice that $\xi+m$ can never be zero since $\xi\in[-1/2,N^{-\theta_1})\cup (N^{-\theta_1},1/2)$ and $m\in\mathbb{Z}$. By van der Corput's lemma we get that
\[
|U_l(t)|=\bigg|\sum_{k\in(P/l, t]}e(F(k))\bigg|\lesssim \big(P/l\big)|\xi+m|^{1/2}l\big(h(P)P^{-2}\big)^{1/2}+|m+\xi|^{-1/2}l^{-1}\big(h(P)P^{-2}\big)^{-1/2}
\]
\[
=|\xi+m|^{1/2}h(P)^{1/2}+|\xi+m|^{-1/2}l^{-1}h(P)^{-1/2}P
\text{,}
\]
and thus we get
\[
|S_1|\lesssim \log(P_1)\sum_{l\le v}\Big(|\xi+m|^{1/2}h(P)^{1/2}+|\xi+m|^{-1/2}l^{-1}h(P)^{-1/2}P\Big)
\]
\begin{equation}\label{S1}
\lesssim\log(P_1)|\xi+m|^{1/2}P^{1/3}h(P)^{1/2}+\log^2(P_1)|\xi+m|^{-1/2}Ph(P)^{-1/2}\text{,}
\end{equation}
where in the last estimate we took into account the fact that $v=P_1^{1/3}/2$. For $S_{2,1}$ the estimate follows from similar considerations. Firstly, notice that
 	\[
 	|S_{2,1}|=\bigg|\sum_{l\le v}\sum_{P/l<k\le P_1/l}\Pi_{v,v}(l)e((\xi+m)h(kl))\bigg|\le \sum_{l\le v}|\Pi_{v,v}(l)||U_l(P_1/l)| 
 	\] and also that $|\Pi_{v,v}(l)|=\bigg|\sum_{\substack{rs=l\\r\le v\text{, }s\le v}}\Lambda(r)\mu(s)\bigg|\le \sum_{r|l}\Lambda(r)=\log(l)\le \log(v)\le \log(P_1)\text{ for }l\le v$. Thus
 \[
 |S_{2,1}|\le\sum_{l\le v}|\Pi_{v,v}(l)||U_l(P_1/l)|\le  
 \log(P_1)\sum_{l\le v}\sup_{P/l<t\le P_1/l}|U_l(t)|
 \]
 and we can bound exactly as in the case of $S_1$, i.e.:
 \begin{equation}\label{S21}
|S_{2,1}| \lesssim \log(P_1)|\xi+m|^{1/2}P^{1/3}h(P)^{1/2}+\log^2(P_1)|\xi+m|^{-1/2}Ph(P)^{-1/2}\text{.}
\end{equation}

We treat $S_{2,2}$ and $S_3$ simultaneously. We have
 	\begin{equation}\label{DyadPc1} 		|S_{2,2}|\lesssim\log^2(P_1)\sup_{L\in[v,v^2]}\sup_{K\in[P/v^2,P_1/v]}\sup_{L'\in(L,2L]}\sup_{K'\in(K,2K]}\bigg|
 		\sum_{L< l\le L'}\sum_{\substack{K<k\le K'\\P<kl\le P_1}}\Pi_{v,v}(l)e((\xi+m)h(kl))\bigg|\text{,}
 	\end{equation}
 and
 	\begin{equation}\label{DyadPc2}
|S_3|\lesssim\log^2(P_1)\sup_{L\in[v,P_1/v]}\sup_{K\in[v,P_1/v]}\sup_{L'\in(L,2L]}\sup_{K'\in(K',2K]}\bigg|\sum_{L<l\le L'}\sum_{\substack{K<k\le K'\\P<kl\le P_1}}\Lambda(k)\Xi_v(l)e((\xi+m)h(kl))\bigg|\text{.}
 	\end{equation}
 	Note that
 	\begin{equation}\label{delta1}
 	\sum_{L<l\le L'}|\Pi_{v,v}(l)|^2 \lesssim L\log^2 L\text{,}\quad\sum_{L< l\le L'}|\Xi_v(l)|^2\lesssim  L\log^3 L\text{,}\quad\text{see for example page 31 in \cite{RHLMLD}.}
 	\end{equation}
We now use the following technical lemma. 
\begin{lemma}\label{techlem} Let $K,L,P,P_1\in\mathbb{N}$ be such that $ \frac{1}{10}K^{1/2}\le L\le 10 K^{2}$ and  $P<P_1\le 2P$, and let $m\in\mathbb{Z}$ and $\xi\in[-1/2,1/2)\setminus\{0\}$. For all $\varepsilon>0$ we have
 	\begin{multline*}
\bigg|\sum_{L<l\le L'\le 2L}\sum_{\substack{K<k\le K'\le 2K\\P<kl\le P_1}}\Delta_1(l)\Delta_2(k)e((\xi+m) h(kl))\bigg|
\\
\lesssim_{\varepsilon}\log^{4}(K)\Big((|m|+1)^{\frac{1}{4}}K^{\frac{5+c}{6}+\varepsilon}L^{\frac{4+c}{6}+\varepsilon}+|m+\xi|^{\frac{-1}{4}}K^{\frac{7-c}{6}+\varepsilon}L^{\frac{5-c}{6}+\varepsilon}\Big)\text{,}
\end{multline*}
for every two sequences of complex numbers $\big(\Delta_1(l)\big)_{l\in(L,2L]}$ and $\big(\Delta_2(k)\big)_{k\in(K,2K]}$ satisfying
 \[
 		\sum_{L<l\le 2L}|\Delta_1(l)|^2\lesssim L\log^3(L)\text{ and }\sum_{K<k\le 2K}|\Delta_2(k)|^2\lesssim K\log^3(K)\text{.}
\]
\end{lemma}
Before we establish Lemma~$\ref{techlem}$, let us see how it will help us estimate $|S_{2,2}|$ and $|S_3|$ and finally appropriately bound the second term of $\ref{1stapprox}$,  finishing the proof of Lemma~
$\ref{Minor}$. We return to $\ref{DyadPc1}$ and $\ref{DyadPc2}$. For any $K,L\in\mathbb{N}$ which make the dyadic piece nonempty, there exist $k,l\in\mathbb{N}$ such that
 	\[
 	KL< kl \le P_1\text{ and }KL\ge \frac{kl}{4}>P/4\ge P_1/8
\text{, and thus }P_1/8\le KL\le P_1\text{.} 
\]
For $S_{2,2}$, notice that
 	\[
 	K\le P_1/v=2P_1^{2/3}\text{ and }K\ge P/v^2\ge 2P_1^{1/3}\text{, and thus }K\in\big[P_1^{1/3}/2,2P_1^{2/3}\big] 
 	\]
 	and similarly
 	\[	L\in[v,v^2]=\big[P_1^{1/3}/2,P_1^{2/3}/4\big]\subseteq \big[P_1^{1/3}/2,2P_1^{2/3}\big] \text{.}
 	\]
For $S_3$, notice that
\[
K,L\in[v,P_1/v]=\big[P_1^{1/3}/2,2P_1^{2/3}\big]\text{.}
\]
Therefore, in either case, $K,L\in \big[P_1^{1/3}/2,2P_1^{2/3}\big] $. Thus it is easy to see that $\frac{1}{10}K^{1/2}\le L\le 10 K^2$ and since we have $\ref{delta1}$, the previous lemma is applicable. Taking into account the previous comments, the bound the lemma gives becomes 
\[
\log^{4}(K)\Big((|m|+1)^{\frac{1}{4}}K^{\frac{5+c}{6}+\varepsilon}L^{\frac{4+c}{6}+\varepsilon}+|m+\xi|^{\frac{-1}{4}}K^{\frac{7-c}{6}+\varepsilon}L^{\frac{5-c}{6}+\varepsilon}\Big)
\]
\[
\lesssim\log^4(P_1)\Big((|m|+1)^{\frac{1}{4}}(KL)^{\frac{5+c}{6}+\varepsilon}L^{\frac{-1}{6}}+|m+\xi|^{\frac{-1}{4}}(KL)^{\frac{7-c}{6}+\varepsilon}L^{\frac{-2}{6}}\Big)
\]
\[
\lesssim\log^4(P_1)\Big((|m|+1)^{\frac{1}{4}}P_1^{\frac{5+c}{6}+\varepsilon}P_1^{\frac{-1}{18}}+|m+\xi|^{\frac{-1}{4}}P_1^{\frac{7-c}{6}+\varepsilon}P_1^{\frac{-1}{9}}\Big)\text{.}
\] 	
Thus we get
\begin{equation}\label{S22nS3}
\max(|S_{2,2}|,|S_3|)\lesssim_{\varepsilon} \log^6(P_1)\Big((|m|+1)^{\frac{1}{4}}P_1^{\frac{5+c}{6}+\varepsilon}P_1^{\frac{-1}{18}}+|m+\xi|^{\frac{-1}{4}}P_1^{\frac{7-c}{6}+\varepsilon}P_1^{\frac{-1}{9}}\Big)  \text{.}
\end{equation}

Taking into account $\ref{S1}$, $\ref{S21}$,  $\ref{S22nS3}$, we can estimate the expression $\ref{LastMinorArc}$ by
\begin{multline*}
\log(N)\bigg(\sup_{1\le P<P_1\le 2P\le N}\Big|\sum_{P<n\le P_1}\Lambda(n)e(h(n)\xi)\Big|
\\
+\sum_{1\le|m|\le M}\frac{1}{|m|}\sup_{1\le P<P_1\le 2P\le N}\Big|\sum_{P<n\le P_1}\Lambda(n)e(h(n)(\xi+m))\Big|\bigg)
\end{multline*}
\begin{multline*}
\lesssim_{\varepsilon}\log^2(N)|\xi|^{1/2}N^{1/3}h(N)^{1/2}+\log^3(N)|\xi|^{-1/2}Nh(N)^{-1/2}
\\
+\log^7(N)\Big(N^{(5+c)/6+\varepsilon}N^{-1/18}+|\xi|^{-1/4}N^{(7-c)/6+\varepsilon}N^{-1/9}\Big) +\sum_{1\le|m|\le M}\frac{1}{|m|}\bigg(\log^2(N)|m|^{1/2}N^{1/3}h(N)^{1/2}
\\
+\log^3(N)|m|^{-1/2}Nh(N)^{-1/2}+\log^7(N)\Big(|m|^{1/4}N^{(5+c)/6+\varepsilon}N^{-1/18}+|m|^{-1/4}N^{(7-c)/6+\varepsilon}N^{-1/9}\Big)\bigg) 
\end{multline*}
\begin{multline*}
\lesssim\log^2(N)N^{1/3}h(N)^{1/2}+\log^3(N)N^{\theta_1/2}Nh(N)^{-1/2}
\\
+\log^7(N)\Big(N^{(5+c)/6+\varepsilon}N^{-1/18}+N^{\theta_1/4}N^{(7-c)/6+\varepsilon}N^{-1/9}\Big) +\log^2(N)N^{1/3}h(N)^{1/2}M^{1/2}
\\
+\log^3(N)Nh(N)^{-1/2}+\log^7(N)\Big(M^{1/4}N^{(5+c)/6+\varepsilon}N^{-1/18}+N^{(7-c)/6+\varepsilon}N^{-1/9}\Big)\bigg) 
\end{multline*}
\begin{multline*}
\lesssim_{\varepsilon}\log^2(N)N^{1/3}N^{c/2+\varepsilon}+\log^3(N)N^{\theta_1/2}NN^{-c/2+\varepsilon}
\\
+\log^7(N)\Big(N^{(5+c)/6+\varepsilon}N^{-1/18}+N^{\theta_1/4}N^{(7-c)/6+\varepsilon}N^{-1/9}\Big) +\log^2(N)N^{1/3}N^{c/2+\varepsilon}N^{\theta_2/2}
\\
+\log^3(N)NN^{-c/2+\varepsilon}+\log^7(N)\Big(N^{\theta_2/4}N^{(5+c)/6+\varepsilon}N^{-1/18}+N^{(7-c)/6+\varepsilon}N^{-1/9}\Big)\bigg)\text{.}
\end{multline*}
By carefully examining all eight summands we see that we obtain our desired estimate. For the convenience of the reader we provide the tedious details here. It suffices to have
\begin{enumerate}[label=\roman*)]
\item $1/3+c/2<1-\chi$ ($\iff \chi<\frac{4-3c}{6}$, which is true since $\chi<\frac{8-6c}{45}$),
\item $\theta_1/2+1-c/2<1-\chi$ ($\iff \chi<\frac{c-\theta_1}{2}$, which is true since $\frac{c-\theta_1}{2}>\frac{3c+2}{18}>\frac{5}{18}>\frac{8-6c}{45}>\chi$),
\item $5/6+c/6-1/18<1-\chi$ ($\iff \chi<\frac{4-3c}{18}$, which is true),
\item $\theta_1/4+7/6-c/6-1/9<1-\chi$ ($\iff \chi<\frac{6c-2-9\theta_1}{36}$, which is true),
\item $1/3+c/2+\theta_2/2<1-\chi$ ($\iff \chi<\frac{4-3c}{6}-\frac{4-3c}{45}$, which is true since $\chi<\frac{8-6c}{45}$),
\item $1-c/2<1-\chi$ ($\iff \chi<c/2$, which is true since $\chi<1/2<c/2$),
\item $\theta_2/4+5/6+c/6-1/18<1-\chi$ ($\iff \chi<\frac{8-6c}{45}$, which is true),
\item $7/6-c/6-1/9<1-\chi$ ($\iff \chi<\frac{3c-1}{18}$, which is true since $\frac{3c-1}{18}>\frac{2}{18}>\frac{8-6c}{45}$),
\end{enumerate}
which we verified. Since all these exponents are smaller than $1-\chi$, we can fix $\varepsilon>0$ appropriately small, to obtain
\begin{multline*}
\log(N)\bigg(\sup_{1\le P<P_1\le 2P\le N}\Big|\sum_{P<n\le P_1}\Lambda(n)e(h(n)\xi)\Big|
\\
+\sum_{1\le|m|\le M}\frac{1}{|m|}\sup_{1\le P<P_1\le 2P\le N}\Big|\sum_{P<n\le P_1}\Lambda(n)e(h(n)(\xi+m))\Big|\bigg)\lesssim N^{1-\chi}\text{.}
\end{multline*}
We have estimated both summands from $\ref{1stapprox}$, and the proof of Lemma~$\ref{Minor}$ is complete.
\end{proof} 

\begin{proof}[Proof of Lemma~$\ref{techlem}$]
We follow the ideas from of Lemma~2.16 from \cite{MMR}, see also Section~4 from \cite{ExpSumMethods}. We set up some notation first. For every $r\in\mathbb{Z}$ let 
\begin{equation}\label{DefEr}
E_r=\sum_{L<l\le 2L}\sum_{\substack{K<k,k+r\le K'\\P<kl,(k+r)l\le P_1}}\Delta_2(k)\overline{\Delta_2(k+r)}e\big((\xi+m)h(kl)-(\xi+m)h((k+r)l)\big)\text{,}
\end{equation}
and note that 
\begin{equation}\label{Ezerobound}
|E_0|\le \sum_{L<l\le 2L}\sum_{K<k\le K'}|\Delta_2(k)|^2\lesssim L K\log^3(K)\text{.} 
\end{equation}
Using the Cauchy-Schwarz inequality together with the assumption for $\big(\Delta_1(l)\big)_{l\in(L,2L]}$ we get
\[
\bigg|\sum_{L<l\le L'\le 2L}\sum_{\substack{K<k\le K'\le 2K\\P<kl\le P_1}}\Delta_1(l)\Delta_2(k)e((\xi+m)h(kl))\bigg|^2
\]
\begin{equation}\label{Teclem1est}
\lesssim L\log^3(L)\bigg(\sum_{L<l\le 2L}\bigg|\sum_{\substack{K<k\le K'\le 2K\\P<kl\le P_1}}\Delta_2(k)e((\xi+m)h(kl))\bigg|^2\bigg)\text{.}
\end{equation}
To proceed, we use the Weyl and Van der Corput inequality, see Lemma~3.2 in \cite{MMR}, with $H=K$ and $I=I(l,P,P_1)=\{k\in(K,K']:\,P<kl\le P_1\}\subseteq(K,2K]$ which is an interval, $z_k=\Delta_2(k)e((\xi+m)h(kl))$ and $R\in\mathbb{N}$ which will be specified later. With these choices Lemma~3.2 from \cite{MMR}, for any $l\in(L,2L]$ yields the following bound
\[
\bigg|\sum_{\substack{K<k\le K'\le 2K\\P<kl\le P_1}}\Delta_2(k)e((\xi+m)h(kl))\bigg|^2\le
\frac{K+R}{R}\sum_{|r|\le R}\bigg(1-\frac{|r|}{R}\bigg)\sum_{k,k+r\in I}z_k\overline{z_{k+r}}
\]
\[
=\frac{K+R}{R}\sum_{|r|\le R}\bigg(1-\frac{|r|}{R}\bigg)\sum_{\substack{K<k,k+r\le K'\\ P<kl,(k+r)l\le P_1}}\Delta_2(k)\overline{\Delta_2(k+r)}e\big((\xi+m)h(kl)-(\xi+m)h((k+r)l)\big)\text{.}
\]
Taking into account the estimate above together with $\ref{DefEr}$ and $\ref{Ezerobound}$, we can bound the expression  $\ref{Teclem1est}$ by
\[
L\log^3(L)\frac{K+R}{R}\sum_{|r|\le R}\bigg(1-\frac{|r|}{R}\bigg)E_r
\]
\[
\lesssim\log^3(L)\log^3(K)L^2K\frac{K+R}{R}+L\log^3(L)\frac{K+R}{R}\sum_{1\le|r|\le R}|E_r|\text{.}
\]
Now we focus on bounding the last sum appearing in the above expression. Firstly, let 
\[
S(k,r)=\sum_{\max\{L,\frac{P}{k},\frac{P}{k+r}\}<l\le \min\{2L,\frac{P_1}{k},\frac{P_1}{k+r}\}}e\big((\xi+m)h(kl)-(\xi+m)h((k+r)l)\big)\text{,}
\]
so that for all $r\neq 0$ we have that
\[
E_r=\sum_{\max\{K,K-r\}<k\le \min\{K',K'-r\}}\Delta_2(k)\overline{\Delta_2(k+r)}S(k,r)\text{.}
\]
Then we have that
\[
\sum_{1\le |r|\le R}|E_r|\le\sum_{1\le|r|\le R}\sum_{\max\{K,K-r\}<k\le \min\{K',K'-r\}}|\Delta_2(k)\Delta_2(k+r)S(k,r)|
\]
\[
\lesssim\sum_{1\le|r|\le R}\sum_{\max\{K,K-r\}<k\le \min\{K',K'-r\}}|\Delta_2(k)|^2|S(k,r)|+|\Delta_2(k+r)|^2|S(k,r)|\text{.}
\]
Note that $S(k+r,-r)=\overline{S(k,r)}$ and thus the expression above becomes
\[
\sum_{1\le|r|\le R}\sum_{\max\{K,K-r\}<k\le \min\{K',K'-r\}}|\Delta_2(k)|^2|S(k,r)|+|\Delta_2(k+r)|^2|S(k+r,-r)|
\]
\[
\le
\sum_{1\le|r|\le R}\sum_{K<k,k+r\le K'}|\Delta_2(k)|^2|S(k,r)|+\sum_{1\le|r|\le R}\sum_{K<k,k+r\le K'}|\Delta_2(k+r)|^2|S(k+r,-r)|
\]
\[\lesssim\sum_{1\le|r|\le R}\sum_{K<k,k+r\le K'}|\Delta_2(k)|^2|S(k,r)|=\sum_{K<k\le K'}|\Delta_2(k)|^2\sum_{1\le|r|\le R}1_{(K,K']}(k+r)|S(k,r)|\text{.}
\]
To estimate this we simply apply van der Corput lemma to 
\[
S(k,r)=\sum_{\max\{L,\frac{P}{k},\frac{P}{k+r}\}<l\le \min\{2L,\frac{P_1}{k},\frac{P_1}{k+r}\}}e\big((\xi+m)h(kl)-(\xi+m)h((k+r)l)\text{.}
\]
Let $F(x)=(\xi+m)h(kx)-(\xi+m)h((k+r)x)$, $x\in\big[\max\{L,\frac{P}{k},\frac{P}{k+r}\}, \min\{2L,\frac{P_1}{k},\frac{P_1}{k+r}\}\big]$. We have that
\begin{equation}\label{DIFVan}
|F''(x)|\simeq |r||\xi+m|K\frac{h(KL)}{(KL)^2}=|r||\xi+m|K^{-1}L^{-2}h(KL)\text{.}
\end{equation}
To see why the above estimate holds, we note that 
\[
F''(x)=(\xi+m)k^2h''(kx)-(\xi+m)(k+r)^2h''((k+r)x)=(\xi+m)\big(k^2h''(kx)-(k+r)^2h''((k+r)x)\big)\text{,}
\]
and thus by the Mean Value Theorem applied to the function $f_x(y)=y^2h''(xy)$ on the interval with endpoints $k,k+r$, we get that there exists $\xi_{k,r,x}$ in that interval such that
\[
|k^2h''(kx)-(k+r)^2h''((k+r)x)|=|f_x(k)-f_x(k+r)|=|r||f_x'(\xi_{k,r,x})|
\]
\[
=|r|(2\xi_{k,r,x}h''(x\xi_{k,r,x})+\xi_{k,r,x}^2h'''(x\xi_{k,r,x})x)=|r|(2\xi_{k,r,x}h''(x\xi_{k,r,x})+\xi_{k,r,x}h'''(x\xi_{k,r,x})x\xi_{k,r,x})
\]
\[
=|r|(2\xi_{k,r,x}h''(x\xi_{k,r,x})+\xi_{k,r,x}h''(x\xi_{k,r,x})(c-2+\vartheta_3(x\xi_{k,r,x})))
\]
\[
=|r|\xi_{k,r,x}h''(x\xi_{k,r,x})(c+\vartheta_3(x\xi_{k,r,x}))\simeq|r|Kh''(KL)\simeq |r|K\frac{h(KL)}{(KL)^2}\text{,}
\]
which justifies the estimate $\ref{DIFVan}$, and we note that we have used 2.2 and 2.3 from Lemma~2.1 in \cite{MMR}. Thus
\[
|S(k,r)|\lesssim 
|r|^{1/2}|\xi+m|^{1/2}K^{-1/2}h(KL)^{1/2}+|r|^{-1/2}|\xi+m|^{-1/2}K^{1/2}Lh(KL)^{-1/2}\text{.}
\]
Now we return to 
\[
\sum_{1\le |r|\le R}|E_r|\lesssim \sum_{K<k\le K'}|\Delta_2(k)|^2\sum_{1\le|r|\le R}1_{(K,K']}(k+r)|S(k,r)| 
\]
\[
\lesssim\sum_{K<k\le K'}|\Delta_2(k)|^2\sum_{1\le|r|\le R}\Big(|r|^{1/2}|\xi+m|^{1/2}K^{-1/2}h(KL)^{1/2}+|r|^{-1/2}|\xi+m|^{-1/2}K^{1/2}Lh(KL)^{-1/2}\Big)
\]
\[
\lesssim K\log^3(K)\Big(R^{3/2}|\xi+m|^{1/2}K^{-1/2}h(KL)^{1/2}+R^{1/2}|\xi+m|^{-1/2}K^{1/2}Lh(KL)^{-1/2}\Big)
\]
\[
=K^{1/2}\log^3(K)R^{3/2}|\xi+m|^{1/2}h(KL)^{1/2}+K^{3/2}\log^3(K)R^{1/2}|\xi+m|^{-1/2}Lh(KL)^{-1/2}\text{.}
\]
Combining everything yields
\[
\bigg|\sum_{L<l\le L'\le 2L}\sum_{\substack{K<k\le K'\le 2K\\P<kl\le P_1}}\Delta_1(l)\Delta_2(k)e((\xi+m)h(kl))\bigg|^2
\]
\[
\lesssim\log^3(L)\log^3(K)L^2K\frac{K+R}{R}+L\log^3(L)\frac{K+R}{R}\sum_{1\le|r|\le R}|E_r| 
\]
\begin{multline*}
\lesssim\log^3(L)\log^3(K)L^2K\frac{K+R}{R}
\\
+L\log^3(L)\frac{K+R}{R}\bigg(K^{1/2}\log^3(K)R^{3/2}|\xi+m|^{1/2}h(KL)^{1/2}
\\
+K^{3/2}\log^3(K)R^{1/2}|\xi+m|^{-1/2}Lh(KL)^{-1/2}\bigg)
\end{multline*}
\begin{multline*}
\lesssim\log^3(L)\log^3(K)L^2K\frac{K+R}{R}
\\
+K^{1/2}L\log^3(K)\log^3(L)(K+R)R^{1/2}|\xi+m|^{1/2}h(KL)^{1/2}
\\
+K^{3/2}L^2\log^3(K)\log^3(L)(K+R)R^{-1/2}|\xi+m|^{-1/2}h(KL)^{-1/2}\text{.}
\end{multline*}
Here we choose $R=\lfloor K^{(1-c)/3}L^{(2-c)/3}\rfloor+1$. Note that since $K^{1/2}\lesssim L\lesssim K^2$ we get 
\[
K^{\frac{2}{3}-\frac{c}{2}}=K^{\frac{1}{3}-\frac{c}{3}}K^{\frac{1}{3}-\frac{c}{6}}\lesssim K^{(1-c)/3}L^{(2-c)/3}\lesssim K^{\frac{1}{3}-\frac{c}{3}}K^{\frac{4}{3}-\frac{2c}{3}}=K^{\frac{5}{3}-c}\text{.}
\]
Thus $1\lesssim R\lesssim K$, since $2/3-c/2>0$ and $5/3-c<1$. For that choice of $R$ and using the bound $x^{c-\varepsilon}\lesssim_{\varepsilon}h(x)\lesssim_{\varepsilon} x^{c+\varepsilon}$ (see (2.9) from Lemma~2.6 in \cite{MMR}) the estimate becomes
\[
\bigg|\sum_{L<l\le L'\le 2L}\sum_{\substack{K<k\le K'\le 2K\\P<kl\le P_1}}\Delta_1(l)\Delta_2(k)e((\xi+m) h(kl))\bigg|^2
\]
\[
\lesssim_{\varepsilon}\log^3(L)\log^3(K)\Big(K^{(5+c)/3}L^{(4+c)/3}+|\xi+m|^{1/2}K^{(5+c)/3+\varepsilon}L^{(4+c)/3+\varepsilon}+|\xi+m|^{-1/2}K^{(7-c)/3+\varepsilon}L^{(5-c)/3+\varepsilon}\Big)
\]
\[
\lesssim\log^3(L)\log^3(K)\Big((|m|+1)^{1/2}K^{(5+c)/3+\varepsilon}L^{(4+c)/3+\varepsilon}+|\xi+m|^{-1/2}K^{(7-c)/3+\varepsilon}L^{(5-c)/3+\varepsilon}\Big)\text{,}
\]
and thus for any $\varepsilon>0$ we get
\[
\bigg|\sum_{L<l\le L'\le 2L}\sum_{\substack{K<k\le K'\le 2K\\P<kl\le P_1}}\Delta_1(l)\Delta_2(k)e((\xi+m) h(kl))\bigg|
\]
\[
\lesssim_{\varepsilon}\log^{2}(L)\log^{2}(K)\Big((|m|+1)^{1/4}K^{(5+c)/6+\varepsilon}L^{(4+c)/6+\varepsilon}+|\xi+m|^{-1/4}K^{(7-c)/6+\varepsilon}L^{(5-c)/6+\varepsilon}\Big)\text{.}
\]
Remembering that $ L\lesssim K^2$ we get that
$\log L\lesssim \log K$ for large $K$ and the proof is complete.
\end{proof}

\subsection{Major arc estimates}
Here we prove Lemma~$\ref{Major}$. Through elementary reductions we derive the desired estimate as a corollary of the following lemma. 

\begin{lemma}\label{majarckey}
Let $c\in(1,2)$, $h\in\mathcal{R}_c$, $\theta_1\in(c-1,1)$ and $\varepsilon>0$. Then there exists a positive constant $C=C(\varepsilon,h,\theta_1)$ such that for all $t\ge 1$ and $\xi \in[-t^{-\theta_1},t^{-\theta_1}]$ we have
\begin{equation}\label{dyadapprox}
\Big|\sum_{t/2<n\le t}\Lambda(n)e(h(n)\xi)-\int_{t/2}^{t}e(h(s)\xi)ds\Big|\le Cte^{-(\log t)^{1/3-\varepsilon}}\text{.}
\end{equation}
\end{lemma}

We briefly discuss how Lemma~$\ref{majarckey}$ implies Lemma~$\ref{Major}$.
\begin{proof}[Proof of Lemma~$\ref{Major}$ assuming Lemma~$\ref{majarckey}$]\label{Steps}
We fix $c\in(1,2)$, $h\in\mathcal{R}_c$ and $\theta_1\in(c-1,1)$ and all the implicit constants may depend on these choices. Without loss of generality we may assume that $\xi>0$ since after identifying $\mathbb{T}$ with $[-1/2,1/2)$ one may simply conjugate to get the result for negative values of $\xi$. \\
\textbf{Step 1: Removing the floor.} Note that 
\begin{equation}\label{Nofloor}
\sum_{p\le N}\log(p)e(\lfloor h(p)\rfloor \xi)=\sum_{p\le N}\log(p)e( h(p)\xi)+O(N^{1-\theta_1})\text{,}
\end{equation}
since
\[
\Big|\sum_{p\le N}\log(p)e(\lfloor h(p)\rfloor \xi)-\sum_{p\le N}\log(p)e( h(p)\xi)\Big|\lesssim \sum_{p\le N}\log(p)|\lfloor h(p)\rfloor -h(p)|\xi \lesssim N N^{-\theta_1}=N^{1-\theta_1}\text{.}
\]
\textbf{Step 2: Replacing $\log$ with $\Lambda$.} Note that 
\[
\sum_{p\le N}\log(p)e(h(p)\xi)=\sum_{n\le N}\Lambda(n)e( h(n)\xi)+O(N^{1/2})\text{,}\quad\text{see $\ref{Ltolog}$.}
\]
\textbf{Step 3: Sums over dyadic scales.} According to the previous steps, to establish the estimate $\ref{mainmajarc}$ for $0\le \xi\le N^{-\theta_1}$, it suffices to prove that
\[
\Big|\sum_{n\le N}\Lambda(n)e(h(n)\xi)-\sum_{n\le h(N)}\varphi'(n)e(n\xi)\Big|\lesssim_{\varepsilon} Ne^{-(\log N)^{1/3-\varepsilon}}\text{.}
\]
We achieve this using the estimate $\ref{dyadapprox}$; we have
\[
\sum_{n\le N}\Lambda(n)e(h(n)\xi)=\sum_{0\le i\le \sqrt{\log N}}S_{N,i}(\xi)+O_{\varepsilon}(Ne^{-(\log N)^{1/3-\varepsilon}})
\]
where $S_{N,i}(\xi)=\sum_{N/2^{i+1}< n\le N/2^{i}}\Lambda(n)e(h(n)\xi)$. This is clear since
\[
\Big|\sum_{n\le N}\Lambda(n)e(h(n)\xi)-\sum_{0\le i\le \sqrt{\log N}}S_{N,i}(\xi)\Big|\le\sum_{n\le \frac{N}{2^{\lfloor \sqrt{\log N}\rfloor+1}} }\Lambda(n)\lesssim N 2^{-\lfloor \sqrt{\log N}\rfloor-1} \lesssim_{\varepsilon} Ne^{-(\log N)^{1/3-\varepsilon}}\text{.}
\]
If we assume that Lemma~$\ref{majarckey}$ holds, then for all $i\in\{0,\dotsc,\lfloor \sqrt{\log N}\rfloor\}$ we get  
\begin{equation}\label{usefulstep}
S_{N,i}(\xi)=\int_{N/2^{i+1}}^{N/2^i}e(h(s)\xi)ds+O_{\varepsilon}(N/2^ie^{-(\log (N/2^i))^{1/3-\varepsilon/2}})\text{,}
\end{equation}
since for all such $i$ we have $1\le N 2^{-\sqrt{\log N}} \le N/2^i\le N$, and $|\xi|\le N^{-\theta_1}\le (N/2^i)^{-\theta_1}$. Thus 
\[
\sum_{0\le i\le \sqrt{\log N}}S_{N,i}(\xi)=\int_{N/2^{\lfloor \sqrt{\log N}\rfloor+1}}^Ne(h(s)\xi)ds+O_{\varepsilon}\bigg(\sum_{0\le i\le \sqrt{\log N}}N/2^ie^{-(\log (N/2^i))^{1/3-\varepsilon/2}}\bigg)
\]
\[
=\int_1^Ne(h(s)\xi)ds+O\big(N2^{-\sqrt{\log N}}\big)+O_{\varepsilon}\Big(\sqrt{\log N} N e^{-(\log N)^{1/3-\varepsilon/2}}\Big)=
\int_1^Ne(h(s)\xi)ds+O_{\varepsilon}(Ne^{-(\log N)^{1/3-\varepsilon}})\text{.}
\]
Combining everything we get
\[
\sum_{n\le N}\Lambda(n)e(h(n)\xi)=\int_{1}^{N}e(h(s)\xi)ds+O_{\varepsilon}\big(Ne^{-(\log N)^{1/3-\varepsilon}}\big)\text{.}
\]
\\
\textbf{Step 4: A standard estimate.} The final observation here is the following.
\[
\Big|\sum_{n\le h(N)}\varphi'(n)e(n\xi)-\int_{1}^Ne(h(s)\xi)ds\Big|=\sum_{n=2}^{\lfloor h(N)\rfloor}\int_{n-1}^n\Big|\varphi'(n)e(n\xi)-\varphi'(u)e(u\xi)\Big|du+O(1)
\]
\[
\le\sum_{n=2}^{\lfloor h(N)\rfloor}\int_{n-1}^n\Big|\varphi'(n)e(n\xi)-\varphi'(u)e(n\xi)\Big|du+\sum_{n=2}^{\lfloor h(N)\rfloor}\int_{n-1}^n\Big|\varphi'(u)e(n\xi)-\varphi'(u)e(u\xi)\Big|du+O(1)
\]
\[
\le\sum_{n=2}^{\lfloor h(N)\rfloor}\int_{n-1}^n \big(\varphi'(n-1)-\varphi'(n)\big)du+\sum_{n=2}^{\lfloor h(N)\rfloor}\int_{n-1}^n\varphi'(u)|n\xi-u\xi|du+O(1)
\]
\[
\le O(1)-\varphi'(\lfloor h(N)\rfloor )+|\xi|\int_1^{h(N)}\varphi'(u)du= O(1)+|\xi| N\lesssim N^{1-\theta_1}\text{,}
\]
where we have used that $\varphi'$ is decreasing and $\lim_{x\to
\infty}\varphi'(x)=0$. The proof is complete.
\end{proof}
We now focus on Lemma~$\ref{majarckey}$. Following \cite{Tolev}, the argument relies on estimating certain sums over zeros of the zeta function in the critical strip, see Lemma~$\ref{Keyroots}$. Instead of beginning our analysis collecting all the necessary estimates, we start with the proof of Lemma~$\ref{majarckey}$, letting these expressions naturally appear. The first step is an application of the truncated explicit formula, which we state here for the sake of completeness.

\begin{theorem}\label{explform} Let $2\le T\le x$. Then 
\[
\Psi(x):=\sum_{n\le x}\Lambda(n)=x-\sum_{|\Imm \rho|\le T}\frac{x^{\rho}}{\rho}+O\bigg(\frac{x\log^2 x}{T}\bigg)\text{,}
\]
where $\rho$ are the zeros of the zeta function in the critical strip.
\end{theorem}

\begin{proof} The result is standard and can be found in \cite{Kar}, see page 69, Theorem~3. 
\end{proof}

\begin{proof}[Proof of Lemma~$\ref{majarckey}$]
Again without loss of generality we may assume that $\xi\in[0,t^{-\theta_1}]$. We wish to establish
\[
\Big|\sum_{t/2<n\le t}\Lambda(n)e(h(n)\xi)-\int_{t/2}^{t}e(h(s)\xi)ds\Big|\lesssim_{\varepsilon} te^{-(\log t)^{1/3-\varepsilon}}\text{.}
\]
Let $\Psi(x)=\sum_{n\le x}\Lambda(n)$ and note that by summation by parts we get
\[
\sum_{t/2<n\le t}\Lambda(n)e(h(n)\xi)=\big(\Psi(t)-\Psi(t/2)\big)e(h(t)\xi)-\int_{t/2}^t\big(\Psi(s)-\Psi(t/2)\big)\frac{d}{ds}\big(e(h(s)\xi)\big)ds\text{.}
\]
Now we use Theorem~$\ref{explform}$ with $T=t^{\theta_3}$ where $\theta_3$ is a real number in $(0,1)$ which will be specified later and can only depend on $c$ and $\theta_1$. Since $T\le t/2$ for large $t$, the previous expression becomes
\[
\Big(t/2-\sum_{|\Imm\rho|\le T}\frac{t^{\rho}-(t/2)^{\rho}}{\rho}\Big)e(h(t)\xi)+O\bigg(\frac{t\log^2 t}{T}\bigg)-
\]
\[
-\int_{t/2}^t\Big(s-t/2-\sum_{|\Imm\rho|\le T}\frac{s^{\rho}-(t/2)^{\rho}}{\rho}\Big)\frac{d}{ds}\big(e(h(s)\xi)\big)ds-\int_{t/2}^tO\bigg(\frac{t\log^2 t}{T}\bigg)\frac{d}{ds}\big(e(h(s)\xi)\big)ds\text{.}
\]
Note that
\[
\int_{t/2}^t\bigg|\frac{t\log^2 t}{T}\frac{d}{ds}\big(e(h(s)\xi)\big)\bigg|ds\lesssim\frac{t\log^2 t}{T}\int_{t/2}^t \xi h'(s)ds\lesssim \frac{t\log^2 t}{T} \xi h(t)\text{.}
\]
Thus we have
\[
\sum_{t/2<n\le t}\Lambda(n)e(h(n)\xi)=\Big(t/2-\sum_{|\Imm\rho|\le T}\frac{t^{\rho}-(t/2)^{\rho}}{\rho}\Big)e(h(t)\xi)-
\]
\[
-\int_{t/2}^t\Big(s-t/2-\sum_{|\Imm\rho|\le T}\frac{s^{\rho}-(t/2)^{\rho}}{\rho}\Big)\frac{d}{ds}\big(e(h(s)\xi)\big)ds+O\bigg(\frac{t\log^2 t}{T}\big(1+\xi h(t)\big)\bigg)\text{.}
\]
We use integration by parts to obtain
\[
\sum_{t/2<n\le t}\Lambda(n)e(h(n)\xi)=\int_{t/2}^t\Big(1-\sum_{|\Imm\rho|\le T}s^{\rho-1}\Big)e(h(s)\xi)ds+O\bigg(\frac{t\log^2 t}{T}\big(1+\xi h(t)\big)\bigg)
\]
\[
=\int_{t/2}^te(h(s)\xi)ds-\sum_{|\Imm\rho|\le T}\int_{t/2}^ts^{\rho-1}e(h(s)\xi)ds+O\bigg(\frac{t\log^2 t}{T}\big(1+\xi h(t)\big)\bigg)\text{.}
\]
For the error term, firstly note that since $0\le\xi\le t^{-\theta_1}$ we get that for any $\varepsilon>0$
\[
\frac{t\log^2 t}{T}\big(1+\xi h(t)\big)\lesssim_{\varepsilon} t^{(1-\theta_3)+\varepsilon}+t^{(c-\theta_1)+(1-\theta_3)+\varepsilon}\text{.}
\]
Since $c-\theta_1<1$, it is clear that we can choose $\varepsilon$ close to $0$ and $\theta_3$ close to $1$ in order to get a polynomial saving. For the sake of completeness, let us mention that for the following choices  
\[
\theta_3=1-\frac{1-(c-\theta_1)}{4}\in(0,1)\quad\text{and}\quad \varepsilon=\frac{1-(c-\theta_1)}{4}>0
\]
we get that there exists $\chi=\chi(c,\theta_1)>0$ such that 
\[
\frac{t\log^2 t}{T}\big(1+\xi h(t)\big)\lesssim t^{1-\chi}\lesssim_{\varepsilon} te^{-(\log t)^{1/3-\varepsilon}}\text{.}
\] 
To conclude the proof, it suffices to establish 
\begin{equation}\label{weirdtruncroot}
\bigg|\sum_{|\Imm\rho|\le T}\int_{t/2}^ts^{\rho-1}e(h(s)\xi)ds\bigg|\lesssim_{\varepsilon} te^{-(\log t)^{1/3-\varepsilon}}\text{.}
\end{equation}
To do this, we begin by collecting some useful technical lemmas.

\begin{lemma}\label{zerofreerg}
Let $1\le R\le Y$ and define
\[
\theta(R)=\frac{1}{\log^{2/3}(R+10)\log\log(R+10)}\text{.}
\]
There exists a constant $C>0$ such that
\[
\sum_{0<\gamma\le R}Y^{\beta}\lesssim 
\left\{
\begin{array}{ll}
      Y^{1/2}R\log^6(Y)\text{,} & Y^{3/4}\le R \le Y\text{,} \\
      \exp(2\log(Y)+3\log(R)-2\sqrt{3\log(Y)\log(R)}\log^6 (Y)\text{,} & Y^{1/3}\le R\le Y^{3/4}\text{,}\\
      Y^{1-C\theta(R)}R^{\frac{12C}{5}\theta(R)}\log^{45}Y\text{,}&1\le R \le Y^{1/3}\text{,} \\
\end{array} 
\right. 
\]
where the implied constant does not depend on $R$ or $Y$. 
\end{lemma}
\begin{proof}
See Lemma~5 in \cite{Tolev}. 
\end{proof}
\begin{lemma}\label{MLforSrho}
Let $\theta_3\in(0,1)$ and $\varepsilon>0$. Then there exists a positive constant $C=C(\theta_3,\varepsilon)$, such that for all $t\ge 1$ and $1\le T_1\le t^{\theta_3}$ we have
\[
\frac{1}{\sqrt{T_1}}\sum_{0<\gamma\le T_1}t^{\beta}\le C te^{-(\log t)^{1/3-\varepsilon}}\text{.}
\]
\end{lemma}
\begin{proof}
The proof is identical to the one given for Lemma~11 in \cite{Tolev}. For the convenience of the reader, and to make sure that our parameters allow for the argument to work here, we provide the details. We use Lemma~$\ref{zerofreerg}$ with $R=T_1$ and $Y=t$. There exists $C>0$ such that
\[
\sum_{0<\gamma\le T_1}t^{\beta}\lesssim 
\left\{
\begin{array}{ll}
      t^{1/2}T_1\log^6(t)\text{,} & t^{3/4}\le T_1 \le t\text{,} \\
      \exp(2\log(t)+3\log(T_1)-2\sqrt{3\log(t)\log(T_1)}\log^6 (t)\text{,} & t^{1/3}\le T_1\le t^{3/4}\text{,}\\
      t^{1-C\theta(T_1)}T_1^{\frac{12C}{5}\theta(T_1)}\log^{45}t\text{,}&1\le T_1 \le t^{1/3}\text{,} \\
\end{array} 
\right. 
\]
where $\theta(T_1)$ is defined in Lemma~$\ref{zerofreerg}$, and the analysis is naturally split into three cases.\\\,\\
\textbf{Case 1.} If $1\le T_1 \le t^{1/3}$, we have
\begin{equation}\label{Case1}
\sum_{0<\gamma\le T_1}t^{\beta}\lesssim t^{1-C\theta(T_1)}T_1^{\frac{12C}{5}\theta(T_1)}\log^{45}t\le t^{1-C\theta(T_1)}t^{\frac{12C}{15}\theta(T_1)}\log^{45}t =t^{1-\frac{C}{5}\theta(T_1)}\log^{45}t\le t^{1-\frac{C}{5}\theta(t^{1/3})}\log^{45}t\text{,}
\end{equation}
since $\theta$ is decreasing. Finally, it is easy to see that for large $t$ we have
\[
\frac{C}{5}\frac{\log t}{\log^{2/3}(t^{1/3}+10)\log\log(t^{1/3}+10)}-\log\log^{45}t\gtrsim_{\varepsilon} (\log t)^{1/3-\varepsilon}\text{,}
\]
thus
\[
-\frac{C}{5}\frac{\log t}{\log^{2/3}(t^{1/3}+10)\log\log(t^{1/3}+10)}+\log\log^{45}t\lesssim_{\varepsilon} -(\log t)^{1/3-\varepsilon}\text{,}
\]
and this implies that
\[
t^{-\frac{C}{5}\theta(t^{1/3})}\log^{45}t=t^{-\frac{C}{5}\frac{1}{\log^{2/3}(t^{1/3}+10)\log\log(t^{1/3}+10)}}\log^{45}t\lesssim_{\varepsilon} e^{-(\log t)^{1/3-\varepsilon}}\text{.}
\]
Using this together with $\ref{Case1}$, we get that
\[
\frac{1}{\sqrt{T_1}}\sum_{0<\gamma\le T_1}t^{\beta}\le \sum_{0<\gamma\le T_1}t^{\beta}\lesssim_{\varepsilon}te^{-(\log t)^{1/3-\varepsilon}} \text{.}
\]
\\
\textbf{Case 2.} If $t^{1/3}\le T_1\le t^{3/4}$, we have that
\[
\frac{1}{\sqrt{T_1}}\sum_{0<\gamma\le T_1}t^{\beta}\lesssim \frac{\log^6 t}{\sqrt{T_1}}\exp(2\log(t)+3\log(T_1)-2\sqrt{3\log(t)\log(T_1)} 
\]
\[
=\log^6 t\exp(2\log(t)+5/2\log(T_1)-2\sqrt{3\log(t)\log(T_1)}
\]
and 
\[
\max\Big\{\frac{5}{2}\log(T_1)-2\sqrt{3\log(t)\log(T_1)}:\,T_1\in[t^{1/3},t^{3/4}]\Big\}
\]
\[
=\max\Big\{\frac{5}{2}s^2-2\sqrt{3\log(t)}s:\,\sqrt{\log (t^{1/3})}\le s\le\sqrt{\log (t^{3/4})}]\Big\}
\]
\[
=\max\Big\{\frac{5}{2}\log(t^{1/3})-2\sqrt{3\log(t)\log(t^{1/3}}),\frac{5}{2}\log(t^{3/4})-2\sqrt{3\log(t)\log(t^{3/4}})\Big\}
=-\frac{9}{8}\log(t)\text{.}
\]
Thus
\[
\frac{1}{\sqrt{T_1}}\sum_{0<\gamma\le T_1}t^{\beta}\lesssim \log^{6}(t)\exp(2\log(t)-\frac{9}{8}\log(t))=\log^{6}(t)\exp(7/8\log(t))=\log^{6}(t)t^{7/8}\lesssim_{\varepsilon} te^{-(\log t)^{1/3-\varepsilon}}\text{.}
\]
\textbf{Case 3.} If $t^{3/4}\le T_1\le t$, we have that
\[
\frac{1}{\sqrt{T_1}}\sum_{0<\gamma\le T_1}t^{\beta}\lesssim t^{1/2}T_1^{1/2}\log^6(t) \le  t^{1/2+\theta_3/2}\log^6(t)\lesssim_{\varepsilon,\theta_3} te^{-(\log t)^{1/3-\varepsilon}}\text{,}
\]
since $1/2+\theta_3/2<1$. The proof of the lemma is complete.
\end{proof}
\begin{remark}\label{Polsaving}We note that Case 1 in the previous proof is the only part that forces us to lose the polynomial saving in our error terms. Moreover, assuming that $\zeta$ has a zero free strip upgrades the estimate to a polynomially saving one. More precisely, if we assume that there exists $\delta\in(0,1/2)$ such that $\zeta(s)\neq 0$ for $\Ree(s)\in[1-\delta,1]$, then we get
\[
\frac{1}{\sqrt{T_1}}\sum_{0<\gamma\le T_1}t^{\beta}=\frac{1}{\sqrt{T_1}}\sum_{\substack{0<\gamma\le T_1\\0<\beta\le4/5}}t^{\beta}+\frac{1}{\sqrt{T_1}}\sum_{\substack{0<\gamma\le T_1\\4/5\le \beta<1-\delta}}t^{\beta}\lesssim T_1^{-1/2}t^{4/5}N(T_1)+T_1^{-1/2}t^{1-\delta}N(4/5,T_1)\text{,}
\]
where $N(T_1)$ is the number of nontrivial zeros of $\zeta$ with imaginary part in $[-T_1,T_1]$ and $N(4/5,T_1)$ is the number of nontrivial zeros with imaginary part in $[-T_1,T_1]$ and real part in $[4/5,1)$ counted with multiplicity. It is a standard result that $N(T_1)\lesssim T_1\log(T_1)$ and by \cite{Ign} we get that
\[
N(4/5,T_1)\lesssim T_1^{\frac{3(1-\frac{4}{5})}{2-\frac{4}{5}}}\log^5(T_1)=T_1^{1/2}\log^5(T_1)\text{,}
\] 
and thus
\[
\frac{1}{\sqrt{T_1}}\sum_{0<\gamma\le T_1}t^{\beta}\lesssim T_1^{1/2}t^{4/5}\log(T_1)+t^{1-\delta/2}\lesssim t^{1/6}t^{4/5}\log(T_1)+t^{1-\delta/2}\lesssim t^{1-\min\{\delta/2,1/31\}}\text{.}
\]
Finally, we remark that the zero-free strip assumption is in some sense optimal because if we assume that roots with real part arbitrarily close to $1$ exist, then we do not get polynomial saving here.
\end{remark}
We are ready to prove the key estimate $\ref{weirdtruncroot}$. 
\begin{lemma}\label{Keyroots}
Let $c\in(1,2)$, $h\in\mathcal{R}_c$, $\theta_3\in(0,1)$ and $\varepsilon>0$. Then there exists a positive constant $C=C(h,\theta_3,\varepsilon)$ such that for all $t\ge 1$ and $\xi\in\mathbb{R}$, we get that 
\[
\bigg|\sum_{|\Imm\rho|\le t^{\theta_3}}\int_{t/2}^ts^{\rho-1}e(h(s)\xi)ds\bigg|\le C t e^{-(\log t)^{1/3-\varepsilon}}\text{.}
\]
\end{lemma}
\begin{proof}
We may assume without loss of generality that $\xi\ge 0$ and $t$ is large.
For convenience, let $T:=t^{\theta_3}$. For any root of the zeta function in the critical strip $\rho=\beta+i\gamma$, any $t\ge 1$ and $\xi\ge 0$, we define
\[
I_{t,\xi}(\rho):=\int_{t/2}^t s^{\rho-1}e(h(s)\xi)ds=\int_{t/2}^t s^{\beta-1}e^{i\gamma\log s}e^{2\pi i h(s)\xi}ds=\int_{t/2}^t s^{\beta-1}e\Big(h(s)\xi+\frac{\gamma}{2\pi}  \log s\Big)ds\text{.}
\]
For every $s\in[t/2,t]$ we also define the functions
\[
G(s)=s^{\beta-1}\quad\text{and}\quad F(s)=h(s)\xi+\frac{\gamma}{2\pi}\log s\text{,}\quad\text{so that}\quad I_{t,\xi}(\rho)=\int_{t/2}^t G(s)e(F(s))ds\text{.}
\]
It is clear that $G$ and $F$ depend on $t,\xi,h,\rho$ as well, but we suppress that dependence to keep the exposition reasonable. 

Ultimately, $I_{t,\xi}(\rho)$ will be bounded using van der Corput type estimates, and more specifically Lemma~4.3 and 4.5 from \cite{VanMon}, see pages 71-72. To apply such results, one needs to have control of the sign and the size of the derivatives of $F$, and it is therefore natural to break the analysis into three regions. We define
\[
M_T=\{\rho=\beta+i\gamma:\,\beta\in(0,1),\,\gamma\in[-T,T]\text{ and }\zeta(\rho)=0\}=M_1 \cup M_2\cup M_3
\]
where 
\begin{align*}M_1&=\big\{z\in M_T:\,\frac{-\Imm z}{2\pi}\in\Big(\frac{3}{2}\xi th'(t),\infty\Big)\big\}\text{,}\quad M_2=\big\{z\in M_T:\,\frac{-\Imm z}{2\pi}\in\Big[\frac{1}{4}\xi th'(t/2),\frac{3}{2}\xi th'(t)\Big]\big\}\text{,}\\
M_3&=\big\{z\in M_T:\,\frac{-\Imm z}{2\pi}\in\Big(-\infty,\frac{1}{4}\xi th'(t/2)\Big)\big\}\text{.}
\end{align*}
We have that 
\[
\Big|\sum_{|\Imm\rho|\le T}I_{t,\xi}(\rho)\Big|\le \Big|\sum_{\rho\in M_1}I_{t,\xi}(\rho)\Big|+ \Big|\sum_{\rho\in M_2}I_{t,\xi}(\rho)\Big|+\Big|\sum_{\rho\in M_3}I_{t,\xi}(\rho)\Big|\eqqcolon E_1+E_2+E_3\text{,}
\]
and we estimate them separately.\\\,\\
\textbf{Estimates for $E_1$.} Note that for all $s\in[t/2,t]$ we have
\[
F'(s)=h'(s)\xi+\frac{\gamma}{2\pi}s^{-1}< h'(t)\xi-\frac{3}{2}\xi t h'(t)s^{-1}\le h'(t)\xi-\frac{3}{2}\xi  h'(t)=- \frac{\xi}{2}h'(t)\le 0\text{ and } |G(s)|\le (t/2)^{\beta-1}\text{.}
\] 
Note that $G/F'$ is decreasing, since
\[
\frac{G(s)}{F'(s)}=\frac{s^{\beta-1}}{\xi h'(s)+\frac{\gamma}{2\pi s}}=\Big(\xi s^{1-\beta}h'(s)+\frac{\gamma}{2\pi s^{\beta}}\Big)^{-1}
\]
and one can see that both $\xi s^{1-\beta}h'(s)$ and $\frac{\gamma}{2\pi s^{\beta}}$ are increasing in $s$, the latter is decreasing in absolute value and negative. Finally, note that
\[
F'(s)\le h'(t)\xi+\frac{\gamma}{2\pi}t^{-1}<0\text{,}
\]
where we used the fact that $\frac{\gamma}{2\pi }<0$, thus
\[
\frac{F'(s)}{G(s)}\le \frac{h'(t)\xi+\frac{\gamma}{2\pi}t^{-1}}{(t/2)^{\beta-1}}<0 \text{.}
\]
Applying Lemma~4.3 from \cite{VanMon}, we get
\[
|I_{t,\xi}(\rho)|=
\Big|\int_{t/2}^t G(s)e(F(s))ds\Big|\le 4\frac{(t/2)^{\beta-1}}{-h'(t)\xi-\frac{\gamma}{2\pi}t^{-1}} \lesssim \frac{t^{\beta}}{-th'(t)\xi-\frac{\gamma}{2\pi}}\lesssim\frac{t^{\beta}}{|\gamma|}\text{,}
\]
where for the last inequality we used that 
\[
-\frac{3}{2}\bigg(\frac{\gamma}{2\pi}\bigg)-\frac{3}{2}\xi th'(t)>-\frac{1}{2}\bigg(\frac{\gamma}{2\pi}\bigg)\text{.}
\]
We know that for all $\rho\in M_1$, $\gamma=\Imm(\rho)<0$, and thus
\begin{equation}\label{FinE1step}
E_1\le\sum_{\rho\in M_1}|I_{t,\xi}(\rho)|= \sum_{\substack{\rho\in M_1:\\ 0<-\gamma\le T}}\frac{t^{\beta}}{|\gamma|}\le\sum_{0<\gamma\le T}\frac{t^{\beta}}{\gamma}\text{,}
\end{equation}
where we have used the symmetries of the roots. Using summation by parts together with standard estimates for the number of zeros of the zeta function in the critical strip, one may obtain the following estimate
\[
\sum_{0<\gamma\le T}\frac{t^{\beta}}{\gamma}\lesssim  \log t \max_{1\le T_1\le T}\frac{1}{T_1}\sum_{0<\gamma \le T_1}t^{\beta}\text{.}
\]
To see why this is true, let $0<\gamma_1\le\gamma_2\le,\dotsc,\gamma_{N(T)}$ be an enumeration of the imaginary parts of roots of the zeta function in $(0,1)\times (0,T)$. For convenience, let $N:=N(T)$. Then summation by parts yields
\[
\sum_{0<\gamma\le T}\frac{t^{\beta}}{\gamma}=\sum_{k=1}^{N}\frac{t^{\beta_k}}{\gamma_k}=\bigg(\sum_{k=1}^{N}t^{\beta_k}\bigg)\frac{1}{\gamma_{N}}+\sum_{k=1}^{N-1}\bigg(\sum_{m=1}^{k}t^{\beta_m}\bigg)\bigg(\frac{1}{\gamma_k}-\frac{1}{\gamma_{k+1}}\bigg)\text{.}
\]
Letting $U_t(n)=\sum_{k=1}^n t^{\beta_k}$, the last expression becomes
\[
\frac{N}{\gamma_N}\Big(\frac{1}{N}U_t(N)\Big)+\sum_{k=1}^{N-1}\Big(\frac{1}{k}U_t(k)\Big)\bigg(\frac{k}{\gamma_k}-\frac{k}{\gamma_{k+1}}\bigg)\text{.}
\]
Now we use the fact that $\gamma_n\simeq \frac{2\pi n}{\log n}$, $N=N(T)\simeq \frac{T\log T}{2\pi}$ to obtain on the one hand that $\frac{N}{\gamma_N}\lesssim \log N$, and on the other hand that for all $k\ge 2$ such that $\gamma_k\ge 1$ we get
\[
\frac{1}{k}U_t(k)=\frac{1}{k}\sum_{m=1}^{k}t^{\beta_m}\le(k^{-1}\gamma_k)\frac{1}{\gamma_k}\sum_{0<\gamma\le\gamma_k}t^{\beta}\lesssim \frac{1}{\log k}\max_{1\le T_1\le T}\frac{1}{T_1}\sum_{0<\gamma\le T_1}t^{\beta}\text{.}
\]
Fix an integer $k_0\ge 2$ such that $\gamma_{k_0}\ge 1$. Combining everything we get that
\begin{multline*}
\sum_{0<\gamma\le T}\frac{t^{\beta}}{\gamma}\lesssim_{\varepsilon}
\\
\lesssim_{\varepsilon}\log N\frac{1}{\log (N)}\max_{1\le T_1\le T}\frac{1}{T}\sum_{0<\gamma<T_1}t^{\beta}+\sum_{k=2}^{N-1}\frac{1}{\log(k)}\bigg(\frac{k}{\gamma_k}-\frac{k}{\gamma_{k+1}}\bigg)\max_{1\le T_1\le T}\frac{1}{T}\sum_{0<\gamma<T_1}t^{\beta}+te^{-(\log t)^{1/3-\varepsilon}}
\\
=\bigg(1+\sum_{k=2}^{N-1}\frac{1}{\log(k)}\Big(\frac{k}{\gamma_k}-\frac{k}{\gamma_{k+1}}\Big)\bigg)\cdot \bigg(\max_{1\le T_1\le T}\frac{1}{T}\sum_{0<\gamma<T_1}t^{\beta}\bigg)+te^{-(\log t)^{1/3-\varepsilon}}\text{,}
\end{multline*}
note that we have absorbed all the $k< k_0$ in the error term. In fact, the error term here has polynomial saving in $t$ since there are finitely many roots in $(0,1)\times (0,\gamma_{k_0}]$, and each one of them has real part less than $1$. Finally, we have
\[
\sum_{k=2}^{N-1}\frac{1}{\log(k)}\Big(\frac{k}{\gamma_k}-\frac{k}{\gamma_{k+1}}\Big)=\sum_{k=2}^{N-1}\frac{1}{\log(k)}\Big(\frac{k}{\gamma_k}-\frac{k+1}{\gamma_{k+1}}\Big)+\sum_{k=2}^{N-1}\frac{1}{\log (k)\gamma_{k+1}}\text{,}
\]
and the second sum is bounded by a multiple of $\sum_{k=1}^N \frac{1}{k}\lesssim \log (N)\lesssim \log(T)$. For the first sum we note that
\[
\bigg|\sum_{k=2}^{N-1}\frac{1}{\log(k)}\Big(\frac{k}{\gamma_k}-\frac{k+1}{\gamma_{k+1}}\Big)\bigg|\le \frac{1}{\log N}\bigg|\frac{1}{\gamma_2}-\frac{N+1}{\gamma_{N+1}}\bigg|+\sum_{k=2}^{N-1}\bigg|\frac{1}{\gamma_1}-\frac{k+1}{\gamma_{k+1}}\bigg|\bigg(\frac{1}{\log(k)}-\frac{1}{\log(k+1)}\bigg)
\]
\[
\lesssim 1+\sum_{k=2}^{N-1}\log(k)\bigg(\frac{1}{\log(k)}-\frac{1}{\log(k+1)}\bigg)\lesssim\log(N)\lesssim \log (T)\text{.}
\]
Thus we have shown that
\begin{equation}\label{IntEst}
\sum_{0<\gamma\le T}\frac{t^{\beta}}{\gamma}\lesssim_{\varepsilon}  \log T \max_{1\le T_1\le T}\frac{1}{T_1}\sum_{0<\gamma \le T_1}t^{\beta}+te^{-(\log t)^{1/3-\varepsilon}}\lesssim \log t \max_{1\le T_1\le T}\frac{1}{\sqrt{T_1}}\sum_{0<\gamma \le T_1}t^{\beta} +te^{-(\log t)^{1/3-\varepsilon}}\text{,}
\end{equation}
and one can easily see that using Lemma~$\ref{MLforSrho}$ with $\varepsilon'=\varepsilon/2$, we obtain
\begin{equation}\label{SBPgamma}
\sum_{0<\gamma\le T}\frac{t^{\beta}}{\gamma}\lesssim_{\varepsilon} te^{-(\log t)^{1/3-\varepsilon}}\text{.}
\end{equation}
Returning to $\ref{FinE1step}$, we get that for all $\varepsilon>0$ we have that $E_1\lesssim_{\varepsilon}te^{-(\log t)^{1/3-\varepsilon}}$ and the bound for $E_1$ is complete.\\\,\\
\textbf{Estimates for $E_3$.} Let us firstly assume that $h(t)^{-1}\le \xi$. For all $\rho=\beta+i\gamma \in M_3$ and all $s\in[t/2,t]$ we have
\[
F'(s)=s^{-1}\Big(sh'(s)\xi+\frac{\gamma}{2\pi}\Big)>s^{-1}\Big(sh'(s)\xi-\frac{1}{4}\xi t h'(t/2)\Big)
\]
\[
\ge s^{-1}\xi\Big(\frac{t}{2}h'(t/2)-\frac{t}{4} h'(t/2)\Big)=s^{-1}\xi \frac{t}{4}h'(t/2)\ge \frac{\xi h'(t/2)}{4}>0\text{.}
\]
Let 
\[
L(s)=\int_{t/2}^s e^{iF(w)}dw\text{,}
\]
we wish to apply Lemma~4.2 in \cite{VanMon}, so we need to check the monotonicity of $F'$. We have that $F''(s)=\xi h''(s)-\frac{\gamma}{2\pi}s^{-2}$. We show that $F''$ has at most one root and thus changes sign at most once. This means that we may apply Lemma~4.2 in \cite{VanMon} at most twice, to obtain
\[
|L(s)|\lesssim \frac{1}{h'(t/2)\xi+\frac{\gamma}{\pi t}} \text{.}
\]
Notice that $F''(s)=0\iff s^2h''(s)=\frac{\gamma}{2\pi \xi}$. For large $t$, this has at most one solution in $[t/2,t]$, since $s^2h''(s)$ is eventually increasing. To see this, note that by Lemma~2.1 in \cite{MMR} we have that
\[
s^2h''(s)=sh'(s)(c-1+\vartheta_2(s))
\]
where $\vartheta_2(s)=\vartheta_1(s)+\frac{s\vartheta_1'(s)}{c+\vartheta_1(s)}$. Now note that 
\[
\frac{d}{ds}\Big(s^2h''(s)\Big)=h'(s)(c-1+\vartheta_2(s))+sh''(s)(c-1+\vartheta_2(s))+h'(s)s\vartheta_2'(s)\text{,}
\]
and for large $t$ and $s\in[t/2,2]$ we have that 
\[
h'(s)(c-1+\vartheta_2(s))+sh''(s)(c-1+\vartheta_2(s))\gtrsim_c h'(s)>0\text{.}
\]
On the other hand, 
\[
s\vartheta_2'(s)=s\vartheta_1'(s)+s\bigg(\frac{(\vartheta_1'(s)+s\vartheta_1''(s))(c+\vartheta_1(s))-\vartheta_1'(s)s\vartheta_1'(s)}{(c+\vartheta_1(s))^2}\bigg)=
\]
\[
=s\vartheta_1'(s)+\bigg(\frac{(s\vartheta_1'(s)+s^2\vartheta_1''(s))(c+\vartheta_1(s))-(s\vartheta_1'(s))^2}{(c+\vartheta_1(s))^2}\bigg)\to 0\text{ as }s\to\infty
\]
since $\vartheta_1=\vartheta$ from the definition of $\mathcal{R}_c$, and we have that $\vartheta(x),x\vartheta'(x),x^2\vartheta''(x)\to 0$ as $x\to0$. This means that after fixing $h$, we can choose $x_0=x_0(h)$ such that for all $t\ge x_0$ we have that
\[
\frac{d}{ds}\Big(s^2h''(s)\Big)=h'(s)(c-1+\vartheta_2(s))+sh''(s)(c-1+\vartheta_2(s))+h'(s)s\vartheta_2'(s)\gtrsim h'(s)>0
\] 
and thus $s^2h''(s)$ is monotone in $[t/2,2]$. This implies that $F''$ has at most one root in $[t/2,t]$ and we get the desired estimates for $L$. 

By integration by parts we get
\[
\bigg|\int_{t/2}^t G(s)e^{iF(s)}ds\bigg|=\bigg|\int_{t/2}^t G(s)L'(s)ds\bigg|\le \bigg|\bigg(G(s)L(s)\bigg)_{s=t/2}^{s=t}\bigg|+\bigg|\int_{t/2}^t G'(s)L(s)ds\bigg| 
\]
\[
\lesssim\sup_{s\in{[t/2,t]}}|L(s)|\bigg(|G(t)|+\int_{t/2}^t|G'(s)|ds\bigg)\lesssim \frac{t^{\beta-1}}{h'(t/2)\xi+\frac{\gamma}{\pi t}}\text{.}
\]
Thus 
\[
E_3\le \sum_{\rho\in M_3}\frac{t^{\beta}}{t/2h'(t/2)\xi+\frac{\gamma}{2\pi}}=\sum_{\rho\in M_3 :\,\frac{-2\pi}{4} th'(t/2)\xi<\gamma\le h(t)\xi}\frac{t^{\beta}}{t/2h'(t/2)\xi+\frac{\gamma}{2\pi}}+\sum_{\rho :\,h(t)\xi<\gamma\le T}\frac{t^{\beta}}{t/2h'(t/2)\xi+\frac{\gamma}{2\pi}} 
\]
\[
\lesssim\sum_{\rho\in M_3 :\,\frac{-2\pi}{4} th'(t/2)\xi<\gamma\le h(t)\xi}\frac{t^{\beta}}{t/4h'(t/2)\xi}+\sum_{\rho :\,h(t)\xi<\gamma\le T}\frac{t^{\beta}}{\gamma}\text{,}
\]
where for the first summand we used that for all such $\gamma$'s we have
\[
t/2h'(t/2)\xi+\frac{\gamma}{2\pi}>t/2h'(t/2)\xi-\frac{t}{4}h'(t/2)\xi=\frac{t}{4}h'(t/2)\xi\text{.}
\]
Thus from the basic properties of $h$ together with the symmetry of the roots of the zeta function, we know that there exists $C=C(h)>0$ such that for all $\varepsilon>0$ we have
\[
E_3\lesssim \frac{1}{h(t)\xi}\sum_{\rho\in M_3 :\,0<\gamma\le C h(t)\xi}t^{\beta}+\sum_{0<\gamma\le T}\frac{t^{\beta}}{\gamma}+O_{\varepsilon}\Big(t e^{-(\log t)^{1/3-\varepsilon}}\Big)\text{,}
\]
and since $h(t)\xi\ge 1$, we have
\[
E_3\lesssim \log t \max_{1\le T_1\le T}\frac{1}{T_1}\sum_{0<\gamma \le T_1}t^{\beta}+O_{\varepsilon}\Big(t e^{-(\log t)^{1/3-\varepsilon}}\Big)=O_{\varepsilon}\Big(t e^{-(\log t)^{1/3-\varepsilon}}\Big)\text{.}
\]
For the first sum this is clear and for the second one see $\ref{IntEst}$ and $\ref{SBPgamma}$. 

Now for the case $h(t)\xi\le 1$, we treat $E_3$ differently. Firstly, note that there exists $\tilde{C}=\tilde{C}(h)$ such that $t/2h'(t/2)\xi\le \tilde{C}$. Split $M_3$ as follows
\[
M_3=M_3'\cup M_3'':=\big\{z\in M_3:\,\frac{\Imm z}{2\pi}\in(-2\tilde{C},2\tilde{C})\big\}\cup \big\{z\in M_3:\,\frac{\Imm z}{2\pi}\in\mathbb{R}\setminus (-2\tilde{C},2\tilde{C})\big\}\text{,}
\]
and clearly
\begin{equation}\label{Easysest}
\sum_{\rho\in M_3'}I_{t,\xi}(\rho)\lesssim te^{-(\log t)^{1/3}}
\end{equation}
since $|I_{t,\xi}(\rho)|\lesssim t^{\beta}$, and we have finitely many roots in $M_3'$ with the largest real part of the roots there smaller than $1$. On the other hand, for every $\rho=\beta+i\gamma\in M_3''$ we have
\[
\bigg|\frac{\gamma}{2\pi}+\frac{t}{2}h'(t/2)\xi\bigg|\ge \Big|\frac{\gamma}{2\pi}\Big|-\tilde{C}\ge \Big|\frac{\gamma}{2\pi}\Big|-\frac{1}{2}\Big|\frac{\gamma}{2\pi}\Big|=\frac{1}{2}\Big|\frac{\gamma}{2\pi}\Big|
\]
and an argument identical to the case of $h(t)\xi\ge 1$ (i.e.: applying Lemma~4.2 from \cite{VanMon} in the same manner) yields
\[
E_3\lesssim 
\sum_{\rho\in M_3''}\bigg|\frac{t^{\beta}}{t/2h'(t/2)\xi+\frac{\gamma}{2\pi}}\bigg|\lesssim \sum_{\rho :\,2\tilde{C}<|\gamma|\le T}\frac{t^{\beta}}{|\gamma|} 
\]
\[
\lesssim\log t \max_{1\le T_1\le T}\frac{1}{T_1}\sum_{0<\gamma \le T_1}t^{\beta}+O_{\varepsilon}\Big(t e^{-(\log t)^{1/3-\varepsilon}}\Big)=O_{\varepsilon}\Big(t e^{-(\log t)^{1/3-\varepsilon}}\Big)\text{,}
\]
The estimates for $E_3$ are complete.\\\,\\
\textbf{Estimates for $E_2$.} For all $\rho=\beta+i\gamma \in M_2$ and all $s\in[t/2,t]$ we have 
\[
F''(s)=\xi h''(s)-\frac{\gamma}{2\pi s^2}>0\text{,}
\]
since $\gamma<0$. For large $t$ and $s\in[t/2,t]$, we have that $F''(s)\simeq \xi h(t)t^{-2}+\xi h'(t)t^{-1}\simeq \xi h(t)t^{-2}$, where the implied constants depend only on the function $h$. Letting $L$ be the same function as before, we obtain the same estimate
\[
\bigg|\int_{t/2}^t G(s)e^{iF(s)}ds\bigg|\lesssim 
t^{\beta-1}\sup_{s\in{[t/2,t]}}|L(s)|\text{,}
\]
and to bound $L$ we apply Lemma~4.4 from \cite{VanMon} and we obtain
\[
|L(s)|\lesssim \frac{1}{\sqrt{\xi h(t)t^{-2}}}\text{,}
\]
and thus
\[
E_2\lesssim \sum_{\rho\in M_2}\frac{t^{\beta-1}}{\sqrt{\xi h(t)t^{-2}}}=\sum_{\rho\in M_2}\frac{t^{\beta}}{\sqrt{\xi h(t)}}=\frac{1}{\sqrt{\xi h(t)}}\sum_{\rho\in M_2}t^{\beta}\lesssim\frac{1}{\sqrt{2\xi t h'(t)}}\sum_{\rho\in M_2:\,0<-\gamma\le 2\xi t h'(t)}t^{\beta}\text{.}
\]
Now we distinguish two cases similarly to the treatment of $E_3$.

Firstly, assume that $2\xi  th'(t)\ge 1$. Using the estimate above and the symmetry of the roots we have 
\[
E_2\lesssim \max_{1\le T_1\le T}\frac{1}{\sqrt{T_1}}\sum_{0<\gamma\le T_1}t^{\beta}\lesssim_{\varepsilon}t e^{-(\log t)^{1/3-\varepsilon}}
\] 
by Lemma~$\ref{MLforSrho}$  and the estimate is complete.

If $2\xi t h'(t)<1$, then a treatment like the one for $M_3'$ suffices to conclude. More precisely, we have
\[
E_2\le \sum_{\rho\in M_2}|I_{t,\xi}(\rho)|\lesssim \sum_{0<-\gamma<2\xi th'(t)}|I_{t,\xi}(\rho)|\lesssim \sum_{|\gamma|<1}|I_{t,\xi}(\rho)|\lesssim \sum_{|\gamma|<1}t^{\beta}\lesssim_{\varepsilon}te^{-(\log t)^{1/3-\varepsilon}}\text{,}
\] 
since there are finitely many roots with $\gamma \in (-1,1)$ and the largest real part of the roots there is smaller than $1$. This concludes the estimate of $E_2$, and thus the proof of the lemma.
\end{proof}
Returning back to the estimate $\ref{weirdtruncroot}$, we see that Lemma~\ref{Keyroots} is applicable, it yields the desired estimate, and the proof of Lemma~$\ref{majarckey}$ is complete.
\end{proof}
\section{Applications to Waring-type problems}\label{waring}
In this section we prove Theorem~$\ref{nondiag}$. We begin by collecting some useful lemmas.
\begin{lemma}\label{FlL2norm}
Assume $c\in(1,2)$ and $h\in\mathcal{R}_c$ with $\varphi$ its  compositional inverse. Then there exists $C=C(h)$ such that for all $\lambda\in\mathbb{N}$ we have
\[
\bigg\|\sum_{n\le \lambda}\varphi'(n)e(n\xi)\bigg\|_{L^2_{d\xi}(\mathbb{T})}\le C\lambda^{-\frac{1}{2}}\varphi(\lambda) \text{.}
\]
\end{lemma}
\begin{proof}
Plancherel yields
\begin{equation}\label{L2normestnew1}
\bigg\|\sum_{n\le \lambda}\varphi'(n)e(n\xi)\bigg\|_{L^2_{d\xi}(\mathbb{T})}=\Big(\sum_{n\le \lambda}\varphi'(n)^2\Big)^{1/2}\text{.}
\end{equation}
We firstly prove that
\begin{equation}\label{telescopenew1}
\varphi'(m+1)^2\lesssim \frac{\varphi^2(m+1)}{m+1}-\frac{\varphi^2(m)}{m}\text{.}
\end{equation}
By the Mean Value Theorem there exists $\xi_m\in(m,m+1)$ such that $\frac{\varphi^2(m+1)}{m+1}-\frac{\varphi^2(m)}{m}=\frac{\varphi(\xi_m)^2}{\xi_m^2}(2\gamma-1+2\theta_1(\xi_m))$ since by Lemma~2.14 from \cite{MMR} we get
\[
(x^{-1}\varphi^2(x))'=x^{-1}2\varphi(x)\varphi'(x)-x^{-2}\varphi(x)^2=x^{-2}2\varphi(x)^2(\gamma+\theta_1(x))-x^{-2}\varphi(x)^2=\frac{\varphi(x)^2}{x^2}(2\gamma-1+2\theta_1(x))\text{.}
\]
For large $x$ we have that $|2\theta_1(x)|\le\frac{2\gamma-1}{2}$, since $\lim_{x\to\infty}\theta_1(x)=0$,  and thus for large $m$ we have
\[
\varphi'(m+1)^2\le \varphi'(\xi_m)^2\lesssim \frac{\varphi(\xi_m)^2}{\xi_m^2}\lesssim   \frac{\varphi(\xi_m)^2}{\xi_m^2}(2\gamma-1+2\theta(\xi_m))=\frac{\varphi^2(m+1)}{m+1}-\frac{\varphi^2(m)}{m}\text{.}
\] 
which establishes $\ref{telescopenew1}$. Finally, note that $\ref{L2normestnew1}$ becomes
\[
\bigg\|\sum_{n\le \lambda}\varphi'(n)e(n\xi)\bigg\|^2_{L^2_{d\xi}(\mathbb{T})}\le O(1)+\sum_{n\le\lambda}\Big(\frac{\varphi^2(n+1)}{n+1}-\frac{\varphi^2(n)}{n}\Big)\lesssim \frac{\varphi^2(\lambda+1)}{\lambda+1}\lesssim \frac{\varphi^2(\lambda)}{\lambda}\text{,}
\]
which implies the desired estimate.
\end{proof}
\begin{lemma}\label{geometricgrowth}
Assume $c\in(1,2)$ and $h\in\mathcal{R}_c$ with $\varphi$ its  compositional inverse. Then there exists $C=C(h)$ such that for all $\lambda\in\mathbb{N}$ we have
\[
\sum_{0\le i\le\log_2(\varphi(\lambda))}\frac{2^i}{\sqrt{h(2^i)}}\le C\lambda^{-1/2} \varphi(\lambda)\text{.}
\]
\end{lemma}
\begin{proof}
A variant of the argument below has already appeared in page 10. We claim that there exists $d=d(h)\in(0,1)$ and $x_0=x_0(h)$ such that $\frac{x}{\sqrt{h(x)}}\le d\frac{2x}{\sqrt{h(2x)}}$ for $x\ge x_0$. To see this note that by the definition of $h$ we get that
\[
\frac{\frac{x}{\sqrt{h(x)}}}{\frac{2x}{\sqrt{h(2x)}}}=2^{-1}\frac{\sqrt{h(2x)}}{\sqrt{h(x)}}=2^{-1+\frac{c}{2}}\exp\Big(\frac{1}{2}\int_{x}^{2x}\frac{\vartheta(t)}{t}dt\Big)\text{,}
\]
where $\vartheta(t)\to 0$ as $t\to\infty$. Thus there exists $x_0=x_0(h)$ such
that $|\vartheta(t)|\le 
\frac{2-c}{2}$ for all $t\ge x_0$. Thus for all $x\ge x_0$ we get
\[
\frac{\frac{x}{\sqrt{h(x)}}}{\frac{2x}{\sqrt{h(2x)}}}\le2^{-1+\frac{c}{2}}\exp\Big(\frac{2-c}{4}\int_{x}^{2x}\frac{1}{t}dt\Big)=2^{-1+\frac{c}{2}}e^{\frac{2-c}{4}\log 2}=2^{-1+\frac{c}{2}+\frac{2-c}{4}}=2^{-\frac{1}{2}+\frac{c}{4}}=2^{\frac{c-2}{4}}\eqqcolon d\in(0,1)\text{.}
\]
Let $L(x)\coloneqq \frac{x}{\sqrt{h(x)}}$ and note that for large $x\ge x_0$ we have $L(x)\le dL(2x)\iff L(x)\le \frac{d}{1-d}(L(2x)-L(x))$. Thus
\[
\sum_{0\le i\le \log_2(\varphi(\lambda))}\frac{2^i}{\sqrt{h(2^i)}}\lesssim O(1)+\sum_{\log_2 x_0 \le i\le\log_2(\varphi(\lambda))}(L(2^{i+1})-L(2^i))
\]
\[
\lesssim O(1)+L(2\varphi(\lambda))\lesssim O(1)+L(\varphi(\lambda))=O(1)+\lambda^{-1/2}\varphi(\lambda)\lesssim\lambda^{-1/2}\varphi(\lambda) \text{,}
\]
as desired.
\end{proof}
\begin{lemma}\label{difficultcounting}
Assume $c\in(1,2)$, $h\in\mathcal{R}_c$ with $\varphi$ its  compositional inverse, and assume $\theta\in(c-1,1)$. Then there exists $C=C(h,\theta)$ such that for all $\lambda\in\mathbb{N}$ we have
\[
\int_{|\xi|<\varphi(\lambda)^{-\theta}}\Big|\sum_{p\le\varphi(\lambda)}\log(p)e(\lfloor h(p)\rfloor)\Big|^2d\xi\le C\log^3(\lambda)\lambda^{-1}\varphi^2(\lambda)\text{.}
\]
\end{lemma}
\begin{proof}We focus on a dyadic piece. For convenience let 
\[
B_{\lambda}(\xi)=\sum_{p\le\varphi(\lambda)}\log(p)e(\lfloor h(p)\rfloor)\quad\text{and}\quad B_{\lambda,i}(\xi)=\sum_{2^i< p\le\min\{2^{i+1},\varphi(\lambda)\}}\log(p)e(\lfloor h(p)\rfloor)\text{,}
\]
for $i\in\{0,\dotsc,\lfloor\log_2(\varphi(\lambda))\rfloor\}$. Note that 
\begin{multline}\label{Returntothis}
\bigg(\int_{|\xi|<\varphi(\lambda)^{-\theta}}|B_\lambda(\xi)|^2d\xi\bigg)^{1/2}=\bigg(\int_{|\xi|<\varphi(\lambda)^{-\theta}}\Big|\sum_{i=0}^{\lfloor\log_2(\varphi(\lambda))\rfloor}B_{\lambda,i}(\xi)\Big|^2d\xi\bigg)^{1/2}\le
\\
\le\sum_{i=0}^{\lfloor\log_2(\varphi(\lambda))\rfloor}\bigg(\int_{|\xi|<\varphi(\lambda)^{-\theta}}|B_{\lambda,i}(\xi)|^2d\xi\bigg)^{1/2}
\text{.}
\end{multline}
We have that
\[
\int_{|\xi|<\varphi(\lambda)^{-\theta}}|B_{\lambda,i}(\xi)|^2d\xi=\int_{|\xi|<\varphi(\lambda)^{-\theta}}\sum_{2^i<p_1,p_2\le\min\{ 2^{i+1},\varphi(\lambda)\}}\log(p_1)\log(p_2)e\big((\lfloor h(p_1)\rfloor-\lfloor h(p_2)\rfloor)\xi\big)d\xi
\]
\[
\lesssim\sum_{2^i<p_1,p_2\le \min\{ 2^{i+1},\varphi(\lambda)\}}\log(p_1)\log(p_2)\min\bigg\{\varphi(\lambda)^{-\theta},\frac{1}{\big|\lfloor h(p_1)\rfloor-\lfloor h(p_2)\rfloor\big|}\bigg\}
\] 
\begin{equation}\label{2summandsfor1}
\lesssim \log^2(\lambda)\varphi(\lambda)^{-\theta}\sum_{\substack{2^i<p_1,p_2\le 2^{i+1}\\|\lfloor h(p_1)\rfloor-\lfloor h(p_2)\rfloor|<\varphi(\lambda)^{\theta}}}1+\log^2(\lambda)\sum_{\substack{2^i<p_1,p_2\le 2^{i+1}\\ | \lfloor h(p_1)\rfloor-\lfloor h(p_2)\rfloor|\ge\varphi(\lambda)^{\theta}}}\frac{1}{\big|\lfloor h(p_1)\rfloor-\lfloor h(p_2)\rfloor\big|}\text{.}
\end{equation}
For large $\lambda$ we estimate as follows
\[
|\{(p_1,p_2)
\in\big(\mathbb{P}\cap\big(2^i, 2^{i+1}\big]\big)^2:\,|\lfloor h(p_1)\rfloor-\lfloor h(p_2)\rfloor|<\varphi(\lambda)^{\theta}\}|
\]
\[
\lesssim|\{(n_1,n_2)
\in\big(2^i,2^{i+1}\big]^2:\, h(n_1)-\varphi(\lambda)^{\theta}-1< h(n_2)< h(n_1)+\varphi(\lambda)^{\theta}+1\}|
\] 
\begin{equation}\label{countingrough1}
\lesssim\sum_{2^i<n_1\le2^{i+1}}\sum_{2^i < n_2\le 2^{i+1}}1_{(h(n_1)-\varphi(\lambda)^{\theta}-1, h(n_1)+\varphi(\lambda)^{\theta}+1)}\big(h(n_2)\big)\text{.}
\end{equation}
For any interval of length $2\varphi(\lambda)^{\theta}+2$, since the gaps of $h(n_2)$ when $2^i< n_2\le 2^{i+1}$ are greater than $h(n_2+1)-h(n_2)\ge h'(2^i)$ by the Mean Value Theorem together with the fact that $h'$ is increasing, we will have that for any $n_1$, $h(n_2)$ may belong to this interval for at most $(2\varphi(\lambda)^{\theta}+2)\frac{1}{h'(2^i)}+1$ values of $n_2\in[2^i,2^{i+1}]$. Thus the expression $\ref{countingrough1}$ is bounded by
\[
\sum_{2^i<n_1\le2^{i+1}}\frac{2\varphi(\lambda)^{\theta}+2}{h'(2^i)}+1\lesssim 2^i +2^i\frac{\varphi(\lambda)^{\theta}}{h'(2^i)}\lesssim  2^i +2^{2i}\frac{\varphi(\lambda)^{\theta}}{h(2^i)}\text{,}
\]
and thus the first summand in $\ref{2summandsfor1}$ is bounded by
\begin{equation}\label{1stone1}
\log^2(\lambda)\varphi(\lambda)^{-\theta}\Big(2^i+\frac{2^{2i}\varphi(\lambda)^{\theta}}{h(2^i)}\Big)=\log^2(\lambda)\varphi(\lambda)^{-\theta}2^i+\log^2(\lambda)\frac{2^{2i}}{h(2^i)}\text{.}
\end{equation}
For the second summand of $\ref{2summandsfor1}$ we note that for large $\lambda$ we get
\[
\sum_{\substack{2^i<p_1,p_2\le 2^{i+1}\\ | \lfloor h(p_1)\rfloor-\lfloor h(p_2)\rfloor|\ge\varphi(\lambda)^{\theta}}}\frac{1}{\big|\lfloor h(p_1)\rfloor-\lfloor h(p_2)\rfloor\big|}\lesssim \sum_{\substack{2^i<n_1<n_2\le 2^{i+1}\\ h(n_2)-h(n_1)\ge\frac{1}{2}\varphi(\lambda)^{\theta}}}\frac{1}{h(n_2)-h(n_1)}
\]
\[
\lesssim\sum_{\substack{2^i<n_1<n_2\le 2^{i+1}\\ h(n_2)-h(n_1)\ge\frac{1}{2}\varphi(\lambda)^{\theta}}}\sum_{l=1}^{2\lambda}\frac{1}{h(n_2)-h(n_1)}1_{[\frac{1}{2}\varphi(\lambda)^{\theta}l,\frac{1}{2}\varphi(\lambda)^{\theta}(l+1)\big)}(h(n_2)-h(n_1))
\]
\begin{equation}\label{exprcount1}
\lesssim\sum_{2^i<n_1\le 2^{i+1}}\sum_{l=1}^{2\lambda}\frac{1}{l\varphi(\lambda)^{\theta}}\sum_{2^i<n_2\le 2^{i+1}}1_{[h(n_1)+\frac{1}{2}\varphi(\lambda)^{\theta}l,h(n_1)+\frac{1}{2}\varphi(\lambda)^{\theta}(l+1)\big)}(h(n_2))\text{,}
\end{equation}
and arguing similarly to before we get that every $l$ and $n_1$ 
\[
\sum_{2^i<n_2\le 2^{i+1}}1_{[h(n_1)+\frac{1}{2}\varphi(\lambda)^{\theta}l,h_1(n_1)+\frac{1}{2}\varphi(\lambda)^{\theta}(l+1)\big)}(h(n_2))\lesssim (\varphi(\lambda)^{\theta}+1)\frac{1}{h'(2^i)}+1\lesssim \frac{\varphi(\lambda)^{\theta}}{h'(2^i)}+1\text{.}
\]
Thus the expression $\ref{exprcount1}$ is bounded by
\[
\lesssim2^i \varphi(\lambda)^{-\theta}\log(\lambda)\Big(\frac{\varphi(\lambda)^{\theta}}{h'(2^i)}+1\Big)\text{,}
\]
which in turn implies that the second summand in  $\ref{2summandsfor1}$ will be bounded by 
\begin{equation}\label{2ndone}
\lesssim\log^3(\lambda)2^i\varphi(\lambda)^{-\theta}\Big(\frac{\varphi(\lambda)^{\theta}}{h'(2^i)}+1\Big)\lesssim \log^3(\lambda)2^{2i}\frac{1}{h(2^i)}+\log^3(\lambda)2^{i}\varphi(\lambda)^{-\theta}\text{.}
\end{equation}
Combining everything we get 
\[
\int_{|\xi|<\varphi(\lambda)^{-\theta}}|B_{\lambda,i}(\xi)|^2d\xi\lesssim \log^3(\lambda)2^{2i}\frac{1}{h(2^i)}+\log^3(\lambda)2^{i}\varphi(\lambda)^{-\theta} \text{,}
\]
and returning to $\ref{Returntothis}$ we obtain
\[
\bigg(\int_{|\xi|<\varphi(\lambda)^{-\theta}}|B_\lambda(\xi)|^2d\xi\bigg)^{1/2}\le\sum_{i=0}^{\lfloor\log_2(\varphi(\lambda))\rfloor}\bigg(\int_{|\xi|<\varphi(\lambda)^{-\theta}}|B_{\lambda,i}(\xi)|^2d\xi\bigg)^{1/2}
\]
\[
\le\sum_{i=0}^{\lfloor\log_2(\varphi(\lambda))\rfloor}\Big(\log^3(\lambda)2^{2i}\frac{1}{h(2^i)}+\log^3(\lambda)2^{i}\varphi(\lambda)^{-\theta}\Big)^{1/2}\lesssim \log^{3/2}(\lambda)\bigg(\sum_{i=0}^{\lfloor\log_2(\varphi(\lambda))\rfloor}\frac{2^{i}}{\sqrt{h(2^i)}}+2^{i/2}\varphi(\lambda)^{-\theta/2}\bigg) 
\]
\[
\lesssim\log^{3/2}(\lambda)\big(\lambda^{-1/2}\varphi(\lambda)+\varphi(\lambda)^{1/2}\varphi(\lambda)^{-\theta/2}\bigg)\lesssim\ \log^{3/2}(\lambda)\lambda^{-1/2}\varphi(\lambda)
\text{,}
\]
where for the second to last estimate we have used Lemma~$\ref{geometricgrowth}$ and for the last estimate we have used the fact that $\varphi(\lambda)^{1/2-\theta/2}\lesssim\lambda^{-1/2}\varphi(\lambda)$. To see this note that since $c<1+\theta$ we get that $h(t)\lesssim t^{1+\theta}$. Applying for $t=\varphi(\lambda)$ yields $\lambda\lesssim \varphi(\lambda)^{1+\theta}$ and thus $\varphi(\lambda)^{1-\theta}\lesssim \lambda^{-1}\varphi(\lambda)^2$, which yields the desired estimate. Thus we get
\[
\int_{|\xi|<\varphi(\lambda)^{-\theta}}\Big|\sum_{p\le\varphi(\lambda)}\log(p)e(\lfloor h(p)\rfloor)\Big|^2d\xi\lesssim \log^3(\lambda)\lambda^{-1}\varphi^2(\lambda)\text{,} 
\]
as desired.
\end{proof}
\begin{proof}[Proof of Theorem~$\ref{nondiag}$] Let us start by noting that according to Theorem~1.2 from \cite{LPE}, the latter assertion of the present theorem implies the former since 
\begin{equation}
r_{h_1,h_2,h_3}(\lambda)=\frac{\Gamma(\gamma_1)\Gamma(\gamma_2)\Gamma(\gamma_3)}{\Gamma(\gamma_1+\gamma_2+\gamma_3)}\lambda^{2}\varphi_1'(\lambda)\varphi_2'(\lambda)\varphi_3'(\lambda)+o(\lambda^{2}\varphi_1'(\lambda)\varphi_2'(\lambda)\varphi_3'(\lambda))\text{,}
\end{equation}
and it is easy to see that our assumption $4(1-\gamma_1)+\frac{45}{4}(1-\gamma_2)+\frac{45}{4}(1-\gamma_3)<1$ implies the corresponding ones in the statement of Theorem~1.2 from \cite{LPE}. It therefore suffices to establish the estimate $\ref{SecondEstimatekey}$.

For convenience for any $i\in[3]$ let 
\[
B^i_{\lambda}(\xi)=\sum_{p\le \varphi_i(\lambda)}\log(p)e(\lfloor h_i(p)\rfloor\xi)\text{,}\quad F^i_{\lambda}(\xi)=\sum_{n\le \lambda}\varphi_i'(n)e(n\xi)\quad\text{and}\quad
G^i_{\lambda}(\xi)=\sum_{m\le \varphi_i(\lambda)}e(\lfloor h_i(m)\rfloor)\text{.}
\]
Note that 
\[
R_{h_1,h_2,h_3}(\lambda)=\int_{\mathbb{T}}B^1_{\lambda}(\xi)B^2_{\lambda}(\xi)B^3_{\lambda}(\xi)e(-\lambda\xi)d\xi
\]
\[
=\int_{|\xi|<\varphi_1(\lambda)^{-\theta_1}}B^1_{\lambda}(\xi)B^2_{\lambda}(\xi)B^3_{\lambda}(\xi)d\xi+\int_{\varphi_1(\lambda)^{-\theta_1}\le |\xi|\le 1/2}B^1_{\lambda}(\xi)B^2_{\lambda}(\xi)B^3_{\lambda}(\xi)e(-\lambda\xi)d\xi\eqqcolon I_{\text{maj}}(\lambda)+I_{\text{min}}(\lambda)\text{,}
\]
where $\theta_1=\frac{6c_1}{5}-\frac{14}{15}\in(c_1-1,1)$.\\\,\\
\textbf{Estimates for }$I_{\text{min}}$\textbf{.} Note that 
\begin{equation}\label{Iminest}
|I_{\text{min}}(\lambda)|\lesssim \max_{\varphi_1(\lambda)^{-\theta_1}\le |\xi|\le 1/2}|B^1_{\lambda}(\xi)|\|B^2_{\lambda}\|_{L^2(\mathbb{T})}\|B^3_{\lambda}\|_{L^2(\mathbb{T})}\text{.}
\end{equation}
By Plancherel we get that for $i\in[3]$
\[
\|B_{\lambda}^i\|^2_{L^2(\mathbb{T})}=\sum_{p\le \varphi_i(\lambda)}\log^2p\le\log(\lambda)\sum_{p\le \varphi_i(\lambda)}\log p\lesssim \log(\lambda)\varphi_i(\lambda)\text{,}
\]
and by Lemma~$\ref{Minor}$ we get for every $\chi\in\big(0,\frac{8-6c_1}{45}\big)$ we have
\[
\max_{\varphi_1(\lambda)^{-\theta_1}\le |\xi|\le 1/2}|B^1_{\lambda}(\xi)|\lesssim_{\chi} \varphi_1^{1-\chi}(\lambda)\text{.}
\]
Hence
\[
|I_{\text{min}}(\lambda)|\lesssim\log(\lambda)\varphi_1(\lambda)^{1-\chi}\big(\varphi_2(\lambda)\varphi_3(\lambda)\big)^{1/2}=
\]
\[
\varphi_1(\lambda)\varphi_2(\lambda)\varphi_3(\lambda)\lambda^{-1}\big(\log(\lambda)\lambda\varphi_1(\lambda)^{-\chi}\varphi_2(\lambda)^{-1/2}\varphi_3(\lambda)^{-1/2}\big)\text{.}
\]
Now we note that 
\begin{equation}\label{assumptions}
1<\frac{\gamma_2}{2}+\frac{\gamma_3}{2}+\frac{\gamma_1(8-6c_1)}{45}
\end{equation}
since the inequality above is equivalent to $4(1-\gamma_1)+\frac{45}{4}(1-\gamma_2)+\frac{45}{4}(1-\gamma_3)<1$. Since $\ref{assumptions}$ holds, by continuity there exists $\chi\in\big(0,\frac{8-6c_1}{45}\big)$, $0<\gamma_1'<\gamma_1$, $0<\gamma_2'<\gamma_2$ and $0<\gamma_3'<\gamma_3$ such that
\[
1<\frac{\gamma_2'}{2}+\frac{\gamma_3'}{2}+\chi \gamma_1'\text{,} 
\]
but then we get 
\[
\lambda\varphi_1(\lambda)^{-\chi}\varphi_2(\lambda)^{-1/2}\varphi_3(\lambda)^{-1/2}\lesssim \lambda \lambda^{-\chi\gamma_1'}\lambda^{-\gamma_2'/2}\lambda^{-\gamma_3'/2}=\lambda^{1-\chi \gamma_1'-\frac{\gamma_2'}{2}-\frac{\gamma_3'}{2}}\text{,}
\]
but since $1-\chi \gamma_1'-\frac{\gamma_2'}{2}-\frac{\gamma_3'}{2}<0$ we have that
\[
|I_{\text{min}}(\lambda)|\lesssim \varphi_1(\lambda)\varphi_2(\lambda)\varphi_3(\lambda)\lambda^{-1}e^{-(\log \lambda)^{\frac{1}{3}}}\lesssim \lambda^{2}\varphi_1'(\lambda)\varphi_2'(\lambda)\varphi_3'(\lambda)e^{-(\log \lambda)^{\frac{1}{3}}}\text{,}
\]
and the estimates for $I_{\text{min}}$ are completed.\\\,\\
\textbf{Estimates for }$I_{\text{maj}}$\textbf{.} We have
\begin{align}\label{SPLITnew}
I_{\text{maj}}(\lambda)=\int_{|\xi|<\varphi_1(\lambda)^{-\theta_1}}\big(B^1_{\lambda}(\xi)B^2_{\lambda}(\xi)B^3_{\lambda}(\xi)-F^1_{\lambda}(\xi)F^2_{\lambda}(\xi)F^3_{\lambda}(\xi)\big)e(-\lambda\xi)d\xi
\\
+\int_{|\xi|<\varphi_1(\lambda)^{-\theta_1}}F^1_{\lambda}(\xi)F^2_{\lambda}(\xi)F^3_{\lambda}(\xi)e(-\lambda\xi)d\xi-\int_{\mathbb{T}}F^1_{\lambda}(\xi)F^2_{\lambda}(\xi)F^3_{\lambda}(\xi)(\xi)e(-\lambda\xi)d\xi
\\
+\int_{\mathbb{T}}\big(F^1_{\lambda}(\xi)F^2_{\lambda}(\xi)F^3_{\lambda}(\xi)-G^1_{\lambda}(\xi)G^2_{\lambda}(\xi)G^3_{\lambda}(\xi)\big)e(-\lambda\xi)d\xi
\\
+\int_{\mathbb{T}}G^1_{\lambda}(\xi)G^2_{\lambda}(\xi)G^3_{\lambda}(\xi)e(-\lambda\xi)d\xi
\text{.}
\end{align}
We note that
\[
\int_{\mathbb{T}}G^1_{\lambda}(\xi)G^2_{\lambda}(\xi)G^3_{\lambda}(\xi)e(-\lambda\xi)d\xi=r_{h_1,h_2,h_3}(\lambda)\text{,}
\]
and with an argument identical to the one given in the proof of Theorem~1.2 in \cite{LPE} we get that there exists $\delta>0$ such that
\[
\Big|\int_{\mathbb{T}}\big(F^1_{\lambda}(\xi)F^2_{\lambda}(\xi)F^3_{\lambda}(\xi)-G^1_{\lambda}(\xi)G^2_{\lambda}(\xi)G^3_{\lambda}(\xi)\big)e(-\lambda\xi)d\xi\Big|\lesssim \varphi_1(\lambda)\varphi_2(\lambda)\varphi_3(\lambda)\lambda^{-1-\delta}
\] 
\[
 \lesssim\lambda^{2}\varphi_1'(\lambda)\varphi_2'(\lambda)\varphi_3'(\lambda)e^{-(\log \lambda)^{\frac{1}{3}}}\text{.}
\]
Therefore to conclude our proof it suffices to establish the following two estimates
\begin{equation}\label{1summandtobound}
\Big|\int_{|\xi|<\varphi_1(\lambda)^{-\theta_1}}\big(B^1_{\lambda}(\xi)B^2_{\lambda}(\xi)B^3_{\lambda}(\xi)-F^1_{\lambda}(\xi)F^2_{\lambda}(\xi)F^3_{\lambda}(\xi)\big)e(-\lambda\xi)d\xi\Big|\lesssim_{\varepsilon}\lambda^{2}\varphi_1'(\lambda)\varphi_2'(\lambda)\varphi_3'(\lambda)e^{-(\log \lambda)^{\frac{1}{3}-\varepsilon}} 
\end{equation}
\begin{equation}\label{2summandtobound}
\Big|\int_{\varphi_1(\lambda)^{-\theta_1}\le|\xi|\le 1/2}F^1_{\lambda}(\xi)F^2_{\lambda}(\xi)F^3_{\lambda}(\xi)e(-\lambda\xi)d\xi\Big|\lesssim_{\varepsilon}
\lambda^{2}\varphi_1'(\lambda)\varphi_2'(\lambda)\varphi_3'(\lambda)e^{-(\log \lambda)^{\frac{1}{3}-\varepsilon}} \text{,}
\end{equation}
since according to the previous discussion we will then have that
\[
R_{h_1,h_2,h_3}(\lambda)=I_{\text{maj}}(\lambda)+O\Big(\lambda^{2}\varphi_1'(\lambda)\varphi_2'(\lambda)\varphi_3'(\lambda)e^{-(\log \lambda)^{\frac{1}{3}}}\Big)
\]
\[
=r_{h_1,h_2,h_3}(\lambda)+O_{\varepsilon}\Big(\lambda^{2}\varphi_1'(\lambda)\varphi_2'(\lambda)\varphi_3'(\lambda)e^{-(\log \lambda)^{\frac{1}{3}-\varepsilon}}\Big)\text{,}
\]
which is the desired estimate. The estimate $\ref{2summandtobound}$ can be obtained as follows 
\[
\Big|\int_{\varphi_1(\lambda)^{-\theta_1}\le|\xi|\le 1/2}F^1_{\lambda}(\xi)F^2_{\lambda}(\xi)F^3_{\lambda}(\xi)e(-\lambda\xi)d\xi\Big|\le\max_{\varphi_1(\lambda)^{-\theta_1}\le|\xi|\le1/2}|F^1_{\lambda}(\xi)|\|F^2_{\lambda}\|_{L^2(\mathbb{T})}\|F^3_{\lambda}\|_{L^2(\mathbb{T})}
\]
\[
\lesssim_{\varepsilon}\varphi_1(\lambda)e^{-(\log\varphi_1(\lambda))^{\frac{1}{3}-\frac{\varepsilon}{2}}}\varphi_2(\lambda)\varphi_3(\lambda)\lambda^{-1}\lesssim\lambda^2\varphi_1'(\lambda)\varphi_2'(\lambda)\varphi_3'(\lambda)e^{-(\log\lambda)^{\frac{1}{3}-\varepsilon}}\text{,}
\]
where we have used Lemma~$\ref{FlL2norm}$ together with Lemma~$\ref{Minor}$. For the estimate $\ref{1summandtobound}$, let 
\[
J_{\lambda}\coloneqq\big(-\varphi_1(\lambda)^{-\theta_1},\varphi_1(\lambda)^{-\theta_1}\big)\text{,}
\]
and note that
\[
\Big|\int_{J_{\lambda}}\big(B^1_{\lambda}(\xi)B^2_{\lambda}(\xi)B^3_{\lambda}(\xi)-F^1_{\lambda}(\xi)F^2_{\lambda}(\xi)F^3_{\lambda}(\xi)\big)e(-\lambda\xi)d\xi\Big|
\]
\begin{multline*}
\le\int_{J_{\lambda}}|B^1_{\lambda}(\xi)B^2_{\lambda}(\xi)B^3_{\lambda}(\xi)-F^1_{\lambda}(\xi)B^2_{\lambda}(\xi)B^3_{\lambda}(\xi)|d\xi
\\
+\int_{J_{\lambda}}|F^1_{\lambda}(\xi)B^2_{\lambda}(\xi)B^3_{\lambda}(\xi)-F^1_{\lambda}(\xi)F^2_{\lambda}(\xi)B^3_{\lambda}(\xi)|d\xi
\\
+\int_{J_{\lambda}}|F^1_{\lambda}(\xi)F^2_{\lambda}(\xi)B^3_{\lambda}(\xi)-F^1_{\lambda}(\xi)F^2_{\lambda}(\xi)F^3_{\lambda}(\xi)|d\xi
\end{multline*}
\begin{multline}\label{1tobound}
\le\|B^1_{\lambda}-F^1_{\lambda}\|_{L^{\infty}(J_{\lambda})}\|B^2_{\lambda}\|_{L^{2}(J_{\lambda})}\|B^3_{\lambda}\|_{L^{2}(J_{\lambda})}+\|B^2_{\lambda}-F^2_{\lambda}\|_{L^{\infty}(J_{\lambda})}\|F^1_{\lambda}\|_{L^{2}(\mathbb{T})}\|B^3_{\lambda}\|_{L^{2}(J_{\lambda})}
\\
+\|B^3_{\lambda}-F^3_{\lambda}\|_{L^{\infty}(J_{\lambda})}\|F^1_{\lambda}\|_{L^{2}(\mathbb{T})}\|F^2_{\lambda}\|_{L^{2}(\mathbb{T})}\text{.}
\end{multline}
By Lemma~$\ref{Major}$ and Lemma~$\ref{FlL2norm}$ we immediately get that
\[
\|B^1_{\lambda}-F^1_{\lambda}\|_{L^{\infty}(J_{\lambda})}\lesssim_{\varepsilon}\varphi_1(\lambda)e^{-(\log \varphi_1(\lambda))^{\frac{1}{3}-\frac{\varepsilon}{4}}}\lesssim \varphi_1(\lambda)e^{-(\log \lambda)^{\frac{1}{3}-\frac{\varepsilon}{2}}}
\]
and $\|F_{\lambda}^i\|_{L^2(\mathbb{T})}\lesssim \lambda^{-1/2}\varphi_i(\lambda)$, $i\in[3]$. For the remaining expressions in $\ref{1tobound}$ we proceed as follows. Observe that
\[
1-c_2<\frac{\gamma_1}{\gamma_2}\theta_1\text{,}
\]
since this inequality is equivalent to $\gamma_2+\frac{14\gamma_1}{15}<\frac{11}{5}$. We may therefore choose $\theta_2\in\big(1-c_2,\min\{\frac{\gamma_1}{\gamma_2}\theta_1,1\}\big)$. We can apply Lemma~$\ref{Major}$ and Lemma~$\ref{difficultcounting}$ with $\theta_2$ and obtain
 \[
 \|B^2_{\lambda}-F^2_{\lambda}\|_{L^{\infty}\big(-\varphi_2(\lambda)^{-\theta_2},\varphi_2(\lambda)^{-\theta_2}\big)}\lesssim_{\varepsilon}\varphi_2(\lambda)e^{-(\log \varphi_2(\lambda))^{\frac{1}{3}-\frac{\varepsilon}{4}}}\lesssim \varphi_2(\lambda)e^{-(\log \lambda)^{\frac{1}{3}-\frac{\varepsilon}{2}}}
\]
and
\[
\int_{|\xi|<\varphi_2(\lambda)^{-\theta_2}}|B^2_\lambda(\xi)|^2d\xi
\lesssim\log^3(\lambda)\lambda^{-1}\varphi_2(\lambda)^2\text{.}
\]
Our choice of $\theta_2$ implies that $\varphi_2(\lambda)^{-\theta_2}> \varphi_1(\lambda)^{-\theta_1}$ for large $\lambda$. To see this, note that since $\gamma_2\theta_2<\gamma_1\theta_1$, we get that there exists $\gamma'_2>\gamma_2$ and $\gamma_1'<\gamma_1$ such that $\gamma'_2\theta_2<\gamma'_1\theta_1$, and then there exists $\varepsilon>0$ sufficiently small such that $\varphi_2(\lambda)^{\theta_2}\lesssim 
 \lambda^{\gamma_2'\theta_2}\le \lambda^{\gamma_1'\theta_1}\lesssim_{\varepsilon}  \varphi_1(\lambda)^{\theta_1}\lambda^{-\varepsilon}$, which implies $\varphi_2(\lambda)^{-\theta_2}> \varphi_1(\lambda)^{-\theta_1}$ for large $\lambda$. We therefore get
\[
\|B^2_{\lambda}-F^2_{\lambda}\|_{L^{\infty}(J_{\lambda})}\le \|B^2_{\lambda}-F^2_{\lambda}\|_{L^{\infty}\big(-\varphi_2(\lambda)^{-\theta_2},\varphi_2(\lambda)^{-\theta_2}\big)}\lesssim_{\varepsilon}\varphi_2(\lambda)e^{-(\log \lambda)^{\frac{1}{3}-\frac{\varepsilon}{2}}}
\]
and
\[
\int_{|\xi|<\varphi_1(\lambda)^{-\theta_1}}|B^2_\lambda(\xi)|^2d\xi\lesssim \int_{|\xi|<\varphi_2(\lambda)^{-\theta_2}}|B^2_\lambda(\xi)|^2d\xi
\lesssim\log^3(\lambda)\lambda^{-1}\varphi_2(\lambda)^2\text{.}
\]
One may choose $\theta_3$ in an analogous manner, proceed with the same reasoning, and conclude
\[
\|B^3_{\lambda}-F^3_{\lambda}\|_{L^{\infty}(J_{\lambda})}\lesssim_{\varepsilon}\varphi_3(\lambda)e^{-(\log \lambda)^{\frac{1}{3}-\frac{\varepsilon}{2}}}\quad\text{
and}\quad
\int_{|\xi|<\varphi_1(\lambda)^{-\theta_1}}|B^3_\lambda(\xi)|^2d\xi 
\lesssim\log^3(\lambda)\lambda^{-1}\varphi_3(\lambda)^2\text{.}
\]
We are now ready to bound the expression $\ref{1tobound}$. We get that
\[
\Big|\int_{|\xi|<\varphi_1(\lambda)^{-\theta}}\big(B^1_{\lambda}(\xi)B^2_{\lambda}(\xi)B^3_{\lambda}(\xi)-F^1_{\lambda}(\xi)F^2_{\lambda}(\xi)F^3_{\lambda}(\xi)\big)e(-\lambda\xi)d\xi\Big|
\]
\[
\lesssim_{\varepsilon}\varphi_1(\lambda)\varphi_2(\lambda)\varphi_3(\lambda)\lambda^{-1}\log^3(\lambda)e^{-(\log 
\lambda)^{\frac{1}{3}-\frac{\varepsilon}{2}}}\lesssim \lambda^{2}\varphi'_1(\lambda)\varphi'_2(\lambda)\varphi'_3(\lambda)e^{-(\log 
\lambda)^{\frac{1}{3}-\varepsilon}}\text{,}
\]
and the proof is complete.
\end{proof}

\end{document}